\documentclass[a4paper]{amsart}
\usepackage{amsmath, amsthm, amsfonts, mathtools, tikz-cd, amssymb}
\usepackage{enumitem}
\usepackage{stmaryrd}
\usepackage[T1]{fontenc}

\usepackage[pagebackref,pdfborder={0 0 0},urlcolor=black]{hyperref}
\hypersetup{pdfauthor={Kevin Li, Clara L\"oh, Marco Moraschini, Roman Sauer, Matthias Uschold},%
  pdftitle={The algebraic cheap rebuilding property}}

\title
      {The algebraic cheap rebuilding property}
\author[K.~Li]{Kevin Li}
\address{Institut f\"ur Mathematik, Freie Universit\"at Berlin, 14195 Berlin, Germany}
\email{kevin.li@fu-berlin.de}

\author[C.~L\"oh]{Clara L\"oh}
\address{Fakult\"at f\"ur Mathematik, Universit\"at Regensburg, 93040 Regensburg, Germany}
\email{clara.loeh@ur.de}

\author[M.~Moraschini]{Marco Moraschini}
\address{Dipartimento di Matematica, Universit{\`a} di Bologna, 40126 Bologna, Italy}
\email{marco.moraschini2@unibo.it}

\author[R.~Sauer]{Roman Sauer}
\address{Fakult\"at f\"ur Mathematik, Karlsruher Institut f\"ur Technologie, 76131 Karlsruhe, Germany}
\email{roman.sauer@kit.edu}

\author[M.~Uschold]{Matthias Uschold}
\address{Fakult\"at f\"ur Mathematik, Universit\"at Regensburg, 93040 Regensburg, Germany}
\email{matthias.uschold@ur.de}

\makeatletter
\@namedef{subjclassname@2020}{%
  \textup{2020} Mathematics Subject Classification}
\makeatother

\newcounter{commentcounter}

\theoremstyle{definition}
\newtheorem{defn}{Definition}[section]
\newtheorem{ex}[defn]{Example}
\newtheorem{rem}[defn]{Remark}
\theoremstyle{plain}
\newtheorem{thm}[defn]{Theorem}
\newtheorem{lem}[defn]{Lemma}
\newtheorem{prop}[defn]{Proposition}
\newtheorem{cor}[defn]{Corollary}

\numberwithin{equation}{section}

\newcommand{\IN}{\ensuremath\mathbb{N}}
\newcommand{\IZ}{\ensuremath\mathbb{Z}}
\newcommand{\IQ}{\ensuremath\mathbb{Q}}
\newcommand{\IR}{\ensuremath\mathbb{R}}
\newcommand{\IC}{\ensuremath\mathbb{C}}
\newcommand{\IF}{\ensuremath\mathbb{F}}
\let\N\IN
\let\Z\IZ
\let\Q\IQ

\let\C\IC

\newcommand{\calN}{\ensuremath\mathcal{N}}
\newcommand{\calG}{\ensuremath\mathcal{G}}

\newcommand{\sfF}{\ensuremath\mathsf{F}}
\newcommand{\sfFP}{\ensuremath\mathsf{FP}}
\newcommand{\sfFL}{\ensuremath\mathsf{FL}}
\newcommand{\sfB}{\ensuremath\mathsf{B}}
\newcommand{\sfI}{\ensuremath\mathsf{I}}
\newcommand{\sfA}{\ensuremath\mathsf{A}}
\newcommand{\sfH}{\ensuremath\mathsf{H}}
\newcommand{\sfT}{\ensuremath\mathsf{T}}
\newcommand{\sfCR}{\ensuremath\mathsf{CR}}
\newcommand{\sfCWR}{\ensuremath\mathsf{CWR}}
\newcommand{\sfCD}{\ensuremath\mathsf{CD}}
\newcommand{\sfFG}{\ensuremath\mathsf{FG}}
\newcommand{\sfCM}{\ensuremath\mathsf{CM}}
\newcommand{\sfC}{\ensuremath\mathsf{C}}

\newcommand{\bfX}{\ensuremath\mathbf{X}}
\newcommand{\bfY}{\ensuremath\mathbf{Y}}
\newcommand{\bfC}{\ensuremath\mathbf{C}}

\newcommand{\enum}{\rm{(\roman*)}}
\newcommand{\spann}[1]{{\ensuremath \langle{#1}\rangle}}
\newcommand{\into}{\ensuremath\hookrightarrow}
\newcommand{\onto}{\ensuremath\twoheadrightarrow}
\def\qand{\quad\text{and}\quad}
\newcommand{\qm}{{\ensuremath\mathord{?}}}

\DeclareMathOperator{\ind}{ind}
\DeclareMathOperator{\id}{id}
\DeclareMathOperator{\res}{res}
\DeclareMathOperator{\tors}{tors}
\DeclareMathOperator{\rk}{rk}
\DeclareMathOperator{\pt}{pt}\DeclareMathOperator{\Cone}{Cone}
\DeclareMathOperator{\coker}{coker}
\DeclareMathOperator{\colim}{colim}
\DeclareMathOperator{\interior}{int}
\DeclareMathOperator{\Core}{Core}

\newcommand{\ABFG}{Ab\'ert--Bergeron--Fr\k{a}czyk--Gaboriau}

\makeatletter
\def\l@subsection{\@tocline{2}{0pt}{2.5pc}{5pc}{}}
\makeatother

\makeatletter
\newcommand*{\axiomenum}[1]{%
  \expandafter\@axiomenum\csname c@#1\endcsname%
}

\newcommand*{\@axiomenum}[1]{%
  \ifcase#1\or(B-deg)\or(B-susp)\or(B-cone)%
    \else\@ctrerr\fi%
}
\AddEnumerateCounter{\axiomenum}{\@axiomenum}{(B-cone)}
\makeatother

\begin{document}

\begin{abstract}
	We present an axiomatic approach to combination theorems for various homological properties of groups and, more generally, of chain complexes.
	Examples of such properties include algebraic finiteness properties, $\ell^2$-invisibility, $\ell^2$-acyclicity, lower bounds for Novikov--Shubin invariants, and vanishing of homology growth.

        As a key example, we introduce an algebraic version of \ABFG's cheap rebuilding property that implies vanishing of torsion homology growth and fits into our axiomatic framework for combination theorems. 
	In particular, we obtain that certain graphs of groups with amenable vertex groups and elementary amenable edge groups have vanishing torsion homology growth.
\end{abstract}

\maketitle


\section{Introduction}

\noindent
\textbf{Bootstrapping Theorems.}
Many homological properties of groups enjoy inheritance results, e.g., for graphs of groups (such as amalgamated products and HNN-extensions) and for group extensions.
For a sequence~$\sfB_*=(\sfB_n)_{n\in \IZ}$ of classes of groups, the following are desirable:
\begin{enumerate}[label=(\arabic*)]
	\item\label{item:amalg intro} Let~$\Gamma$ be the fundamental group of a finite graph of groups.
	If all vertex groups lie in~$\sfB_n$ and all edge groups lie in~$\sfB_{n-1}$, then $\Gamma\in \sfB_n$;
	\item\label{item:extensions intro} Let~$1\to N\to \Gamma\to Q\to 1$ be a group extension.
	If~$N\in \sfB_m$ for all~$m\le n$ and~$Q$ is of type~$\sfF_n$, then $\Gamma\in \sfB_n$.
\end{enumerate}
Recall that a group~$\Gamma$ is \emph{of type~$\sfF_n$} if there exists a model for the classifying space~$K(\Gamma,1)$ with finite $n$-skeleton.
Statements~\ref{item:amalg intro} and~\ref{item:extensions intro} are very useful in practice, e.g., for inductive approaches.
For instance, if the group of integers~$\IZ$ lies in~$\sfB_m$ for all~$m\le n$, then it follows from~\ref{item:amalg intro} and~\ref{item:extensions intro} that all infinite elementary amenable groups of type~$\sfFP_\infty$ lie in~$\sfB_n$ (Example~\ref{ex:elementary:ame}).
In fact, both \ref{item:amalg intro} and~\ref{item:extensions intro} are special cases of the more general inheritance result below (Theorem~\ref{thm:bootstrapping intro}).

We introduce the notion of a \emph{bootstrappable property~$\sfB_*$ of
  groups} (Definition~\ref{defn:pre-boot}).
The sequence~$\sfB_*=(\sfB_n)_{n\in \IZ}$ of classes of groups is associated to a family~$(\sfB^\Gamma_n)_{n\in \IZ,\Gamma}$, where~$\Gamma$ ranges over all groups and~$\sfB^\Gamma_n$ is a class of $\IZ\Gamma$-chain complexes that is compatible with suspensions, mapping cones, and induction:
The class~$\sfB_n$ consists of the groups~$\Gamma$ for which the trivial $\IZ\Gamma$-module~$\IZ$ admits a projective resolution lying in~$\sfB^\Gamma_n$.
The basic example of a bootstrappable property of groups is algebraic finiteness properties, where~$\sfB^\Gamma_n$ is the class of $\IZ\Gamma$-chain complexes with finitely generated chain modules in degrees~$\le n$.

\begin{thm}[Bootstrapping Theorem~\ref{thm:bootstrapping}]
\label{thm:bootstrapping intro}
	Let~$\sfB_*$ be a bootstrappable property of groups.
	Let~$\Gamma$ be a group and let~$n\in \IN$.
	Let~$\Omega$ be a $\Gamma$-CW-complex such that the following hold:
	\begin{enumerate}[label=\enum]
		\item $\Omega$ is $(n-1)$-acyclic (i.e., $H_j(\Omega;\IZ)\cong H_j(\pt;\IZ)$ for all~$j\le n-1$);
		\item The quotient $\Gamma\backslash \Omega^{(n)}$ of the $n$-skeleton~$\Omega^{(n)}$ is compact;
		\item For every cell~$\sigma$ of~$\Omega$ with~$\dim(\sigma)\le n$, the stabiliser~$\Gamma_\sigma$ lies in~$\sfB_{n-\dim(\sigma)}$.
	\end{enumerate}
	Then~$\Gamma\in \sfB_n$.
\end{thm}
Theorem~\ref{thm:bootstrapping intro} is called a Bootstrapping Theorem (or Combination Theorem) for the obvious reason:
The property of the group~$\Gamma$ is obtained by ``bootstrapping'' from the properties of the stabilisers~$(\Gamma_\sigma)_\sigma$.
In practice, the axioms for bootstrappable properties are usually 
straightforward to verify.
Examples of bootstrappable properties of groups include the algebraic finiteness properties~$\sfFP_n$, $\ell^2$-invisibility, $\ell^2$-acyclicity, and lower bounds for Novikov--Shubin invariants (Section~\ref{sec:examples}). 
The Bootstrapping Theorems for the first three of these properties are well-known~\cite{Brown87, Sauer-Thumann14, Jo07}.
Our axiomatic treatment provides a unified point of view, highlighting the essential axioms that lead to a Bootstrapping Theorem.
Moreover, while previous proofs involved spectral sequences, our proof is inductive on the level of chain complexes and enjoys a certain simplicity.

In case the bootstrappable property~$\sfB_*$ satisfies an additional restriction axiom, it follows that the intersection~$\bigcap_{m\le n}\sfB_m$ of classes of groups is closed under commensurability and under passing from commensurated subgroups to overgroups~(Corollary~\ref{cor:commensurability}).

\hfill 

\noindent
\textbf{Torsion homology growth.}
In the following, we focus on bootstrappable properties related to the vanishing of (torsion) homology growth.
Recall that a \emph{residual chain}~$\Lambda_*=(\Lambda_i)_{i\in \IN}$ in a group~$\Gamma$ is a nested sequence $\Gamma=\Lambda_0\ge \Lambda_1\ge \cdots$, where each~$\Lambda_i$ is a finite index normal subgroup of~$\Gamma$, such that $\bigcap_{i \in \N} \Lambda_i=\{1\}$.
A group is \emph{residually finite} if it admits a residual chain.

\begin{defn}
\label{defn:growth intro}
	Let~$\Gamma$ be a residually finite group and let~$\Lambda_*$ be a residual chain in~$\Gamma$.
	Let~$j\in \IN$ and let~$\IF$ be a field.
	Define
	\begin{align*}
		\widehat{b}_j(\Gamma,\Lambda_*;\IF)
		&\coloneqq \limsup_{i\to \infty} \frac{\dim_\IF H_j(\Lambda_i;\IF)}{[\Gamma:\Lambda_i]};
		\\
		\widehat{t}_j(\Gamma,\Lambda_*)
		&\coloneqq \limsup_{i\to \infty} \frac{\log\tors H_j(\Lambda_i;\IZ)}{[\Gamma:\Lambda_i]}.
	\end{align*}
	Here~$\tors$ denotes the cardinality of the torsion subgroup. 
	We set~$\log \infty\coloneqq \infty$.
\end{defn}

The motivation to study these invariants comes from L\"uck's Approximation Theorem which states that~$\widehat{b}_j(\Gamma,\Lambda_*;\IQ)$ agrees with the $\ell^2$-Betti number~$b^{(2)}_j(\Gamma)$, provided that~$\Gamma$ is of type~$\sfFP_{j+1}$.
Conjecturally, $\widehat{t}_j(\Gamma,\Lambda_*)$ is related to $\ell^2$-torsion (Section~\ref{sec:L2-torsion}).
We refer to L{\"u}ck's survey~\cite{Lueck16_survey} for more background.
We study the vanishing of these invariants and introduce the following notation.

For~$n\in \IN$, the class~$\sfH_n(\IF)$ (resp.~$\sfT_n$) consists of all residually finite groups~$\Gamma$ satisfying $\widehat{b}_j(\Gamma,\Lambda_*;\IF)=0$ (resp.\ $\widehat{t}_j(\Gamma,\Lambda_*)=0$) for all residual chains~$\Lambda_*$ in~$\Gamma$ and all~$j\le n$.
Then the class~$\sfH_0(\IF)$ is the class of all infinite residually finite groups.
We also define the classes~$\sfH_\infty(\IF)\coloneqq \bigcap_{n\in \IZ} \sfH_n(\IF)$ and~$\sfT_\infty\coloneqq \bigcap_{n\in \IZ}\sfT_n$.

The sequence~$\sfH_*(\IF)$ is a bootstrappable property of residually finite groups (Proposition~\ref{prop:bootstrappable H}), whereas~$\sfT_*$ is not (Remark~\ref{rem:T not bootstrappable}).
While the rank of homology is bounded by the rank of the chain module, the situation is more complicated for torsion in homology.
The main estimate is given by a result of Gabber:
Let~$(X,\partial)$ be a chain complex consisting of finitely generated based free $\IZ$-modules. Then
\begin{equation}\label{eqn:Gabber}
	\log \tors H_j(X)\le \rk(X_j)\cdot \log_+\|\partial_{j+1}\|,
\end{equation}
where~$\|\partial_{j+1}\|$ denotes the operator-norm of the map $\partial_{j+1}\colon X_{j+1}\to X_j$ with respect to the $\ell^1$-norms on~$X_{j+1}$ and~$X_j$, and~$\log_+\coloneqq \max\{\log,0\}$.
(Gabber's original estimate~\cite{soule1999perfect} is more refined and phrased for $\ell^2$-norms.)
Based on the estimate~\eqref{eqn:Gabber}, the goal is to find bootstrappable properties of residually finite groups that are (degreewise) contained in~$\sfT_*$.

\ABFG~\cite{ABFG21} introduced the sequence of classes of residually finite groups satisfying the \emph{(geometric) cheap rebuilding property}, which is contained in~$\sfT_{*-1}$ and in~$\sfH_*(\IF)$ for every field~$\IF$ and admits a Bootstrapping Theorem.
Their definition is geometric and requires the existence of models for classifying spaces satisfying delicate estimates on the number of cells and norms of boundary maps, homotopy equivalences, and homotopies.
The proof that the (geometric) cheap rebuilding property admits a Bootstrapping Theorem~\cite[Theorem~F]{ABFG21} is quite sophisticated and relies on a quantitative version of Geoghegan's Rebuilding Lemma~\cite{Geoghegan08}.
As the base case, the group of integers~$\IZ$ satisfies the (geometric) cheap rebuilding property.
Then the Bootstrapping Theorem yields a plethora of examples satisfying the (geometric) cheap rebuilding property in appropriate degrees.
For instance, applied to the rational Tits building, one obtains:

\begin{thm}[\ABFG]
\label{thm:SLd}
	Let~$d\in \IN$ with $d\ge 3$ and let~$\IF$ be a field.
	Then $\mathrm{SL}_d(\IZ) \in \sfT_{d-2}\cap \sfH_{d-2}(\IF)$.
\end{thm}

Other applications of the Bootstrapping Theorem for the (geometric) cheap rebuilding property~\cite[Theorem~F]{ABFG21} include mapping class groups, certain Artin groups~\cite{ABFG21}, outer automorphism groups of free products of~$\IZ/2$~\cite{GGH22}, mapping tori of polynomially growing automorphisms~\cite{AHK22,AGHK23}, and inner-amenable groups~\cite{Uschold22}.

The article~\cite{ABFG21} was the main inspiration for our work and we owe much to its ingenuity.
Our starting point was to provide a more algebraic proof of Theorem~\ref{thm:SLd}. We reprove it in Remark~\ref{rem:arithmetic groups}. Apart from that, the motivation for our 
axiomatic approach is two-fold: 

First, to extend the scope of the method to unify existing and obtain new vanishing results. For example, 
while residually finite infinite amenable groups of type~$\sfF_\infty$ lie in~$\sfT_\infty$~\cite{KKN17} and in~$\sfH_\infty(\IF)$, it is not known if they satisfy the (geometric) cheap rebuilding property~\cite[Question~10.21]{ABFG21}.
In particular, no general Bootstrapping Theorem with amenable stabilisers for vanishing torsion homology growth seems to be known.
We make progress in this direction, when the stabilisers of vertices are amenable and the stabilisers of positive dimensional cells are elementary amenable (Theorem~\ref{thm:bootstrapping CWR intro}).

Second, in future work we will combine the algebraic cheap rebuilding properties with the
dynamical approach to (torsion) homology growth~\cite{LLMSU25}. This
dynamical approach goes beyond the residually finite case and
provides estimates for (torsion) homology growth in terms of dynamical
systems. While our algebraic approach seems suitable for this, we would not know how to do this with the geometric approach.

\hfill

\noindent
\textbf{The algebraic cheap rebuilding property.}
We introduce the sequences~$\sfCR_*$ and~$\sfCWR_*$ of classes of residually finite groups satisfying the \emph{algebraic cheap rebuilding property} and the \emph{algebraic cheap weak rebuilding property}, respectively (Definition~\ref{defn:CR}).
The definitions are purely algebraic;
they require the existence of free resolutions and homotopy retracts of free chain complexes with quantitative control on ranks of modules and on norms of maps in terms of the index of subgroups along residual chains.
For all~$n\in \IZ$ and every field~$\IF$ the following inclusions hold (Lemma~\ref{lem:CWR implies T}):
\[
	\sfCR_n\ \subset\  \sfCWR_n\ \subset\ \sfT_{n-1}\cap \sfH_n(\IF).
\]
The classes~$\sfCR_0$ and~$\sfCWR_0$ coincide and consist of all infinite residually finite groups (Lemma~\ref{lem:deg0 infinite}).
We denote $\sfCR_\infty\coloneqq \bigcap_{n\in \IZ}\sfCR_n$ and~$\sfCWR_\infty\coloneqq \bigcap_{n\in \IZ}\sfCWR_n$.

We show that~$\sfCR_*$ is a bootstrappable property of residually finite groups (Proposition~\ref{prop:bootstrappable CR}).
The key observation is that the mapping cone of homotopy retracts is again a homotopy retract, and that all maps can be made explicit allowing for norm estimates.
Our algebraic approach has some similarity with the topological one by Okun--Schreve~\cite{Okun-Schreve21}.
Then we obtain the Bootstrapping Theorem~\ref{thm:bootstrapping intro} for~$\sfCR_*$ immediately from our axiomatic framework.
While the algebraic cheap rebuilding property is heavily inspired by the (geometric) cheap rebuilding property of \ABFG~\cite{ABFG21}, there is no obvious implication (Remark~\ref{rem:rebuilding} and Remark~\ref{rem:CR}).

The basic example of a group that lies in~$\sfCR_\infty$ is the integers~$\IZ$ (Example~\ref{ex:integers}).
Then an application of the Bootstrapping Theorem~\ref{thm:bootstrapping intro} for~$\sfCR_*$ recovers Theorem~\ref{thm:SLd} (Remark~\ref{rem:arithmetic groups}).
Moreover, we deduce that all residually finite infinite elementary amenable groups of type~$\sfFP_\infty$ lie in~$\sfCR_\infty$ (Example~\ref{ex:elementary:ame}).
While we do not know if residually finite infinite (not necessarily elementary) amenable groups of type~$\sfFP_\infty$ lie in~$\sfCR_\infty$, we show that they lie in~$\sfCWR_\infty$ (Theorem~\ref{thm:amenable CWR}).

The sequence~$\sfCWR_*$ does not seem to be a bootstrappable property of residually finite groups. 
Nevertheless, we obtain the following modified Bootstrapping Theorem for~$\sfCWR_*$ under stronger assumptions on stabilisers of cells of dimension~$\ge 1$.

\begin{thm}[Theorem~\ref{thm:bootstrapping CWR}]
\label{thm:bootstrapping CWR intro}
	Let~$\Gamma$ be a residually finite group and let~$n\in \IN$.
	Let~$\Omega$ be a $\Gamma$-CW-complex such that the following hold:
	\begin{enumerate}[label=\enum]
		\item $\Omega$ is $(n-1)$-acyclic;
		\item $\Gamma\backslash \Omega^{(n)}$ is compact;
		\item For every vertex~$v$ of~$\Omega$, the stabiliser~$\Gamma_v$ lies in~$\sfCWR_n$;
		\item For every cell~$\sigma$ of~$\Omega$ with~$\dim(\sigma)\in \{1,\ldots,n\}$, the stabiliser~$\Gamma_\sigma$ lies in~$\sfCR_{n-\dim(\sigma)}$.
	\end{enumerate}
	Then~$\Gamma\in \sfCWR_n$ (and hence~$\Gamma\in \sfT_{n-1}$).
\end{thm}

We point out again that for every~$n\in \IN$ the classes~$\sfCWR_n$ and~$\sfCR_n$ contain only infinite groups.
In particular, they do not contain the trivial group.
As a special case, Theorem~\ref{thm:bootstrapping CWR intro} applies to certain graphs of amenable groups: 
\begin{cor}[{Corollary~\ref{cor:amenable graph}}]
  Let~$\Gamma$ be the fundamental group of a finite graph of groups that is residually finite
  and let $n \in \N$.
	If all vertex groups are infinite amenable of type~$\sfFP_n$ and all edge groups are infinite elementary amenable of type~$\sfFP_\infty$, then~$\Gamma\in \sfT_{n-1}$.
\end{cor}

\tableofcontents


\section{Preliminaries on chain complexes}

Bootstrappable properties involve mapping cones of chain
maps and the key ingredient for proving the Bootstrapping Theorem
(Theorem~\ref{thm:bootstrapping intro}) is to understand projective
replacements in terms of mapping cones.  In this section, we will
provide the corresponding background from homological algebra.

We fix some notational conventions.
For the following generalities,
we work over an arbitrary ring that will be left implicit.
A chain complex~$X$ consists of chain modules~$X_j$ and differentials~$\partial_j\colon X_j\to X_{j-1}$ for all~$j\in \IZ$.
Let~$f,g\colon X\to Y$ be chain maps.
For a chain homotopy~$H\colon X_*\to Y_{*+1}$ between~$f$ and~$g$, we write $H\colon f\simeq g$ to indicate that $\partial^Y_{j+1}\circ H_j+H_{j-1}\circ \partial^X_j=f_j-g_j$ for all~$j\in \IZ$.

Let~$X$ be a chain complex.
We say that~$X$ is \emph{concentrated in degrees~$\le n$} if $X_j=0$ for all~$j>n$.
A \emph{weak equivalence} of chain complexes is a quasi-isomorphism (i.e., a chain map that induces an isomorphism on homology in all degrees).

The \emph{suspension} of~$X$ is the chain complex~$\Sigma X$ with chain modules $(\Sigma X)_j\coloneqq X_{j-1}$ and differentials $\partial^{\Sigma X}_j\coloneqq -\partial^X_{j-1}$.
The suspension $\Sigma f\colon \Sigma X\to \Sigma Y$ of a chain map~$f\colon X\to Y$ is given by $(\Sigma f)_j\coloneqq f_{j-1}$.
We denote by~$\Sigma^kX$ and~$\Sigma^k f$ the $k$-fold suspension of~$X$ and~$f$, respectively. 
We denote by~$\Sigma^{-1}$ the desuspension functor, which is inverse to~$\Sigma$.

\subsection{Mapping cones}
We review the mapping cone construction for chain complexes and its functoriality.

A \emph{homotopy commutative square} consists of chain maps
	\begin{equation}
	\label{eqn:square}
	\begin{tikzcd}
		X\ar[r, "f"]\ar[d, "a" left]
		& Y\ar[d, "b"] 
		\\
		Z\ar[r, "g" below]
		& W
		\arrow[from=2-2, to=1-1, phantom, "\text{\Large$\circlearrowleft$}_{H}" description]
	\end{tikzcd}
	\end{equation}
	together with a chain homotopy $H\colon g\circ a\simeq b\circ f$.
	If~$H=0$, then the square is (strictly) commutative. 
	
\begin{defn}[Mapping cone]
	Let~$f\colon X\to Y$ be a chain map.
	The \emph{mapping cone}~$\Cone(f)$ is the chain complex with chain modules $\Cone(f)_j\coloneqq X_{j-1}\oplus Y_j$ and differentials $\partial_j\colon \Cone(f)_j\to \Cone(f)_{j-1}$ given by
	\[
		\partial_j(x,y)\coloneqq \bigl(-\partial^X_{j-1}(x),\partial^Y_j(y)+f_{j-1}(x)\bigr).
	\]
	The mapping cone is functorial with respect to the chain map in the following sense:
	A homotopy commutative square as in diagram~\eqref{eqn:square} induces a chain map \[(a,b;H)\colon \Cone(f)\to \Cone(g)\] given by
	\[
		(a,b;H)_j(x,y)\coloneqq \bigl(a_{j-1}(x),b_j(y)-H_{j-1}(x)\bigr).
	\]
\end{defn}

For example, we have $\Sigma X\cong \Cone(X\to 0)$ and, more generally, $X\oplus Y\cong \Cone(0\colon \Sigma^{-1}X\to Y)$.
We include the proof of the following basic lemma~\cite{weibel} for completeness.

\begin{lem}
\label{lem:split ses}
	Let $f\colon X\to Y$ be a chain map.
	The following hold:
	\begin{enumerate}[label=\enum]
		\item\label{item:split ses} There exists a short exact sequence of chain complexes
		\begin{equation}
		\label{eqn:split ses}
			0\to Y\xrightarrow{\iota} \Cone(f)\xrightarrow{\pi} \Sigma X\to 0
		\end{equation}
		which splits degreewise;
		\item There exists a natural long exact sequence of homology groups
		\begin{equation}
		\label{eqn:les}
			\cdots\to H_j(X)\to H_j(Y)\to H_j\bigl(\Cone(f)\bigr)\to H_{j-1}(X)\to \cdots;
		\end{equation}
		\item\label{item:2outof3} Let a homotopy commutative square as in diagram~\eqref{eqn:square} be given.
		Consider the diagram of mapping cone sequences
		\[\begin{tikzcd}
			X\ar{r}{f}\ar{d}[swap]{a} & Y\ar{r}\ar{d}{b} & \Cone(f)\ar{d}{(a,b;H)}
			\\
			Z\ar{r}[swap]{g} & W\ar{r}\arrow[from=2-2, to=1-1, phantom, "\text{\Large$\circlearrowleft$}_{H}" description] & \Cone(g)\arrow[from=2-3, to=1-2, phantom, "\text{\Large$\circlearrowleft$}_{0}" description]
		\end{tikzcd}\]
		If the chain maps~$a$ and~$b$ are weak equivalences, then so is~$(a,b;H)$;
				\item\label{item:truncation} For every~$n\in \IZ$, there exists a subcomplex~$Y^n$ of~$Y$ that is concentrated in degrees~$\ge n$ such that $\tau_{<n}Y\coloneqq \Cone(Y^n\to Y)$ satisfies
		\[
				H_j(\tau_{<n}Y)\cong
				\begin{cases}
					0 & \text{for }j\ge n; \\
					H_j(Y) & \text{for }j\le n-1.				
				\end{cases}
		\]
		In particular, $(\tau_{<n}Y)_j = Y_j$ for all~$j\le n$.
	\end{enumerate}
	\begin{proof}
		(i) Define the chain maps~$\iota$ and~$\pi$ by $\iota(y)\coloneqq (0,y)$ and $\pi(x,y)\coloneqq x$. 
		Then a degreewise split $\sigma$ of~$\iota$ is given by $\sigma(x,y)=y$.
		
		(ii) Consider the long exact homology sequence associated to the short exact sequence~\eqref{eqn:split ses} and apply the suspension isomorphism $H_j(\Sigma X)\cong H_{j-1}(X)$.
		
		(iii) Consider the (strictly) commutative diagram with short exact rows
		\[\begin{tikzcd}
			0\ar{r} & Y\ar{r}\ar{d}{b} & \Cone(f)\ar{r}\ar{d}{(a,b;H)} & \Sigma X\ar{r}\ar{d}{\Sigma a} & 0
			\\
			0\ar{r} & W\ar{r} & \Cone(g)\ar{r} & \Sigma Z\ar{r} & 0
		\end{tikzcd}\]
		If the chain maps~$a$ (or equivalently,~$\Sigma a$) and~$b$ are weak equivalences, then so is~$(a,b;H)$ by the naturality of the long exact sequence~\eqref{eqn:les} and the Five Lemma.
				
		(iv) Take~$Y^n$ to be the subcomplex of~$Y$ given by
		\[
			Y^n_j\coloneqq \begin{cases}
				Y_j & \text{for }j\ge n+1; \\
				\ker(\partial^Y_n) & \text{for }j=n; \\
				0 & \text{for }j\le n-1.
			\end{cases}
		\]
		It follows from the long exact sequence~\eqref{eqn:les} that~$\tau_{<n}Y=\Cone(Y^n\to Y)$ is as desired.
	\end{proof}
\end{lem}

In Section~\ref{sec:CR}, we will use a ``higher'' functoriality of the mapping cone with respect to cubes that are coherently homotopy commutative.
	A \emph{homotopy commutative cube} consists of six homotopy commutative squares
	\begin{equation}
	\label{eqn:cube}
	\begin{tikzcd}[row sep=large, column sep=huge]
		X'\ar[rrr, "f'"]\ar[ddd, "a'" left] 
		& \arrow[from=1-2, to=2-3, phantom, "\text{\Large$\circlearrowleft$}_{F}" description]
		&& Y'\ar[ddd, "b'"]
		& \arrow[from=1-5, to=4-4, phantom, "\text{\Large$\circlearrowleft$}_{H'}" description]
		\\
		\arrow[from=2-1, to=3-2, phantom, "\text{\Large$\circlearrowleft$}_{A}" description]
		& X\ar[r, "f"]\ar[lu, "\xi\ " below]\ar[d, "a" left] 
		& Y\ar[ru, "\ \upsilon" below]\ar[d, "b"]
		\\
		& Z\ar[r, "g" below]\ar[ld, "\zeta\ " above] & W\ar[dr, "\ \omega" above]
		& \arrow[from=3-4, to=2-3, phantom, "\text{\Large$\circlearrowleft$}_{B}" description]
		\\
		Z'\ar[rrr, "g'" below] 
		& \arrow[from=4-2, to=3-3, phantom, "\text{\Large$\circlearrowleft$}_{G}" description]
		&& W'
		\arrow[from=2-2, to=3-3, phantom, "\text{\Large$\circlearrowleft$}_{H}" description]
	\end{tikzcd}
	\end{equation}
	where 
	\begin{align*}
		H\colon & g\circ a\simeq b\circ f; \\
		H'\colon & g'\circ a'\simeq b'\circ f'; \\
		A\colon & \zeta\circ a\simeq a'\circ \xi; \\
		B\colon & \omega\circ b\simeq b'\circ \upsilon; \\
		F\colon & f'\circ \xi\simeq \upsilon\circ f; \\
		G\colon & g'\circ \zeta\simeq \omega\circ g;
	\end{align*}
	together with a map~$\Phi\colon X_*\to W'_{*+2}$ such that
	\begin{equation}
	\label{eqn:filler}
	\begin{aligned}
		\partial^{W'}_{j+2} &\circ \Phi_j -\Phi_{j-1}\circ \partial^X_j=
		\\
		& \omega_{j+1}\circ H_j-H'_j\circ \xi_j + B_j\circ f_j-g'_{j+1}\circ A_j + G_j\circ a_j-b'_{j+1}\circ F_j. 
	\end{aligned}
	\end{equation}
The right hand side in equation~\eqref{eqn:filler} is the sum of all six maps~$X_*\to W'_{*+1}$ in diagram~\eqref{eqn:cube}, where opposite faces of the cube contribute with opposite signs.
In this sense, the map~$\Phi$ is a \emph{filler} of the cube.

\begin{lem}
\label{lem:cube}
	Let a homotopy commutative cube as in diagram~\eqref{eqn:cube} be given.
	Then the following square is homotopy commutative
	\begin{equation*}
	\begin{tikzcd}[column sep=huge]
		\Cone(f)\ar[r,"{(\xi,\upsilon;F)}"]\ar[d,"{(a,b;H)}" left]
		& \Cone(f')\ar[d, "{(a',b';H')}"]
		\\
		\Cone(g)\ar[r, "{(\zeta,\omega;G)}" below]
		& \Cone(g')
		\arrow[from=2-2, to=1-1, phantom, "\text{\Large$\circlearrowleft$}_{\Psi}" description]
	\end{tikzcd}
	\end{equation*}
	where $\Psi\colon \Cone(f)_*\to \Cone(g')_{*+1}$ is defined as
	\[
		\Psi_j(x,y)\coloneqq \bigl(-A_{j-1}(x),B_j(y)-\Phi_{j-1}(x)\bigr).
	\]
	\begin{proof}
		For~$(x,y)\in X_{j-1}\oplus Y_j=\Cone(f)_j$, we have
		\begin{align*}
			& \partial^{\Cone(g')}_{j+1}\circ \Psi_j (x,y) + \Psi_{j-1}\circ \partial^{\Cone(f)}_j (x,y)
			\\
			&= \partial^{\Cone(g')}_{j+1}\bigl(-A_{j-1}(x),B_j(y)-\Phi_{j-1}(x)\bigr) + \Psi_{j-1}\bigl(-\partial^X_{j-1}(x),\partial^Y_j(y)+f_{j-1}(x)\bigr)
			\\
			&= \bigl(\partial^{Z'}_{j}\circ A_{j-1}(x),\partial^{W'}_{j+1}\circ B_j(y)-\partial^{W'}_{j+1}\circ \Phi_{j-1}(x)-g'_j\circ A_{j-1}(x)\bigr)
			+ 
			\\
			&\qquad\qquad \bigl(A_{j-2}\circ \partial^X_{j-1}(x),B_{j-1}\circ \partial^Y_j(y)+B_{j-1}\circ f_{j-1}(x)+\Phi_{j-2}\circ \partial^X_{j-1}(x)\bigr)
			\\
			&= \bigl(\partial^{Z'}_{j}\circ A_{j-1}(x)+A_{j-2}\circ \partial^X_{j-1}(x),\partial^{W'}_{j+1}\circ B_j(y)+B_{j-1}\circ \partial^Y_j(y)\bigr)+
			\\
			&\qquad\qquad \bigl(0,-\partial^{W'}_{j+1}\circ \Phi_{j-1}(x)+\Phi_{j-2}\circ \partial^X_{j-1}(x)+B_{j-1}\circ f_{j-1}(x)-g'_j\circ A_{j-1}(x)\bigr)
			\\
			&= \bigl(\zeta_{j-1}\circ a_{j-1}(x)-a'_{j-1}\circ\xi_{j-1}(x), \omega_j\circ b_j(y)-b'_j\circ \upsilon_j(y)\bigr)+
			\\
			&\qquad\qquad \bigl(0,-\omega_j\circ H_{j-1}(x)+H'_{j-1}\circ \xi_{j-1}(x)-G_{j-1}\circ a_{j-1}(x)+b'_j\circ F_{j-1}(x)\bigr)
			\\
			&= \bigl(\zeta_{j-1}\circ a_{j-1}(x),\omega_j\circ b_j(y)-\omega_j\circ H_{j-1}(x)-G_{j-1}\circ a_{j-1}(x)\bigr) -
			\\
			&\qquad\qquad \bigl(a'_{j-1}\circ \xi_{j-1}(x),b'_j\circ \upsilon_j(y)-b'_j\circ F_{j-1}(x)-H'_{j-1}\circ \xi_{j-1}(x)\bigr)
			\\
			&= (\zeta,\omega;G)_j\bigl(a_{j-1}(x),b_j(y)-H_{j-1}(x)\bigr) - (a',b';H')_j\bigl(\xi_{j-1}(x),\upsilon_j(y)-F_{j-1}(x)\bigr)
			\\
			&= (\zeta,\omega;G)_j\circ (a,b;H)_j(x,y)-(a',b';H')_j\circ (\xi,\upsilon;F)_j (x,y).
		\end{align*}
		For the fourth equality, we used that~$A$ and~$B$ are chain homotopies and equation~\eqref{eqn:filler}.
		The other steps consist merely of expanding the definitions of the maps and rearranging the terms.
	\end{proof}
\end{lem}

\subsection{Projective replacements}
\label{sec:projective replacements}
Let~$X$ be a chain complex and let~$k\in \IZ$.
The \emph{$k$-skeleton} of~$X$ is the chain subcomplex~$X^{(k)}$ with chain modules
\[
	X^{(k)}_j\coloneqq \begin{cases}
		0 & \text{if }j\ge k+1; \\
		X_j & \text{if }j\le k; 
	\end{cases}
\]
and differentials
\[
	\partial^{X^{(k)}}_j\coloneqq \begin{cases}
		0 & \text{if }j\ge k+1; \\
		\partial^X_j & \text{if }j\le k. 
	\end{cases}
\]
If we regard the module~$X_k$ as a chain complex concentrated in degree zero, then~$X^{(k)}$ is the mapping cone of the chain map $\Sigma^{k-1}X_k\to X^{(k-1)}$ given by~$\partial^X_k$ in degree~$k-1$ and zero otherwise.

We recall a construction of projective replacements for chain complexes~\cite[Lemma~1.5]{Brown82_euler}. 
Given a chain complex~$X$ and projective resolutions of its chain modules, it produces a \emph{projective} chain complex~$\widehat{X}$ that is weakly equivalent to~$X$, obtained by iterated mapping cones from (suspensions of) the given projective resolutions.

\begin{prop}[Projective replacement]
\label{prop:blow-up}
	Let $X$ be a chain complex that is concentrated in degrees~$\ge 0$. 
	For each~$j\ge 0$, let $P^j = (P_i^j)_{i \geq 0}$ be a projective resolution of the module~$X_j$.
	Then there exists a projective chain complex~$\widehat{X}$ together with a filtration $(\widehat{X}^{[k]})_{k \geq 0}$ and a chain map $q\colon \widehat{X}\to X$ such that for every~$k\ge 0$ the following hold:
	\begin{enumerate}[label=\enum]
		\item $\widehat{X}^{[k]}$ is the mapping cone of a chain map~$\Sigma^{k-1}P^k\to \widehat{X}^{[k-1]}$;
		\item The restriction $q^k\coloneqq q|_{\widehat{X}^{[k]}}\colon \widehat{X}^{[k]}\to X^{(k)}$ is a weak equivalence.
	\end{enumerate}
	In particular, the chain modules of~$\widehat{X}$ are of the form 
	\[
		\widehat{X}_n=\bigoplus_{j+i=n} P^j_i
	\]
	and $q\colon \widehat{X}\to X$ is a weak equivalence.
	\medskip
	
	\noindent
	\textup{We call the chain map~$q\colon \widehat{X}\to X$ a \emph{projective replacement} of~$X$ with respect to the projective resolutions~$(P^j)_j$.}
\begin{proof}
	We proceed by induction on~$k \geq 0$.
	Set $\widehat{X}^{[0]}\coloneqq P^0$ and let~$q^0\colon P^0\to X^{(0)}$ be given by the augmentation $P_0^0\onto X_0$ in degree~$0$. 
	Since~$P^0$ is a resolution of~$X_0$, the chain map~$q^0$ is a weak equivalence.
	 For the inductive step, assume that the chain complex $\widehat{X}^{[k-1]}$ and a weak equivalence $q^{k-1}\colon \widehat{X}^{[k-1]}\to X^{(k-1)}$ have been constructed.
	
	Consider the following diagram
	\[\begin{tikzcd}
		\Sigma^{k-1}P^k\ar{d}[swap]{\Sigma^{k-1}\varepsilon^k} & \widehat{X}^{[k-1]}\ar{d}{q^{k-1}} \\
		\Sigma^{k-1}X_k\ar{r}{\partial^X_k} & X^{(k-1)}
	\end{tikzcd}\]
	where $\varepsilon^k\colon P^k\to X_k$ is given by the augmentation $P^k_0\onto~X_k$ in degree~$0$.
	Since the chain complex $\Sigma^{k-1}P^k$ consists of projective modules and the chain map~$q^{k-1}$ is a weak equivalence, there exists a chain map~$\widehat{\partial^X_k}\colon \Sigma^{k-1}P^k\to \widehat{X}^{[k-1]}$ making the square homotopy commutative~\cite[Lemma~1.1]{Brown82_euler}.
	We set $\widehat{X}^{[k]}\coloneqq \Cone(\widehat{\partial^X_k})$ and let $q^k\colon \widehat{X}^{[k]}\to X^{(k)} = \Cone(\partial_k^X)$ be the induced chain map on mapping cones.
	Thus we have a diagram of mapping cone sequences
	\[\begin{tikzcd}
		\Sigma^{k-1}P^k\ar{r}{\widehat{\partial^X_k}}\ar{d}[swap]{\Sigma^{k-1}\varepsilon^k} & \widehat{X}^{[k-1]}\ar{d}{q^{k-1}}\ar{r} & \widehat{X}^{[k]}\ar{d}{q^k} \\
		\Sigma^{k-1}X_k\ar{r}{\partial^X_k} & X^{(k-1)}\ar{r} & X^{(k)}
	\end{tikzcd}\]
	Since the chain maps~$\Sigma^{k-1}\varepsilon^k$ and~$q^{k-1}$ are weak equivalences, so is the chain map~$q^k$ by Lemma~\ref{lem:split ses}~\ref{item:2outof3}. 
	
	Taking the limit, we set $\widehat{X}\coloneqq \colim_k \widehat{X}^{[k]}$ and $q\coloneqq \colim_k q^{k}$, which are as desired.
	By construction, we have
	\[
		\widehat{X}^{[k]}_n=\bigoplus_{\substack{j+i=n,\\ 0 \le j\le k}} P^j_i
	\]
	and $\widehat{X}^{[k]}_n=\widehat{X}_n$ for~$k\ge n$, and thus $\widehat{X}_n = \bigoplus_{j+i=n} P^j_i$.
	Since homology commutes with directed colimits and all maps~$(q^k)_k$ are weak equivalences, the map~$q$ is a weak equivalence.
\end{proof}
\end{prop}

\section{Bootstrapping}
The goal of this section is to prove the Bootstrapping Theorem~\ref{thm:bootstrapping}.
First, in Section~\ref{sec:bootstrapping} we work over a fixed ring (e.g., a group ring).
Then, in Section~\ref{sec:equivariant bootstrapping} we pass to the global equivariant setting over 
varying group rings to obtain corresponding classes of groups. 

\subsection{Bootstrappable properties of chain complexes}
\label{sec:bootstrapping}

We propose a set of axioms on classes of chain complexes that ensure stability under projective replacements (Section~\ref{sec:projective replacements}).
Various examples of classes satisfying these axioms will be discussed in Section~\ref{sec:examples}.
The classical example to have in mind is, for given~$n \in \N$, the class of chain complexes that are finitely generated in all degrees~$\le n$ (Section~\ref{sec:FG}).

\begin{defn}[Bootstrappable property of chain complexes]
\label{defn:bootstrappable}
	Let $R$ be a ring.
	A \emph{bootstrappable property of $R$-chain complexes} is a sequence~$\sfB_*=(\sfB_n)_{n\in \IZ}$ of classes of $R$-chain complexes that are closed under isomorphisms and satisfy the following axioms:
	Let~$n\in \IN$ and let~$X$ be an $R$-chain complex.
	\begin{enumerate}[label=\rm{\axiomenum*}]
		\item\label{ax:pos} \emph{Degree.} 
		If~$X$ is concentrated in degrees~$\ge 0$, then~$X\in \sfB_m$ for all~$m<0$;
		\item\label{ax:susp} \emph{Suspension.}
		$X\in \sfB_n$ if and only if~$\Sigma X\in \sfB_{n+1}$;
		\item\label{ax:cone} \emph{Mapping cone.}
		Let~$f\colon X\to Y$ be an $R$-chain map.
		If~$X\in \sfB_{n-1}$ and~\mbox{$Y\in \sfB_n$}, then~$\Cone(f)\in \sfB_n$.
	\end{enumerate}
\end{defn}

We regard axiom~\ref{ax:cone} as the most important of the three, because in our examples axioms~\ref{ax:pos} and~\ref{ax:susp} will hold for trivial reasons.
We record some immediate consequences of the axioms:

\begin{lem}
\label{lem:axioms conseq}
	Let~$\sfB_*$ be a bootstrappable property of $R$-chain complexes.
	Let~$n\in \IZ$ and
	let~$X$ and~$Y$ be $R$-chain complexes. Then the following hold: 
	\begin{enumerate}[label=\enum]
		\item\label{ax:deg} 
		If~$X$ is concentrated in degrees~$\ge n$, then~$X\in \sfB_m$ for all~$m<n$;
		\item\label{ax:finite sums}
		 If~$X,Y\in \sfB_n$, then~$X\oplus Y\in \sfB_n$.
	\end{enumerate}
	\begin{proof}
		(i) If~$X$ is concentrated in degrees~$\ge n$, then $\Sigma^{-n}X$ is concentrated in degrees~$\ge 0$ and hence lies in~$\sfB_m$ for all~$m<0$ by axiom~\ref{ax:pos}.
		Then $X\cong \Sigma^n\Sigma^{-n}X$ lies in~$\sfB_m$ for all~$m< n$ by axiom~\ref{ax:susp}.
		
		(ii) Since~$X\in \sfB_n$, the desuspension~$\Sigma^{-1}X$ lies in~$\sfB_{n-1}$ by axiom~\ref{ax:susp}. 
		Since~$Y\in \sfB_n$, then the direct sum $X\oplus Y\cong \Cone(0\colon \Sigma^{-1}X\to Y)$ lies in~$\sfB_n$ by axiom~\ref{ax:cone}.
	\end{proof}
\end{lem}

The axioms in Definition~\ref{defn:bootstrappable} are conceived in order to be compatible with projective replacements (Proposition~\ref{prop:blow-up}) in the following sense:

\begin{prop}[Bootstrapping for chain complexes]
\label{prop:bootstrapping}
	Let~$\sfB_*$ be a bootstrappable property of $R$-chain complexes and let~$n\in \IN$.
	Let~$X$ be an $R$-chain complex that is concentrated in degrees~$\ge 0$.
	Suppose that for all~$j\le n$, the $R$-module~$X_j$ admits a projective resolution lying in~$\sfB_{n-j}$.
	Then there exists a projective $R$-chain complex~$\overline{X}\in \sfB_n$ that is weakly equivalent to~$X$.
	\begin{proof}
		\textbf{Case~1:}
		Suppose that~$X$ is concentrated in degrees~$\le n$.
		For~$j\le n$, let~$P^j$ be a projective resolution of~$X_j$ lying in~$\sfB_{n-j}$. 
		Proposition~\ref{prop:blow-up} yields a projective replacement $q\colon \widehat{X}\to X$ with respect to the resolutions~$(P^j)_j$.
		The $R$-chain complex~$\widehat{X}$ is equipped with a filtration
		\[
			0=\widehat{X}^{[-1]}\subset \widehat{X}^{[0]}\subset \cdots\subset \widehat{X}^{[n]}=\widehat{X}
		\]
		such that for all~$k\in \{1,\ldots,n\}$ there exists a mapping cone sequence
		\[
			\Sigma^{k-1}P^k\to \widehat{X}^{[k-1]}\to \widehat{X}^{[k]}.
		\]
		Since $\widehat{X}^{[0]} = P^0\in \sfB_n$, it follows by induction on~$k$ from axioms~\ref{ax:susp} and \ref{ax:cone} that~$\widehat{X}\in \sfB_n$.
		Since the $R$-chain map~$q$ is a weak equivalence, the $R$-chain complex~$\overline{X}\coloneqq \widehat{X}$ is as desired.
		
		\textbf{Case~2:}
		Let~$X$ be arbitrary (possibly of infinite length).
		For~$j\le n$, let~$P^j$ be a projective resolution of~$X_j$ lying in~$\sfB_{n-j}$.
		For~$j>n$, let~$P^j$ be an arbitrary projective resolution of~$X_j$.
		We write~$X_{>n}\coloneqq X/X^{(n)}$ for the quotient complex. 
		Then~$X$ is the mapping cone of the $R$-chain map $\Sigma^{-1}(X_{>n})\to X^{(n)}$ given by~$\partial^X_{n+1}$ in degree~$n$ and zero otherwise.
		Proposition~\ref{prop:blow-up} yields projective replacements $q^{(n)}\colon \widehat{X^{(n)}}\to X^{(n)}$ and $q_{>n}\colon \widehat{X_{>n}}\to X_{>n}$ with respect to the resolutions~$(P^j)_j$.
		Consider the diagram
		\[\begin{tikzcd}
			\Sigma^{-1}(\widehat{X_{>n}})\ar{d}[swap]{\Sigma^{-1}(q_{>n})} & \widehat{X^{(n)}}\ar{d}{q^{(n)}} & \\
			\Sigma^{-1}(X_{>n})\ar{r}{\partial^X_{n+1}} & X^{(n)}
		\end{tikzcd}\]
		Since the $R$-chain complex~$\Sigma^{-1}(\widehat{X_{>n}})$ is projective and the $R$-chain map~$q^{(n)}$ is a weak equivalence, there exists an $R$-chain map~$\widehat{\partial^X_{n+1}}\colon \Sigma^{-1}(\widehat{X_{>n}})\to \widehat{X^{(n)}}$ that lifts~$\partial^X_{n+1}$, up to homotopy~\cite[Lemma~1.1]{Brown82_euler}. 
		We define $\overline{X}\coloneqq \Cone(\widehat{\partial^X_{n+1}})$ and $q\colon \overline{X}\to X$ to be the induced map on mapping cones.
		We have constructed a diagram of mapping cone sequences
		\[\begin{tikzcd}
			\Sigma^{-1}(\widehat{X_{>n}})\ar{r}{\widehat{\partial^X_{n+1}}}\ar{d}[swap]{\Sigma^{-1}(q_{>n})} & \widehat{X^{(n)}}\ar{r}\ar{d}{q^{(n)}} & \overline{X}\ar{d}{q} \\
			\Sigma^{-1}(X_{>n})\ar{r}{\partial^X_{n+1}} & X^{(n)}\ar{r} & X
		\end{tikzcd}\]
		Since the $R$-chain complexes~$\Sigma^{-1}(\widehat{X_{>n}})$ and~$\widehat{X^{(n)}}$ are projective, so is~$\overline{X}$.
		Since the $R$-chain maps~$\Sigma^{-1}(q_{>n})$ and~$q^{(n)}$ are weak equivalences, so is~$q$ by Lemma~\ref{lem:split ses}~\ref{item:2outof3}.
		By Case~1, $\widehat{X^{(n)}}$ lies in~$\sfB_n$.
		Since~$\Sigma^{-1}(\widehat{X_{>n}})$ is concentrated in degrees~$\ge n$, it lies in~$\sfB_{n-1}$ by Lemma~\ref{lem:axioms conseq}~\ref{ax:deg}.
		Hence~$\overline{X}$ lies in~$\sfB_n$ by axiom~\ref{ax:cone}.
	\end{proof}
\end{prop}

\subsection{Bootstrappable properties of groups and the Bootstrapping Theorem}
\label{sec:equivariant bootstrapping}

In order to relate bootstrappable properties of chain complexes over group rings for different groups, we demand a compatibility with the induction functor.

An \emph{is-class of groups} is a class of groups that is closed under isomorphisms and taking subgroups.
We will use the question mark symbol~$\qm$ to denote an is-class of groups (e.g., the class of all groups or the class of residually finite groups).

\begin{defn}[Equivariantly bootstrappable property of chain complexes]
\label{defn:equivariant bootstrappable}
	Let~$R$ be a ring and let~$\qm$ be an is-class of groups.
	An \emph{equivariantly bootstrappable property of chain complexes over~$R$} is a family~$\sfB^\qm_*=(\sfB^\Gamma_n)_{\Gamma\in \qm,n\in \IZ}$, where~$\sfB^\Gamma_n$ is a class of $R\Gamma$-chain complexes, such that for every group~$\Gamma\in \qm$, the sequence~$\sfB^\Gamma_*$ is a bootstrappable property of $R\Gamma$-chain complexes (Definition~\ref{defn:bootstrappable}), and for all~$n\in \IZ$ the following holds:
	\begin{enumerate}[label=\rm{(B-ind)}]
	\item\label{ax:ind} \emph{Induction.} 
		Let~$\Gamma\in \qm$, let~$\Delta$ be a subgroup of~$\Gamma$, and let~$X$ be an $R\Delta$-chain complex.
		If~$X\in \sfB_n^\Delta$, then $\ind_\Delta^\Gamma X \coloneqq R\Gamma\otimes_{R\Delta} X\in \sfB_n^\Gamma$. 
	\end{enumerate}
\end{defn}

\begin{defn}[Bootstrappable property of groups]
\label{defn:pre-boot}
	Let~$R$ be a ring and let~$\qm$ be an is-class of groups.
	Let $\sfB^?_*=(\sfB^\Gamma_n)_{\Gamma\in \qm, n\in \IZ}$ be a family, where~$\sfB^\Gamma_n$ is a class of $R\Gamma$-chain complexes.
	For every~$n\in \IZ$, we define the associated class~$\sfB_n$ of groups~$\Gamma\in \qm$ for which the trivial $R\Gamma$-module~$R$ admits a projective resolution lying in~$\sfB^\Gamma_n$.
	We denote the intersection $\sfB_\infty\coloneqq \bigcap_{n\in \IZ}\sfB_n$ of subclasses of the class of all groups.
	
	A \emph{bootstrappable property of the is-class~$?$ of groups over~$R$} is a sequence~$\sfB_*=(\sfB_n)_{n\in \IZ}$ of classes of groups that is associated to some equivariantly bootstrappable property~$\sfB^?_*$ of chain complexes over~$R$ (Definition~\ref{defn:equivariant bootstrappable}).
\end{defn}

Our axiomatic framework culminates in the Bootstrapping Theorem for groups.

\begin{thm}
\label{thm:bootstrapping algebraic}
	Let~$\sfB_*$ be a bootstrappable property of an is-class~$\qm$ of groups over a ring~$R$.
	Let~$\Gamma\in \qm$ and let~$n\in \IN$.
	Let~$X$ be an $R\Gamma$-chain complex such that the following hold:
	\begin{enumerate}[label=\enum]
		\item $X$ is $(n-1)$-acyclic (i.e., $H_0(X)\cong R$ and $H_j(X)=0$ for all~$j\le n-1$);
		\item For all~$j\le n$, the $R\Gamma$-module~$X_j$ is isomorphic to $\bigoplus_{\sigma\in S_j} R[\Gamma/\Gamma_\sigma]$, where~$S_j$ is a finite set and~$\Gamma_\sigma$ is a subgroup of~$\Gamma$ lying in~$\sfB_{n-j}$. 
	\end{enumerate}
	Then~$\Gamma\in \sfB_n$.
	\begin{proof}
		Lemma~\ref{lem:split ses}~\ref{item:truncation} yields an $R\Gamma$-resolution~$\tau_{<n}X$ of~$R$ with $(\tau_{<n}X)^{(n)}=X^{(n)}$.
		Replacing~$X$ by~$\tau_{<n}X$, we may assume that~$X$ is an $R\Gamma$-resolution of~$R$.
		
		Let~$\sfB^\qm_*$ be an equivariantly bootstrappable property of chain complexes over~$R$ such that the bootstrappable property~$\sfB_*$ of the is-class~$\qm$ is associated to~$\sfB^\qm_*$.
		For~$j\le n$ and~$\sigma\in S_j$, since the group~$\Gamma_\sigma$ lies in~$\sfB_{n-j}$, there exists a projective $R[\Gamma_\sigma]$-resolution~$P^\sigma$ of~$R$ lying in~$\sfB_{n-j}^{\Gamma_\sigma}$.
		Then the projective $R\Gamma$-resolution $P^j\coloneqq \bigoplus_{\sigma\in S_j}\ind_{\Gamma_\sigma}^\Gamma P^\sigma$ of~$X_j$ lies in~$\sfB_{n-j}^\Gamma$ by axiom~\ref{ax:ind} and Lemma~\ref{lem:axioms conseq}~\ref{ax:finite sums}.	
		Finally, Proposition~\ref{prop:bootstrapping} yields a projective $R\Gamma$-chain complex~$\overline{X}\in \sfB_n^\Gamma$ that is weakly equivalent to~$X$.
		Hence~$\overline{X}$ is a projective $R\Gamma$-resolution of~$R$, witnessing that the group~$\Gamma$ lies in~$\sfB_n$.
	\end{proof}
\end{thm}

Suitable $R\Gamma$-chain complexes arise as the cellular chain complex of $\Gamma$-CW-com\-plex\-es.
Recall that a \emph{$\Gamma$-CW-complex} is a $\Gamma$-space that is built inductively by attaching $\Gamma$-orbits of cells.
In particular, the stabiliser of a cell fixes that cell pointwise.

\begin{thm}[Bootstrapping for groups]
\label{thm:bootstrapping}
	Let~$\sfB_*$ be a bootstrappable property of an is-class~$\qm$ of groups over a ring~$R$. 
	Let~$\Gamma\in \qm$ and let~$n\in \IN$.
	Let~$\Omega$ be a $\Gamma$-CW-complex such that the following hold:
	\begin{enumerate}[label=\enum]
		\item $\Omega$ is $(n-1)$-acyclic over~$R$ (i.e., $H_j(\Omega;R)\cong H_j(\pt;R)$ for all~$j\le n-1$);
		\item $\Gamma\backslash \Omega^{(n)}$ is compact;
		\item For every cell~$\sigma$ of~$\Omega$ with~$\dim(\sigma)\le n$, the stabiliser~$\Gamma_\sigma$ lies in~$\sfB_{n-\dim(\sigma)}$.
	\end{enumerate}
	Then $\Gamma\in \sfB_n$.
	\begin{proof}
		Let~$X$ be the cellular $R\Gamma$-chain complex of~$\Omega$.
		For~$j\le n$, the $R\Gamma$-module~$X_j$ is of the form
		\[
			X_j\cong \bigoplus_{\sigma\in S_j} R[\Gamma/\Gamma_\sigma],
		\]
		where~$S_j$ is a finite set of representatives for the $\Gamma$-orbits of $j$-cells of~$\Omega$.
		Then Theorem~\ref{thm:bootstrapping algebraic} applies to the $R\Gamma$-chain complex~$X$ and yields the claim.
		\end{proof}
\end{thm}

\begin{rem}
	The construction in the proof of Theorem~\ref{thm:bootstrapping algebraic} is an algebraic version of a topological ``blow-up'' construction by L\"uck~\cite{Lueck00,MPSS20} and Geoghegan~\cite{Geoghegan08}.
	Let~$\Omega$ be a $\Gamma$-CW-complex and let classifying spaces~$(E\Gamma_{\sigma})_\sigma$ for its cell stabilisers be given.
	Then there exists a \emph{free} $\Gamma$-CW-complex~$\widehat{\Omega}$ that is (non-equivariantly) homotopy equivalent to~$\Omega$ and, roughly speaking, obtained from~$\Omega$ by replacing each $\Gamma$-orbit of cells~$\Gamma/\Gamma_\sigma\times \sigma$ with the free $\Gamma$-CW-complex $\ind_{\Gamma_\sigma}^\Gamma E\Gamma_\sigma\times \sigma$.
\end{rem}

As special cases of Theorem~\ref{thm:bootstrapping algebraic} and Theorem~\ref{thm:bootstrapping}, we obtain inheritance results for classical group-theoretic constructions.
Recall that a group~$\Gamma$ is \emph{of type~$\sfFP_n(R)$} if the trivial $R\Gamma$-module~$R$ admits a projective resolution~$P$ such that the $R\Gamma$-module~$P_j$ is finitely generated for all~$j\le n$.
A group~$\Gamma$ is \emph{of type~$\sfFP_\infty(R)$} if~$\Gamma$ is of type~$\sfFP_n(R)$ for all~$n\in \IZ$.
For~$R=\IZ$, we write~$\sfFP_n$ and~$\sfFP_\infty$ instead of~$\sfFP_n(\IZ)$ and~$\sfFP_\infty(\IZ)$, respectively.

If~$\Gamma$ is of type~$\sfFP_n(R)$ and~$N$ is a normal subgroup of~$\Gamma$ of type~$\sfFP_{n-1}(R)$, then the quotient~$\Gamma/N$ is of type~$\sfFP_n(R)$.
This is well-known and implicit in Bieri's book~\cite[Proposition~2.7]{Bieri81}.
(The analogous statement for the geometric finiteness properties~$\sfF_n$ is proved in Geoghegan's book~\cite[Theorem~7.2.21]{Geoghegan08}.)

\begin{cor}
\label{cor:bootstrapping}
	Let~$\sfB_*$ be a bootstrappable property of an is-class~$\qm$ of groups over a ring~$R$.
	For all~$n\in \IZ$, the following hold:
	\begin{enumerate}[label=\enum]
		\item\label{item:boot gog} (Graphs of groups).
		Let~$\Gamma\in \qm$ be the fundamental group of a finite graph of groups.
		If all vertex groups lie in~$\sfB_n$ and all edge groups lie in~$\sfB_{n-1}$, then $\Gamma\in \sfB_n$; 
		\item\label{item:boot extension} (Extensions). 
		Let~$\Gamma\in \qm$ be a group containing a normal subgroup~$N$ such that the quotient~$\Gamma/N$ is of type~$\sfFP_n(R)$.
		If~$N\in \sfB_m$ for all~$m\le n$, then~$\Gamma\in \sfB_n$;
		\item\label{item:boot finite extension} (Finite extensions).
		Let~$\Gamma\in \qm$ and
		let~$N$ be a finite index normal subgroup of~$\Gamma$.
		If~$N\in \sfB_m$ for all~$m\le n$, then~$\Gamma\in \sfB_n$;
		\item\label{item:boot abelian} (Finitely generated free abelian groups).
		Let~$\Gamma\in \qm$ be a non-trivial finitely generated free abelian group.
		If~$\IZ\in \sfB_m$ for all~$m\le n$, then~$\Gamma\in \sfB_n$;
		\item\label{item:boot abelian normal} (Groups containing a free abelian normal subgroup).
		Let~$\Gamma\in \qm$ be a group of type~$\sfFP_n(R)$ containing a non-trivial finitely generated free abelian normal subgroup.
		If~$\IZ\in \sfB_m$ for all~$m\le n$, then~$\Gamma\in \sfB_n$.
	\end{enumerate}
	\begin{proof}
		(i) Apply Theorem~\ref{thm:bootstrapping} to the Bass--Serre tree.
		
		(ii) Since the quotient~$Q\coloneqq \Gamma/N$ is of type~$\sfFP_n(R)$, there exists a free $RQ$-resolution~$P$ of~$R$ such that the $RQ$-module~$P_j$ is finitely generated for all~$j\le n$.
		We consider the functor $\res_{\Gamma\onto Q}$ that turns an $RQ$-chain complex into an $R\Gamma$-chain complex using the group homomorphism $\Gamma\onto Q$ to define the $R\Gamma$-module structures.
		Then apply Theorem~\ref{thm:bootstrapping algebraic} to the $R\Gamma$-resolution~$\res_{\Gamma\onto Q}P$ of~$R$.
		
		(iii) This follows from part~\ref{item:boot extension}, since finite groups are of type~$\sfF_\infty$ and in particular of type~$\sfFP_\infty(R)$.
		
		(iv) This follows immediately from part~\ref{item:boot extension}.
		
		(v) This follows from parts~\ref{item:boot abelian} and~\ref{item:boot extension}.
	\end{proof}
\end{cor}

Part~\ref{item:boot gog} of Corollary~\ref{cor:bootstrapping} applies in particular to amalgamated products and HNN-extensions.
In parts~\ref{item:boot extension}--\ref{item:boot abelian normal}, since~$N$ (resp.~$\IZ$) lies in~$\sfB_m$ for all~$m\le n$, tautologically it also lies in~$\sfB_m$ for all~$m\le n-1$. Hence the conclusions are in fact that~$\Gamma$ lies in~$\sfB_m$ for all~$m\le n$.
In many of our examples of bootstrappable properties~$\sfB_*$ of groups (Section~\ref{sec:examples}), the class~$\sfB_n$ is contained in~$\sfB_m$ for all~$m\le n$ and the group~$\IZ$ lies in~$\sfB_\infty$.
However, we warn the reader that in many examples the trivial group does not lie in~$\sfB_0$.

Part~\ref{item:boot abelian normal} applies to non-trivial finitely generated torsion-free nilpotent groups because these groups are of type~$\sfF$ and have non-trivial centre.
More generally, properties of the group~$\IZ$ can be bootstrapped to infinite polycyclic-by-finite groups (i.e., groups that admit a subnormal series whose factors are cyclic or finite; Example~\ref{ex:polycyclic}).

In each of the following Examples~\ref{ex:polycyclic}, \ref{ex:elementary:ame}, \ref{ex:RAAG}, and~\ref{ex:arithmetic groups}, let~$\sfB_*$ be a bootstrappable property of an is-class~$\qm$ of groups over a ring~$R$.
Here it is crucial that the class~$\qm$ of groups is an is-class.

\begin{ex}[Polycyclic-by-finite groups]\label{ex:polycyclic}
	Let~$\Gamma\in \qm$ be an infinite polycyclic-by-finite group.
	If~$\IZ\in \sfB_m$ for all~$m\le n$, then~$\Gamma\in \sfB_n$. 

	Indeed, $\Gamma$ contains a finite index normal subgroup~$\Lambda$ that is poly-infinite-cyclic (i.e., $\Lambda$ admits a subnormal series whose factors are infinite cyclic)~\cite[Lem\-ma~5.11]{ceccherini2021topics}. 
	By Corollary~\ref{cor:bootstrapping}~\ref{item:boot finite extension}, it suffices to show that~$\Lambda$ lies in~$\sfB_m$ for all~$m \leq n$. 
	Let 	
	\[	
		\{1\} = \Lambda_0 \lhd \Lambda_1 \lhd \cdots \lhd \Lambda_k = \Lambda	
	\]	
	be a subnormal series witnessing the structure of~$\Lambda$ as a poly-infinite-cyclic group. 
	Since~$\Lambda_1$ is isomorphic to~$\IZ$, the group~$\Lambda_1$ lies in~$\sfB_m$ for all~$m \leq n$ by assumption.
	By induction on the length of the subnormal series, we may assume that~$\Lambda_{k-1}$ lies in~$\sfB_m$ for all~$m \leq n$. 
	Since $\Lambda=\Lambda_k$ is an extension of $\Lambda_{k-1}$ by $\IZ$, Corollary~\ref{cor:bootstrapping}~\ref{item:boot extension} shows that~$\Lambda$ lies in~$\sfB_m$ for all~$m \leq n$.
\end{ex}

In fact, properties of the group~$\IZ$ can be bootstrapped to infinite elementary amenable groups of type $\sfFP_\infty$.
Recall that the class of \emph{elementary amenable} groups is the smallest class of groups that contains all finite groups and all abelian groups, and that is closed under taking subgroups, quotients, extensions, and directed unions.

\begin{ex}[Elementary amenable groups]\label{ex:elementary:ame}
	Let~$\Gamma\in \qm$ be an infinite elementary amenable group of type~$\sfFP_\infty$.
	If~$\IZ\in \sfB_m$ for all~$m\le n$, then~$\Gamma\in \sfB_n$.
	
	Indeed, a characterisation by Kropholler--Mart\'{i}nez-P\'{e}rez--Nucinkis~\cite[Theorem~1.1]{KMPN} yields that either
	\begin{enumerate}
		\item $\Gamma$ is an infinite polycyclic-by-finite group; or
		\item $\Gamma$ contains a normal subgroup~$N$ such that:
		\begin{itemize} 
			\item $N$ is a strictly ascending HNN-extension~$H\ast_{H,t}$ over a finitely generated virtually nilpotent group~$H$, and
			\item $\Gamma/N$ is a Euclidean crystallographic group.
		\end{itemize}
	\end{enumerate}
	In case~(1), $\Gamma$ lies in~$\sfB_n$ by Example~\ref{ex:polycyclic}.
	In case~(2), since Euclidean crystallographic groups are of type~$\sfF_\infty$, it suffices to show that~$N$ lies in~$\sfB_m$ for all~$m\le n$ by Corollary~\ref{cor:bootstrapping}~\ref{item:boot extension}.
	Since the HNN-extension~$N=H\ast_{H,t}$ is strictly ascending, the group~$H$ must be infinite.
	The infinite finitely generated virtually nilpotent group~$H$ is infinite polycyclic-by-finite and hence lies in~$\sfB_m$ for all~$m\le n$ by Example~\ref{ex:polycyclic}.
	Then the HNN-extension~$N$ lies in~$\sfB_m$ for all~$m\le n$ by Corollary~\ref{cor:bootstrapping}~\ref{item:boot gog}.
\end{ex}

Properties of the group~$\IZ$ can also be bootstrapped to certain right-angled Artin groups that are iterated amalgamated products of free abelian groups.
Recall that the \emph{right-angled Artin group}~$A_\calG$ associated to a finite simplicial graph~$\calG$ has one generator for each vertex of~$\calG$ and the only relations are that two generators commute if their corresponding vertices are connected by an edge in~$\calG$.

\begin{ex}[Chordal right-angled Artin groups] 
\label{ex:RAAG}
		Let~$\calG$ be a non-empty connected finite simplicial graph that is chordal (i.e.,~$\calG$ contains no cycles of length~$\ge 4$ as full subgraphs).
		Let~$A_\calG$ be the right-angled Artin group associated to~$\calG$ and suppose that~$A_\calG\in \qm$.
		If~$\IZ\in \sfB_m$ for all~$m\le n$, then~$A_\calG\in \sfB_n$.
		 
		 Indeed, it is a classical result of Dirac~\cite{Dirac61} about chordal graphs that either 
		 \begin{enumerate}
		 	\item $\calG$ is complete; or 
			\item there exist full proper subgraphs~$\calG_1,\calG_2$ of~$\calG$ such that~$\calG=\calG_1\cup \calG_2$ and $\calG_0\coloneqq \calG_1\cap \calG_2$ is complete.
		\end{enumerate}	 
		In case~(1), the group~$A_\calG$ is non-trivial finitely generated free abelian and hence lies in~$\sfB_n$ by Corollary~\ref{cor:bootstrapping}~\ref{item:boot abelian}.
		In case~(2), the group~$A_\calG$ splits as an amalgamated product $A_\calG\cong A_{\calG_1}\ast_{A_{\calG_0}} A_{\calG_2}$.
		Since~$\calG$ is connected, the graph~$\calG_0$ is non-empty.
		 Then the group~$A_{\calG_0}$ is non-trivial finitely generated free abelian and hence lies in~$\sfB_{n-1}$ by Corollary~\ref{cor:bootstrapping}~\ref{item:boot abelian}.
		By induction on the number of vertices of~$\calG$, we may assume that~$A_{\calG_1}$ and~$A_{\calG_2}$ lie in~$\sfB_n$. 
		Hence~$A_\calG$ lies in~$\sfB_n$ by Corollary~\ref{cor:bootstrapping}~\ref{item:boot gog}.
\end{ex}

A more sophisticated application of the Bootstrapping Theorem shows that properties of the group~$\IZ$ can be bootstrapped to arithmetic groups such as~$\mathrm{SL}_d(\IZ)$.

\begin{ex}[Arithmetic groups] \label{ex:arithmetic groups}
		Let $G$ be a connected semisimple algebraic $\Q$-group. Let $\Gamma\subset G(\IQ)$ be an arithmetic group and suppose that~$\Gamma\in \qm$. 
		Let $r$ be the rational rank of~$G$. 
		If~$\IZ\in \sfB_m$ for all~$m\le n$, then~$\Gamma\in \sfB_{\min\{n, r-1\}}$.
		
		We sketch the argument following the proof of~\cite[Theorem~11.1]{ABFG21}. The group $G(\Q)$ and in particular $\Gamma$ acts on the spherical Tits building $\Delta$ of $G(\Q)$, which is an $(r-2)$-connected and cocompact $\Gamma$-CW complex. Let $\sigma$ be a cell of~$\Delta$. 
		The stabiliser $\Gamma_\sigma$ is the intersection of $\Gamma$ with a rational parabolic subgroup~$Q$ of~$G(\Q)$. The group~$\Gamma_\sigma$ 
		contains a non-trivial finitely generated torsion-free nilpotent normal subgroup~$N_\sigma$. By Example~\ref{ex:polycyclic} and our assumption on $\IZ$, we have $N_\sigma\in B_m$ for all~$m\le n$. Hence $\Gamma_\sigma\in \sfB_m$ for all~$m\le n$ by Corollary~\ref{cor:bootstrapping}~\ref{item:boot extension}. 
		By Theorem~\ref{thm:bootstrapping} we obtain that 
		$\Gamma\in \sfB_{\min\{n, r-1\}}$. 
\end{ex}

An equivariantly bootstrappable property~$\sfB^\qm_*$ of chain complexes over~$R$ may satisfy the following additional axiom for all~$n\in \IZ$:
\begin{enumerate}[label=\rm{(B-res)}]
	\item\label{ax:res} \emph{Restriction to finite index subgroups.} Let~$\Gamma\in \qm$, let~$\Delta$ be a finite index subgroup of~$\Gamma$, and let~$X$ be an $R\Gamma$-chain complex.
	If~$X\in \sfB^\Gamma_n$, then~$\res^\Gamma_\Delta X\in \sfB^\Delta_n$.
\end{enumerate}

\begin{cor}
\label{cor:commensurability}
	Let~$R$ be a ring and let~$\qm$ be an is-class of groups.
	Let~$\sfB_*$ be the bootstrappable property of~$\qm$ over~$R$ associated to an equivariantly bootstrappable property~$\sfB^\qm_*$ of chain complexes over~$R$ satisfying axiom~\ref{ax:res}.
	For all~$n\in \IZ$, the following hold:
	\begin{enumerate}[label=\enum]
		\item\label{item:boot finite sub} (Finite index subgroups).
		Let~$\Gamma\in \qm$ and let~$\Delta$ be a finite index subgroup of~$\Gamma$.
		If~$\Gamma\in \sfB_n$, then~$\Delta\in \sfB_n$;
		\item\label{item:boot finite over} (Finite index overgroups). Let~$\Gamma\in \qm$ and let~$\Delta$ be a finite index subgroup of~$\Gamma$.
		If~$\Delta\in \sfB_m$ for all~$m\le n$, then~$\Gamma\in \sfB_n$;
		\item\label{item:boot commensurated} (Commensurated subgroups). 
		Let~$\Gamma\in \qm$ and let $\Lambda$ be a commensurated subgroup of $\Gamma$ (i.e., for all~$\gamma\in \Gamma$, the intersection~$\Lambda\cap \gamma \Lambda\gamma^{-1}$ has finite index in~$\Lambda$ and in~$\gamma \Lambda\gamma^{-1}$).
		Suppose that~$\Gamma$ is of type~$\sfF_n$ and~$\Lambda$ is of type~$\sfF_{n-1}$. 
		If $\Lambda\in \sfB_m$ for all~$m\le n$, then $\Gamma\in \sfB_n$. 
	\end{enumerate}
	\textup{It follows from parts~\ref{item:boot finite sub} and~\ref{item:boot finite over} that the intersection~$\bigcap_{m\le n}\sfB_m$ of classes of groups is closed under commensurability.}
	\begin{proof}
		(i) If~$\Gamma$ lies in~$\sfB_n$, there exists a projective $R\Gamma$-resolution~$P$ of~$R$ lying in~$\sfB^\Gamma_n$.
		Then the restriction~$\res^\Gamma_\Delta P$ is a projective $R\Delta$-resolution of~$R$ lying in~$\sfB^\Delta_n$ by axiom~\ref{ax:res}, witnessing that~$\Delta$ lies in~$\sfB_n$.
		
		(ii) Suppose that~$\Delta$ lies in~$\sfB_m$ for all~$m\le n$.
		The normal core~$\Core_\Gamma(\Delta)\coloneqq \bigcap_{\gamma\in \Gamma}\gamma\Delta\gamma^{-1}$ of~$\Delta$ in~$\Gamma$ is a finite index normal subgroup of~$\Gamma$.
		Since~$\Core_\Gamma(\Delta)$ has finite index in~$\Delta$, we have~$\Core_\Gamma(\Delta)\in \sfB_m$ for all~$m\le n$ by part~\ref{item:boot finite sub}.
		Hence~$\Gamma$ lies in~$\sfB_n$ by Corollary~\ref{cor:bootstrapping}~\ref{item:boot finite extension}.
		
		(iii) We consider the Schlichting completion~$G$ of~$\Gamma$ relative to~$\Lambda$, that is,~$G$ is the closure of the image of the translation action $\tau\colon\Gamma\to\operatorname{Sym}(\Gamma/\Lambda)$ with respect to the topology of pointwise convergence. 
		Then~$G$ is a locally compact totally disconnected group and the closure~$U$ of~$\tau(\Lambda)$ is a compact-open subgroup. Furthermore, we have $U\cap\tau(\Gamma)=\tau(\Lambda)$~\cite[Section~2]{bonn+sauer}.
		Since~$\Gamma$ is of type~$\sfF_n$ and~$\Lambda$ is of type~$\sfF_{n-1}$,
		the group~$G$ is of type~$\sfF_n$ in the sense that there is a contractible $G$-CW-complex~$\Omega$ with compact-open stabilisers such that the $n$-skeleton is cocompact~\cite[Definition~3.1, Theorem~1.2]{bonn+sauer}. Every compact-open subgroup of~$G$ is commensurable with~$U$. Therefore the stabilisers of the $\Gamma$-CW-complex~$\res_\tau \Omega$ are commensurable with~$\Lambda$ and hence lie in~$\sfB_m$ for all~$m\le n$ by parts~\ref{item:boot finite sub} and~\ref{item:boot finite over}.
			    By Theorem~\ref{thm:bootstrapping}, the group~$\Gamma$ lies in~$\sfB_n$. 
	\end{proof}
\end{cor}

\begin{rem}
	Other examples of groups that admit appropriate actions to which the Bootstrapping Theorem~\ref{thm:bootstrapping} can be applied include mapping class groups, chain commuting groups, certain Artin groups~\cite{ABFG21}, outer automorphism groups of free products of~$\IZ/2$~\cite{GGH22}, mapping tori of polynomially growing automorphisms~\cite{AHK22,AGHK23}, and inner-amenable groups~\cite{Uschold22}.
\end{rem}

\section{Examples of bootstrappable properties}
\label{sec:examples}

We discuss several interesting equivariantly bootstrappable properties of chain complexes and relations between them.

\subsection{Algebraic finiteness properties}
\label{sec:FG}
For the basics on algebraic finiteness properties, we refer to Brown's book~\cite[Chapter~VIII]{Brown82}.
We work over an arbitrary ring~$R$.
\begin{defn}
\label{defn:FG}
	Let~$\Gamma$ be a group and let~$n\in \IZ$.
	The class~$\sfFG^\Gamma_n(R)$ consists of all $R\Gamma$-chain complexes~$X$ satisfying for all~$j\le n$ that the $R\Gamma$-module~$X_j$ is finitely generated.
\end{defn}
	Then a group~$\Gamma$ lies in~$\sfFG_n(R)$ (Definition~\ref{defn:pre-boot}) if and only if~$\Gamma$ is of type~$\sfFP_n(R)$.
	For example, it is easy to see that the group~$\IZ$ and finite groups lie in~$\sfFG_\infty(R)$.
	
\begin{lem}
	Let~$\qm$ be the is-class of all groups.
	The family~$\sfFG^\qm_*(R)$ is an equivariantly bootstrappable property of chain complexes over~$R$.
	Moreover, the family~$\sfFG^\qm_*(R)$ satisfies axiom~\ref{ax:res}.
	\begin{proof}
		First, for every group~$\Gamma$, the sequence~$\sfFG_*^\Gamma(R)$ is a bootstrappable property of $R\Gamma$-chain complexes because axioms~\ref{ax:pos}, \ref{ax:susp}, and~\ref{ax:cone} clearly hold.
		
		Second, the family~$\sfFG_*^\qm(R)$ satisfies axiom~\ref{ax:ind}.
		Indeed, let~$\Delta$ be a subgroup of~$\Gamma$ and let~$M$ be an $R\Delta$-module.
		If the $R\Delta$-module~$M$ is finitely generated, then the $R\Gamma$-module~$\ind_\Delta^\Gamma M$ is finitely generated.
		Applying this fact degreewise to the chain modules of chain complexes shows that axiom~\ref{ax:ind} holds.
		
		Third, the family~$\sfFG^\qm_*(R)$ satisfies axiom~\ref{ax:res}.
		Indeed, let~$\Delta$ be a finite index subgroup of~$\Gamma$ and let~$M$ be an $R\Gamma$-module.
		If the $R\Gamma$-module~$M$ is finitely generated, then the $R\Delta$-module~$\res^\Gamma_\Delta M$ is finitely generated.
		Applying this fact degreewise to the chain modules of chain complexes shows that axiom~\ref{ax:res} holds.
	\end{proof}
\end{lem}	
	
	We obtain the Bootstrapping Theorem~\ref{thm:bootstrapping} for~$\sfFG_*(R)$, which is a classical result~\cite[Proposition~1.1]{Brown87}.
	
\subsection{\texorpdfstring{$\ell^2$-Invariants}{l\texttwosuperior-Invariants}}

The basic principle behind $\ell^2$-invariants of groups~$\Gamma$ is
to adapt classical cohomological and spectral invariants to
coefficients in~$\ell^2\Gamma$ or the group von~Neumann algebra~$\calN
\Gamma$, which is the weak closure of~$\C \Gamma$ in~$B(\ell^2
\Gamma)$.  The key property of~$\calN \Gamma$ is that it admits an
accessible module theory and, in particular, a good notion of
dimension~$\dim_{\calN\Gamma}$.  The $\ell^2$-Betti numbers of~$\Gamma$
are then given by~$b^{(2)}_j (\Gamma) \coloneqq
\dim_{\calN \Gamma} H_j (\Gamma;\calN \Gamma)$.  For further
background on $\ell^2$-invariants we refer to L\"uck's
book~\cite{Lueck02}.

We will show that vanishing of $\ell^2$-homology, vanishing of $\ell^2$-Betti numbers, and lower bounds for Novikov--Shubin invariants are bootstrappable properties.
All cases follow from the general Lemma~\ref{lem:L2-chain} using that the respective property is preserved by submodules, quotients, extensions, induction, and restriction to finite index subgroups.
Throughout this section, we work over the ring~$R=\IZ$.

\begin{lem}
\label{lem:L2-chain}
	Let~$\Gamma$ be a group and let~$j\in \IZ$.
	The following hold:
	\begin{enumerate}[label=\enum]
		\item\label{item:L2-chain cone} Let~$f\colon X\to Y$ be a $\IZ\Gamma$-chain map.
		There is a short exact sequence of $\calN\Gamma$-modules
		\[
			0\to M\to H_j\bigl(\calN\Gamma\otimes_{\IZ\Gamma} \Cone(f)\bigr)\to L\to 0,
		\]
		where~$M$ is a quotient module of~$H_j(\calN\Gamma\otimes_{\IZ\Gamma} Y)$ and~$L$ is a submodule of~$H_{j-1}(\calN\Gamma\otimes_{\IZ\Gamma} X)$;
		\item\label{item:L2-chain ind} Let~$\Delta$ be a subgroup of~$\Gamma$ and let~$X$ be a $\IZ\Delta$-chain complex.
		There is an isomorphism of $\calN\Gamma$-modules
		\[
			H_j(\calN\Gamma\otimes_{\IZ\Gamma} \ind_\Delta^\Gamma X)\cong \calN\Gamma\otimes_{\calN\Delta} H_j(\calN\Delta\otimes_{\IZ\Delta} X);
		\]
		\item\label{item:L2-chain res} Let~$\Delta$ be a finite index subgroup of~$\Gamma$ and let~$X$ be a $\IZ\Gamma$-chain complex.
		There is an isomorphism of $\calN\Delta$-modules
		\[
			H_j(\calN\Delta\otimes_{\IZ\Delta} \res^\Gamma_\Delta X) \cong \res^{\calN\Gamma}_{\calN\Delta} H_j(\calN\Gamma\otimes_{\IZ\Gamma} X).
		\]
	\end{enumerate}
	\begin{proof}
		(i) By Lemma~\ref{lem:split ses}~\ref{item:split ses}, we have a short exact sequence of $\IZ\Gamma$-chain complexes
		\[
			0\to Y\to \Cone(f)\to \Sigma X\to 0
		\]
		which splits degreewise over~$\IZ\Gamma$.
		Hence the induced sequence of $\calN\Gamma$-chain complexes
		\[
			0\to \calN\Gamma\otimes_{\IZ\Gamma} Y \to \calN\Gamma\otimes_{\IZ\Gamma}\Cone(f) \to \calN\Gamma\otimes_{\IZ\Gamma}\Sigma X \to 0
		\]
		is (degreewise split) exact.
		Let $\delta_*\colon H_*(\calN\Gamma\otimes_{\IZ\Gamma} \Sigma X)\to H_{*-1}(\calN\Gamma\otimes_{\IZ\Gamma} Y)$ denote the connecting homomorphism of the associated long exact homology sequence. 
		For every~$j \in \Z$, we obtain an induced short exact sequence of $\calN\Gamma$-modules
 		 \[ 
		 	0\to \coker \delta_{j+1}
  			\to H_j\bigl(\calN \Gamma\otimes_{\Z\Gamma} \Cone(f)\bigr)
  			\to \ker \delta_j
 			 \to 0.
  		\]
		In particular, $\coker \delta_{j+1}$ is a quotient of $H_j(\calN\Gamma\otimes_{\IZ\Gamma} Y)$ and~$\ker \delta_j$ is a submodule of~$H_j(\calN\Gamma\otimes_{\Z\Gamma}\Sigma X)\cong H_{j-1}(\calN\Gamma\otimes_{\IZ\Gamma} X)$.
		
		(ii) There are canonical and natural isomorphisms of $\calN\Gamma$-chain complexes
		\[
			\calN\Gamma\otimes_{\IZ\Gamma}\ind_\Delta^\Gamma X
			=
			\calN\Gamma\otimes_{\IZ\Gamma} \IZ\Gamma\otimes_{\IZ\Delta} X
			\cong
			\calN\Gamma\otimes_{\calN\Delta} \calN\Delta\otimes_{\IZ\Delta}X.
		\]
		Since~$\calN\Gamma$ is flat over~$\calN\Delta$~\cite[Theorem~6.29]{Lueck02}, we obtain an isomorphism of $\calN\Gamma$-modules
 	 \[ 
	 	H_j(  \calN\Gamma\otimes_{\IZ\Gamma}\ind_\Delta^\Gamma X )
  		\cong
   		\calN\Gamma\otimes_{\calN\Delta}H_j(\calN\Delta\otimes_{\IZ\Delta}X). 
  	\]
	
	(iii) Since $\Delta$ has finite index in~$\Gamma$, there is an isomorphism of $\calN\Delta$-$\IZ\Gamma$-bimodules \cite[Proof of Theorem~3.7~(1)]{LueckReichSchick}
	\[
		\calN\Delta\otimes_{\IZ\Delta} \res^\Gamma_\Delta \IZ\Gamma \cong \res^{\calN\Gamma}_{\calN\Delta} \calN\Gamma.
	\]
	Then we have isomorphisms of $\calN\Delta$-modules
	\[
		\calN\Delta\otimes_{\IZ\Delta} \res^\Gamma_\Delta X
		\cong \calN\Delta\otimes_{\IZ\Delta} \res^\Gamma_\Delta \IZ\Gamma\otimes_{\IZ\Gamma} X
		\cong \res^{\calN\Gamma}_{\calN\Delta}\calN\Gamma \otimes_{\IZ\Gamma} X
	\]
	and hence
	\[
		H_j(\calN\Delta\otimes_{\IZ\Delta} \res^\Gamma_\Delta X)\cong \res^{\calN\Gamma}_{\calN\Delta}H_j(\calN\Gamma\otimes_{\IZ\Gamma} X).
	\]
	This finishes the proof.
	\end{proof}
\end{lem}

\subsubsection{$\ell^2$-Invisibility}
For an overview on~$\ell^2$-invisibility and its relevance for the Zero-in-the-spectrum Conjecture, we refer to~\cite[Chapter~12]{Lueck02}.
\begin{defn}
\label{defn:invisibility}
	Let~$\Gamma$ be a group and let~$n\in \IZ$.
	The class~$\sfI^\Gamma_n$ consists of all $\IZ\Gamma$-chain complexes~$X$ satisfying for all~$j\le n$ that
	\[
		H_j(\calN\Gamma\otimes _{\IZ\Gamma}X)=0.
	\]
\end{defn}
	Then a group~$\Gamma$ lies in~$\sfI_n$ (Definition~\ref{defn:pre-boot}) if and only if $H_j(\Gamma;\calN\Gamma)=0$ for all~$j\le n$. 
	Moreover, a group~$\Gamma$ of type~$\sfF_\infty$ lies in~$\sfI_\infty$ if and only if $H_j(\Gamma; \ell^2\Gamma) = 0$ for all~$j\in \IZ$ \cite[Lemma~12.3]{Lueck02}.
	The class~$\sfI_0$ is the class of non-amenable groups~\cite[Lemma~12.11~(4)]{Lueck02}. 
	To obtain groups in~$\sfI_n$,
	there is a product formula: If~$\Gamma_1\in \sfI_{n_1}$ and~$\Gamma_2\in \sfI_{n_2}$, then $\Gamma_1\times \Gamma_2\in \sfI_{n_1+n_2+1}$~\cite[Lemma~12.11~(3)]{Lueck02}.
	The groups in~$\sfI_\infty$ are said to be \emph{$\ell^2$-invisible}.
	
\begin{prop}
\label{prop:I}
	Let~$\qm$ be the is-class of all groups.
	The family~$\sfI_*^\qm$ is an equivariantly bootstrappable property of chain complexes over~$\IZ$.
	Moreover, the family~$\sfI^\qm_*$ satisfies axiom~\ref{ax:res}.
	\begin{proof}
	First, for every group~$\Gamma$, the sequence~$\sfI_*^\Gamma$ is a bootstrappable property of $\IZ\Gamma$-chain complexes.
	Indeed, the axioms~\ref{ax:pos} and~\ref{ax:susp} clearly hold and axiom~\ref{ax:cone} follows from Lemma~\ref{lem:L2-chain}~\ref{item:L2-chain cone}.
	Second, the family~$\sfI_*^\qm$ satisfies axiom~\ref{ax:ind} which follows from Lemma~\ref{lem:L2-chain}~\ref{item:L2-chain ind}.
	Third, the family~$\sfI^\qm_*$ satisfies axiom~\ref{ax:res} which follows from Lemma~\ref{lem:L2-chain}~\ref{item:L2-chain res}.
	\end{proof}
\end{prop}

We obtain the Bootstrapping Theorem~\ref{thm:bootstrapping} for~$\sfI_*$, which is implicit in the work of Sauer--Thumann~\cite[Proof of Theorem~1.1]{Sauer-Thumann14}.
This result has been used to provide examples of $\ell^2$-invisible groups of type~$\sfF_\infty$, given by certain local similarity groups~\cite{{Sauer-Thumann14}}.
It remains unknown whether there exists an $\ell^2$-invisible group of type~$\sfF$.

\subsubsection{$\ell^2$-Acyclicity}
For an overview on~$\ell^2$-acyclicity and its plentiful applications, we refer to~\cite[Section~7.1]{Lueck02}.

\begin{defn}
\label{defn:acyclicity}
	Let~$\Gamma$ be a group and let~$n\in \IZ$.
	The class~$\sfA^\Gamma_n$ consists of all $\IZ\Gamma$-chain complexes~$X$ satisfying for all~$j\le n$ that
	\[
		\dim_{\calN \Gamma}H_j(\calN\Gamma\otimes_{\IZ\Gamma} X)=0.
	\]
	Here~$\dim_{\calN\Gamma}$ is the von~Neumann dimension over~$\calN\Gamma$.
\end{defn}
	Then a group~$\Gamma$ lies in~$\sfA_n$ (Definition~\ref{defn:pre-boot}) if and only if~$b^{(2)}_j(\Gamma)=0$ for all~$j\le n$. 
	For example, the class~$\sfA_0$ is the class of infinite groups~\cite[Theorem~6.54~(8)]{Lueck02}.
	All infinite amenable groups lie in~$\sfA_\infty$~\cite[Corollary~6.75]{Lueck02}, and so does Thompson's group~$F$~\cite[Theorem~7.10]{Lueck02}.
	The groups in~$\sfA_\infty$ are said to be \emph{$\ell^2$-acyclic}.
	Clearly, the class~$\sfI_n$ is contained in~$\sfA_n$ and this inclusion is strict, as witnessed by infinite amenable groups. 
	
\begin{prop}
\label{prop:A}
	Let~$\qm$ be the is-class of all groups.
	The family~$\sfA_*^\qm$ is an equivariantly bootstrappable property of chain complexes over~$\IZ$.
	Moreover, the family~$\sfA^\qm_*$ satisfies axiom~\ref{ax:res}. 
	\begin{proof}
		First, for every group~$\Gamma$, the sequence~$\sfA_*^\Gamma$ is a bootstrappable property of $\IZ\Gamma$-chain complexes.
		Indeed, the axioms~\ref{ax:pos} and~\ref{ax:susp} clearly hold.
		Axiom~\ref{ax:cone} follows from Lemma~\ref{lem:L2-chain}~\ref{item:L2-chain cone} using that~$\dim_{\calN\Gamma}$ is additive with respect to short exact sequences~\cite[Theorem~1.12~(2)]{Lueck02}.
		Second, the family~$\sfA_*^\qm$ satisfies axiom~\ref{ax:ind}. This follows from Lemma~\ref{lem:L2-chain}~\ref{item:L2-chain ind} using that the von Neumann dimension is preserved by the induction functor~\cite[Lemma~1.24~(2)]{Lueck02}.
		Third, the family~$\sfA^\qm_*$ satisfies axiom~\ref{ax:res}.
		This follows from Lemma~\ref{lem:L2-chain}~\ref{item:L2-chain res} using that the von Neumann dimension is multiplicative with respect to restriction to finite index subgroups~\cite[Theorem~1.12~(6)]{Lueck02}.
	\end{proof}
\end{prop}

We obtain the Bootstrapping Theorem~\ref{thm:bootstrapping} for~$\sfA_*$, which was proved (in a slightly weaker form) by Jo~\cite[Theorem~3.5]{Jo07}. 

\subsubsection{Novikov--Shubin invariants and capacity}

Novikov--Shubin invariants (and their reciprocals, the capacities)
are spectral invariants that measure
the difference between $\ell^2$-acyclicity and $\ell^2$-in\-vis\-i\-bil\-i\-ty.
Capacity and Novikov--Shubin invariants take values in the
extended ranges of numbers
\begin{align*}
  \llbracket 0, \infty]
  \coloneqq \{0^-\} \sqcup [0,\infty]
  \qand 
  [0,\infty\rrbracket
  \coloneqq [0,\infty] \sqcup \{\infty^+\},
\end{align*}
respectively.
Here $\llbracket 0,\infty\rrbracket \coloneqq \{0^-\} \sqcup [0,\infty]
\sqcup \{\infty^+\}$ carries the ordering that extends the
usual ordering on~$[0,\infty]$ by~$0^- < 0$ and $\infty < \infty^+$.
Moreover, arithmetic operations are extended as expected and
\[ \frac1{0^-} \coloneqq \infty^+
\qand
\frac1{\infty^+} \coloneqq 0^-.
\]
We briefly indicate the definition in the case of groups~$\Gamma$ of
type~$\sfFP_\infty$: For~$j\in \IN$, let $\partial_j$ denote the
boundary operator of~$\calN \Gamma \otimes_{\IC \Gamma} X$, where $X$
is a finite type $\IC\Gamma$-resolution of~$\IC$, and let $f$ be the
spectral density function of~$\partial_j^* \partial_j$. The $j$-th
Novikov--Shubin invariant~$\alpha_j(\Gamma) \in [0,\infty\rrbracket$
  of~$\Gamma$ is defined as
\[ \alpha_j(\Gamma)
\coloneqq
\begin{cases}
\liminf_{\lambda \searrow 0}
\frac{\ln (f(\lambda^2) - f(0))}{\ln(\lambda)}
& \text{if $f(\lambda^2) > f(0)$ for all~$\lambda >0$;}
\\
\infty^+
& \text{otherwise}.
\end{cases}
\]
In particular, if $\Gamma$ is $\ell^2$-acyclic, then
$\alpha_j(\Gamma) = \infty^+$ for all~$j \in \N$ if and only if
$\Gamma$ is $\ell^2$-invisible~\cite[Lemma~6.98]{Lueck02}.

One may extend Novikov--Shubin invariants to $\calN
\Gamma$-modules~\cite{lueck_NS} (which requires the notion of
measurable and cofinal-measurable modules) and then $\alpha_j(\Gamma)
= \alpha_{\calN \Gamma}(H_{j-1}(\Gamma;\calN \Gamma))$. 
In the following,
we will use the setup developed by
L\"uck--Reich--Schick~\cite{LueckReichSchick}.  In particular, we will
phrase the properties and arguments in terms of capacity~$c_{\calN
  \Gamma} \coloneqq 1/\alpha_{\calN \Gamma}$, which is reciprocal to the
Novikov--Shubin invariants.

\begin{defn}
\label{defn:NovikovShubin}
  Let $\Gamma$ be a group, let $n \in \Z$, and let $\kappa \in \llbracket 0,\infty]$.
  \begin{itemize}
  \item The class~$\sfCM_n^\Gamma$ consists of all $\IZ\Gamma$-chain complexes~$X$ satisfying for all~$j\le n$ that the $\calN \Gamma$-module~$H_j(\calN\Gamma
    \otimes_{\Z\Gamma} X)$ is
    cofinal-measurable~\cite[Definition~2.1]{LueckReichSchick};
  \item The class~$\sfC(\leq \kappa)_n^\Gamma$ consists of all $\IZ\Gamma$-chain complexes~$X\in \sfCM^\Gamma_n$
   satisfying for all~$j\le n$ that
    \[ c_{\calN\Gamma} \bigl( H_j(\calN \Gamma\otimes_{\Z\Gamma} X)\bigr)
    \leq \kappa.
    \]
    Here~$c_{\calN\Gamma}$ is the capacity of
    $\calN \Gamma$-modules~\cite[Section~2]{LueckReichSchick},
    i.e., the reciprocal of the Novikov--Shubin invariant;
    \item The class~$\sfC(<\kappa)^\Gamma_n$ is defined analogously to~$\sfC(\le \kappa)^\Gamma_n$ by replacing the inequality with a strict inequality.
  \end{itemize}
\end{defn}

All cofinal-measurable modules have von~Neumann dimension~$0$. 
In particular, the class~$\sfCM^\Gamma_n$ is a subclass of~$\sfA^\Gamma_n$.
Every $\ell^2$-acyclic group of type~$\sfFP_\infty$
lies in~$\sfCM_\infty$ (Definition~\ref{defn:pre-boot}).  
An $\ell^2$-acyclic group~$\Gamma$ of type~$\sfF_\infty$ lies in~$\sfC(<\infty)_n$ if and only if the Novikov--Shubin invariants of~$\Gamma$
satisfy~$\alpha_j(\Gamma)>0$ for all~$j \leq n$ (in the indexing of Novikov--Shubin invariants of boundary operators, not of Laplacians).
For some calculations of the capacity of groups, see~\cite[Section~3]{LueckReichSchick}.
By construction of the capacity/Novikov--Shubin invariants,
for every group~$\Gamma$,
we have~$\sfI_n^\Gamma = \sfC(\le 0^-)_n^\Gamma$.

\begin{prop}
\label{prop:C}
  Let~$\qm$ be the is-class of all groups.
  The following hold:
  \begin{enumerate}[label=\enum]
  \item The family~$\sfCM_*^\qm$ is an equivariantly bootstrappable property of chain complexes over~$\IZ$ and satisfies axiom~\ref{ax:res};  
  \item The families $\sfC(\le 0^-)^\qm_*$, $\sfC(\le 0)^\qm_*$, and~$\sfC(<\infty)^\qm_*$ are equivariantly bootstrappable properties of chain complexes over~$\IZ$ and satisfy axiom~\ref{ax:res}.
  \end{enumerate}
\begin{proof}
		(i) First, for every group~$\Gamma$, the sequence~$\sfCM_*^\Gamma$ is a bootstrappable property of $\IZ\Gamma$-chain complexes.
		Indeed, the axioms~\ref{ax:pos} and~\ref{ax:susp} clearly hold.
		Axiom~\ref{ax:cone} follows from Lemma~\ref{lem:L2-chain}~\ref{item:L2-chain cone} using that for $\calN\Gamma$-modules, being cofinal-measurable is stable under taking submodules, quotients,
  and extensions~\cite[Lemma~2.11]{LueckReichSchick}.
		Second, the family~$\sfCM_*^\qm$ satisfies axiom~\ref{ax:ind}. This follows from Lemma~\ref{lem:L2-chain}~\ref{item:L2-chain ind} using that being cofinal-measurable is preserved by the induction functor~\cite[Lemma~2.12~(1)]{LueckReichSchick}.
		Third, the family~$\sfCM^\qm_*$ satisfies axiom~\ref{ax:res}.
		This follows from Lemma~\ref{lem:L2-chain}~\ref{item:L2-chain res} using that being cofinal-measurable is preserved by restriction to finite index subgroups~\cite[Lemma~2.12~(3)]{LueckReichSchick}.

		(ii) Let~$\kappa\in \{0^-,0\}$. We will only treat the family~$\sfC(\leq \kappa)^\qm_*$.
		The family~$\sfC(<\infty)^\qm_*$ can be handled analogously.
		
		First, for every group~$\Gamma$, the sequence~$\sfC(\le \kappa)_*^\Gamma$ is a bootstrappable property of $\IZ\Gamma$-chain complexes.
		Indeed, the axioms~\ref{ax:pos} and~\ref{ax:susp} clearly hold.
		Axiom~\ref{ax:cone} follows from Lemma~\ref{lem:L2-chain}~\ref{item:L2-chain cone} using that for cofinal-measurable $\calN\Gamma$-modules, capacity is monotone for submodules and quotients, and subadditive with respect to short exact sequences~\cite[Theorem~2.7~(1)]{LueckReichSchick}.
		Here we use that our choices of~$\kappa$ satisfy~$\kappa+\kappa\le \kappa$.
		Second, the family~$\sfC(\le \kappa)_*^\qm$ satisfies axiom~\ref{ax:ind}. This follows from Lemma~\ref{lem:L2-chain}~\ref{item:L2-chain ind} using that capacity is preserved by the induction functor~\cite[Lemma~2.12~(1)]{LueckReichSchick}.
		Third, the family~$\sfC(\leq \kappa)^\qm_*$ satisfies axiom~\ref{ax:res}.
		This follows from Lemma~\ref{lem:L2-chain}~\ref{item:L2-chain res} using that capacity is preserved by restriction to finite index subgroups~\cite[Lemma~2.12~(3)]{LueckReichSchick}.
\end{proof}
\end{prop}

We obtain corresponding instances of the
Bootstrapping Theorem~\ref{thm:bootstrapping} for each of the families~$\sfC(\le 0^-)_*$, $\sfC(\le 0)_*$, and~$\sfC(<\infty)_*$. 
The family~$\sfC(\le 0^-)_*$ recovers the case of~$\sfI_*$. 
Translated
to Novikov--Shubin invariants, these three cases correspond to the constraints of having cofinal-measurable $\ell^2$-homology whose Novikov--Shubin invariants are
\begin{itemize}
\item equal to~$\infty^+$;
\item equal to $\infty$ or~$\infty^+$;
\item positive;
\end{itemize}
respectively.
It would be interesting if also the constraint ``the Novikov--Shubin
invariants are~\mbox{$\geq 1$}'' could be proved to be (equivariantly)
bootstrappable. However, the known inheritance results for 
Novikov--Shubin invariants of short exact sequences are not
strong enough to handle this case.

\begin{rem}
It was conjectured by Lott--L\"uck~\cite[Conjecture~7.2]{lott-lueck} that, for every group~$\Gamma$, the Novikov--Shubin invariants of every $\Z\Gamma$-chain complex of finitely generated free $\Z\Gamma$-modules are positive. 
Grabowski found a counterexample by con\-structing a solvable group~$\Gamma$, which is not of type~$\sfF_\infty$, and an element~$c\in\Z\Gamma$ with vanishing Novikov--Shubin invariant~\cite{grabowski}. The corresponding chain complex is concentrated in two degrees with the differential being multiplication by~$c$. However, the following weakening of Lott--L\"uck's conjecture remains open: Are the Novikov--Shubin invariants of a group of type~$\sfF_\infty$ positive? 

We provide some evidence:
Since the group~$\IZ$ lies in~$ \sfC(<\infty)_\infty$~\cite[Theorem~3.7]{LueckReichSchick},
it follows from the Bootstrapping Theorem~\ref{thm:bootstrapping} for~$\sfC(<\infty)_*$ that infinite elementary amenable groups of type~$\sfFP_\infty$ lie in~$\sfC(<\infty)_\infty$ (Example~\ref{ex:elementary:ame}).
In particular, the Novikov--Shubin invariants of these groups are positive.
\end{rem}

\begin{rem}
	In Definition~\ref{defn:invisibility}, Definition~\ref{defn:acyclicity}, and Definition~\ref{defn:NovikovShubin} we consider vanishing resp.\ upper bounds in all degrees~$\le n$.
	In particular, for every sequence~$\sfB_*\in\{\sfI_*,\sfA_*,\sfCM_*,\sfC(\le \kappa)_*,\sfC(<\kappa)_*\}$, the class~$\sfB_n$ is contained in~$\sfB_m$ for all~$m\le n$.
	Instead, one may also consider the vanishing of $\ell^2$-homology, vanishing of von Neumann dimension, and upper bounds for capacity in a fixed degree~$n$ only.
	One obtains corresponding equivariantly bootstrappable properties of chain complexes over~$\IZ$ by the same proofs as Proposition~\ref{prop:I}, Proposition~\ref{prop:A}, and Proposition~\ref{prop:C}.
\end{rem}

\subsubsection{$\ell^2$-Torsion}
\label{sec:L2-torsion}
	For an overview on $\ell^2$-torsion, we refer to~\cite[Chapter~3]{Lueck02}.
The $\ell^2$-torsion is an $\ell^2$-version of Reidemeister torsion.
	It is a secondary invariant that is defined for $\ell^2$-acyclic chain complexes.
	Let~$X$ be a Hilbert $\calN\Gamma$-chain complex of finite length satisfying for every~$j\in \IZ$ that $\dim_{\calN\Gamma}(X_j)<\infty$, $b^{(2)}_j(X)=0$, and $\det_{\calN\Gamma}(\partial_j)>0$.
	Then the \emph{$\ell^2$-torsion~$\rho^{(2)}(X)$ of~$X$} is defined as
	\[
		\rho^{(2)}(X;\calN\Gamma)\coloneqq -\sum_{j\in \IZ}(-1)^j\cdot \log\bigl(\operatorname{det}_{\calN\Gamma}(\partial_j)\bigr) \in \mathbb{R}.
	\]
	Here~$\det_{\calN\Gamma}$ is the Fuglede--Kadison determinant~\cite[Definition~3.11]{Lueck02}.
	
	A group~$\Gamma$ is said to be \emph{of~$\det\ge 1$-class} if for every matrix~$A\in M(k\times l;\IZ\Gamma)$, the map $r_A^{(2)}\colon (\ell^2\Gamma)^k\to (\ell^2\Gamma)^l$ given by (right-)multiplication with~$A$ satisfies 
	\[
		\operatorname{det}_{\calN\Gamma}(r_A^{(2)})\ge 1.
	\]
	Conjecturally, all groups are of $\det\ge 1$-class~\cite[Conjecture~13.2]{Lueck02}. It is known that sofic groups are of $\det\ge 1$-class~\cite{Elek-Szabo05}.
	
	Recall that a group~$\Gamma$ is \emph{of type~$\mathsf{FL}$} if the trivial $\IZ\Gamma$-module~$\IZ$ admits a finite free $\IZ\Gamma$-resolution.
	The $\ell^2$-torsion of groups is well-defined for groups~$\Gamma$ that are of type~$\mathsf{FL}$, of $\det\ge 1$-class, and $\ell^2$-acyclic as
	\[
		\rho^{(2)}(\Gamma)\coloneqq \rho^{(2)}(\ell^2\Gamma\otimes_{\IZ\Gamma} P;\calN\Gamma),
	\]
	where~$P$ is a finite free $\IZ\Gamma$-resolution of~$\IZ$.
	Non-trivial amenable groups of type~$\mathsf{FL}$ have vanishing $\ell^2$-torsion~\cite[Theorem~1.3]{Li-Thom}.
	
	L\"uck conjectured an approximation result for~$\ell^2$-torsion in terms of torsion homology growth (Definition~\ref{defn:growth intro}) ~\cite[Conjecture~1.11~(3)]{Lueck13}:
	\[
		\rho^{(2)}(\Gamma)=\sum_{j\in \IZ}(-1)^j\cdot  \widehat{t}_j(\Gamma,\Lambda_*).
	\]
	The vanishing of $\ell^2$-torsion does not fit into our framework of bootstrappable properties because $\ell^2$-torsion depends on the entire chain complex, which causes complications for axiom~\ref{ax:pos}.
	Nevertheless, a suitable Bootstrapping Theorem for cocompact actions on acyclic spaces can be proved.
	We learnt about a Bootstrapping Theorem for (not necessarily vanishing) $\ell^2$-torsion from Hughes--L\"uck~\cite{Hughes-Lueck25},
	who also give many applications to groups and group automorphisms.
	The result of Hughes--L\"uck applies to groups that are not necessarily of type~$\sfFL$ and their proof relies on topological methods.
	Below we give an alternative proof of their result for groups of type~$\sfFL$ that is completely algebraic.	
	
	\begin{thm}[{Hughes--L\"uck~\cite[Theorem~3.7]{Hughes-Lueck25}}]
		Let~$\Gamma$ be a group and let~$\Omega$ be a $\Gamma$-CW-complex such that the following hold:
		\begin{enumerate}[label=\enum]
			\item $\Omega$ is acyclic;
			\item The quotient~$\Gamma\backslash \Omega$ is compact;
			\item The group~$\Gamma$ is of $\det\ge 1$-class;
			\item For every cell~$\sigma$ of~$\Omega$, the stabiliser~$\Gamma_\sigma$ is of type~$\mathsf{FL}$, of $\det\ge 1$-class, and $\ell^2$-acyclic.
		\end{enumerate}
		Then~$\Gamma$ is of type~$\sfFL$, $\ell^2$-acyclic, and we have
		\[
			\rho^{(2)}(\Gamma)=\sum_{j\in \IN}\sum_{\sigma\in S_j} (-1)^j\cdot \rho^{(2)}(\Gamma_\sigma),
		\]
		where~$S_j$ is a finite set of representatives for all $\Gamma$-orbits of~$j$-cells of~$\Omega$.
		\begin{proof}
			The group~$\Gamma$ is $\ell^2$-acyclic by the Bootstrapping Theorem~\ref{thm:bootstrapping} for~$\sfA_*$ (Proposition~\ref{prop:A}).
			To prove the other statements, we closely follow the proof of Theorem~\ref{thm:bootstrapping}.
			Let~$X$ be the cellular $\IZ\Gamma$-chain complex of~$\Omega$.
			For all~$j\in \IZ$, the $\IZ\Gamma$-chain module~$X_j$ is of the form
			\[
				X_j\cong \bigoplus_{\sigma\in S_j}\IZ[\Gamma/\Gamma_\sigma].	
			\]
			For every~$\sigma\in S_j$, let~$P^\sigma$ be a finite free $\IZ[\Gamma_\sigma]$-resolution of the trivial $\IZ[\Gamma_\sigma]$-module~$\IZ$.
			Then
			\[
				P^j\coloneqq \bigoplus_{\sigma\in S_j}\operatorname{ind}_{\Gamma_\sigma}^\Gamma P^\sigma
			\]
			is a finite free $\IZ\Gamma$-resolution of~$X_j$.
			We obtain a \emph{free} replacement~$\widehat{X}$ of~$X$ with respect to the resolutions~$(P^j)_j$ analogously to Proposition~\ref{prop:blow-up}.
			Then~$\widehat{X}$ is a finite free $\IZ\Gamma$-resolution of the trivial~$\IZ\Gamma$-module~$\IZ$.
			This proves that~$\Gamma$ is of type~$\mathsf{FL}$.
			
			Moreover, $\widehat{X}$ admits a filtration~$(\widehat{X}^{[k]})_{k\ge 0}$ of finite length such that~$\widehat{X}^{[k]}$ is the mapping cone of a chain map $\Sigma^{k-1}P^k\to \widehat{X}^{[k-1]}$.
			In particular, for ever~$k\ge 0$, we have a short exact sequence of $\IZ\Gamma$-chain complexes
			\[
				0\to \widehat{X}^{[k-1]}\to \widehat{X}^{[k]}\to \Sigma^kP^k\to 0
			\]
			which splits degreewise over~$\IZ\Gamma$. 
			Hence the sequence of Hilbert $\calN\Gamma$-chain complexes
			\begin{equation}
			\label{eqn:L2torsion}
				0\to \ell^2\Gamma\otimes_{\IZ\Gamma} \widehat{X}^{[k-1]}\to \ell^2\Gamma\otimes_{\IZ\Gamma} \widehat{X}^{[k]}\to \ell^2\Gamma\otimes_{\IZ\Gamma} \Sigma^kP^k\to 0
			\end{equation}
			is exact and splits degreewise over~$\calN\Gamma$. 
			Since the Hilbert $\calN(\Gamma_\sigma)$-chain complex $\ell^2\Gamma_\sigma\otimes_{\IZ\Gamma_\sigma}P^\sigma$ is $\det$-$\ell^2$-acyclic~\cite[Definition~3.29]{Lueck02} by assumption on~$\Gamma_\sigma$, the right term in the sequence~\eqref{eqn:L2torsion} is $\det$-$\ell^2$-acyclic and satisfies
			\begin{align*}
				\rho^{(2)}(\ell^2\Gamma\otimes_{\IZ\Gamma} \Sigma^k P^k;\calN\Gamma)
				&= (-1)^k\cdot \rho^{(2)}(\ell^2\Gamma\otimes_{\IZ\Gamma} P^k;\calN\Gamma)
				\\
				&= (-1)^k\cdot \sum_{\sigma\in S_k} \rho^{(2)}\bigl(\ell^2(\Gamma_\sigma)\otimes_{\IZ[\Gamma_\sigma]} P^\sigma;\calN(\Gamma_\sigma)\bigr)
				\\
				&= (-1)^k\cdot \sum_{\sigma\in S_k} \rho^{(2)}(\Gamma_\sigma). 
			\end{align*}
			For the second equality we used that $\ell^2$-torsion is compatible with induction and direct sums~\cite[Theorem~3.35~(1) and~(8)]{Lueck02}.
			By induction on~$k$, the left term in the sequence~\eqref{eqn:L2torsion} is $\det$-$\ell^2$-acyclic and satisfies
			\[
				\rho^{(2)}(\ell^2\Gamma\otimes_{\IZ\Gamma} \widehat{X}^{[k-1]};\calN\Gamma)
				= \sum_{j=0}^{k-1}\sum_{\sigma\in S_j} (-1)^j\cdot \rho^{(2)}(\Gamma_\sigma).
			\]
			Hence by~\cite[Theorem~3.35~(1)]{Lueck02}, the middle term in the sequence~\eqref{eqn:L2torsion} is $\det$-$\ell^2$-acyclic and satisfies
			\begin{align*}
				\rho^{(2)}(\ell^2\Gamma\otimes_{\IZ\Gamma}\widehat{X}^{[k]};\calN\Gamma)
				&= \rho^{(2)}(\ell^2\Gamma\otimes_{\IZ\Gamma} \widehat{X}^{[k-1]};\calN\Gamma) + \rho^{(2)}(\ell^2\Gamma\otimes_{\IZ\Gamma}\Sigma^kP^k;\calN\Gamma)
				\\
				&= \sum_{j=0}^k\sum_{\sigma\in S_j} (-1)^j\cdot \rho^{(2)}(\Gamma_\sigma).
			\end{align*}
			Since the filtration~$(\widehat{X}^{[k]})_k$ is of finite length and~$\widehat{X}$ is a finite free $\IZ\Gamma$-resolution of~$\IZ$, this yields the desired computation of the $\ell^2$-torsion of~$\Gamma$.
		\end{proof}
	\end{thm}
	
\subsection{Vanishing of (torsion) homology growth}
For an overview on (torsion) homology growth, we refer to L{\"u}ck's survey~\cite{Lueck16_survey}.
We work over the ring~$R=\IZ$.

Let~$\Gamma$ be a group and let~$\Lambda$ be a subgroup of~$\Gamma$.
We consider the functor~$(-)_\Lambda$ that associates to a $\IZ\Gamma$-module~$M$ the $\IZ$-module
$M_\Lambda\coloneqq \IZ\otimes_{\IZ\Lambda}\res^\Gamma_\Lambda M$ of $\Lambda$-coinvariants.
A $\IZ\Gamma$-basis of a free $\IZ\Gamma$-module~$M$ induces a $\IZ$-basis of the free $\IZ$-module~$M_\Lambda$. 
We also denote by~$(-)_\Lambda$ the induced functor of chain complexes that associates to a $\IZ\Gamma$-chain complex~$X$ the $\IZ$-chain complex 
\[
	X_\Lambda\coloneqq \IZ\otimes_{\IZ\Lambda} \res^\Gamma_\Lambda X
\]
of $\Lambda$-coinvariants.
We record some basic properties:
\begin{lem}
\label{lem:coinvariants}
	Let~$\Gamma$ be a group.
	The following hold:
	\begin{enumerate}[label=\enum]
		\item\label{item:coinvariants module} 
		Let~$\Lambda$ be a finite index normal subgroup of~$\Gamma$.
		Let~$\Delta$ be a subgroup of~$\Gamma$ and let~$M$ be a $\IZ\Delta$-module.
		Then there is an isomorphism of $\IZ$-modules
		\[
			(\ind_\Delta^\Gamma M)_\Lambda\cong \bigoplus_{\frac{[\Gamma:\Lambda]}{[\Delta:\Lambda\cap \Delta]}} M_{\Lambda\cap \Delta}
		\]
		which is natural in~$M$;
		\item\label{item:coinvariants chain} 
		Let~$\Lambda$ be a finite index normal subgroup of~$\Gamma$.
		Let~$\Delta$ be a subgroup of~$\Gamma$ and let~$X$ be a $\IZ\Delta$-chain complex.
		Then there is a natural isomorphism of $\IZ$-chain complexes
		\[
			(\ind_\Delta^\Gamma X)_\Lambda\cong \bigoplus_{\frac{[\Gamma:\Lambda]}{[\Delta:\Lambda\cap \Delta]}} X_{\Lambda\cap \Delta};
		\]
		\item\label{item:coinvariants map} Let~$f\colon M\to L$ be a $\IZ\Gamma$-chain map between based free $\IZ\Gamma$-modules.
		Let~$\Lambda$ be a subgroup of~$\Gamma$.
		Then the induced map~$f_\Lambda\colon M_\Lambda\to L_\Lambda$ of based free $\IZ$-modules satisfies $\|f_\Lambda\|\le \|f\|$.
		Here~$\|\cdot\|$ denotes the operator-norm of maps between based free $\IZ$-modules equipped with the $\ell^1$-norm.
	\end{enumerate}
	\begin{proof}
		(i) We have natural isomorphisms of $\IZ$-modules 
		\begin{align*}
			(\ind_\Delta^\Gamma M)_\Lambda 
			&\cong \IZ[\Lambda\backslash \Gamma]\otimes_{\IZ\Gamma} (\ind_\Delta^\Gamma M)
			\cong \IZ[\Lambda\backslash \Gamma]\otimes_{\IZ\Gamma}(\IZ\Gamma\otimes_{\IZ\Delta} M)
			\\
			&\cong \IZ[\Lambda\backslash \Gamma]\otimes_{\IZ\Delta} M
			\cong \bigoplus_{\frac{[\Gamma:\Lambda]}{[\Delta:\Lambda\cap \Delta]}}\IZ[\Lambda\cap \Delta\backslash \Delta]\otimes_{\IZ\Delta}M
			\cong \bigoplus_{\frac{[\Gamma:\Lambda]}{[\Delta:\Lambda\cap \Delta]}} M_{\Lambda\cap \Delta}.
		\end{align*}
		For the second to last isomorphism it is crucial that~$\Lambda$ is normal in~$\Gamma$.
		
		(ii) This follows directly from part~\ref{item:coinvariants module}.
		
		(iii) Consider the diagram of based free $\IZ$-modules
		\[\begin{tikzcd}
			M\ar{r}{f}\ar{d}[swap]{p_M} & L\ar{d}{p_L} \\
			M_\Lambda\ar{r}[swap]{f_\Lambda} & L_\Lambda
		\end{tikzcd}\]
		where~$p_M$ and~$p_L$ are the obvious projections.
		For~$m\in M_\Lambda$, there exists~$\widetilde{m}\in M$ with~$p_M(\widetilde{m})=m$ and~$|\widetilde{m}|_1=|m|_1$. 
		Then we have
		\[
			|f_\Lambda(m)|_1=|f_\Lambda\circ p_M(\widetilde{m})|_1=|p_L\circ f(\widetilde{m})|_1\le \|p_L\|\cdot \|f\|\cdot |\widetilde{m}|_1\le \|f\|\cdot |m|_1
		\]
		and hence~$\|f_\Lambda\|\le \|f\|$.
	\end{proof}
\end{lem}

In the following, we consider the is-class of residually finite groups.

\begin{defn}
\label{defn:H}
	Let~$\Gamma$ be a residually finite group, let~$n\in \IZ$, and let~$\IF$ be a field.
	The class~$\sfH^\Gamma_n(\IF)$ consists of all $\IZ\Gamma$-chain complexes~$X$ satisfying for all residual chains~$\Lambda_*$ in~$\Gamma$ and all~$j\le n$ that the Betti number gradient 
	\[
		\widehat{b}_j(X,\Lambda_*;\IF)\coloneqq
		\limsup_{i\to \infty}\frac{\dim_\IF H_j(\IF\otimes_\IZ X_{\Lambda_i})}{[\Gamma:\Lambda_i]}
	\]
        is equal to~$0$.        
\end{defn}

	Then the class~$\sfH_n(\IF)$ (Definition~\ref{defn:pre-boot}) consists of all residually finite groups~$\Gamma$ satisfying that $\widehat{b}_j(\Gamma, \Lambda_*;\IF)=0$ for all residual chains~$\Lambda_*$ in~$\Gamma$ and all~$j\le n$ (Definition~\ref{defn:growth intro}).
	Clearly, the class~$\sfH_0(\IF)$ is the class of residually finite infinite groups.
	It is easy to see that the group~$\IZ$ lies in~$\sfH_\infty(\IF)$.
	More generally, residually finite infinite amenable groups of type~$\sfFP_{n}$ lie in~$\sfH_n(\IF)$ (Corollary~\ref{cor:amenable T}).
	
	For a group~$\Gamma$ of type~$\sfFP_n$, the Universal Coefficient Theorem shows that if~$\Gamma$ lies in~$\sfH_n(\IF)$ for some field~$\IF$, then~$\Gamma$ lies in~$\sfH_n(\IQ)$.
	If~$\Gamma$ is residually finite and of type~$\sfFP_{j+1}$, 
	then $\widehat{b}_j(\Gamma,\Lambda_*;\IQ)$ coincides with the $\ell^2$-Betti number~$b^{(2)}_j(\Gamma)$ by L\"uck's Approximation Theorem~\cite{Lueck94_approx}. 
	In particular, a residually finite group~$\Gamma$ of type~$\sfFP_{n+1}$ lies in~$\sfH_n(\IQ)$ if and only if~$\Gamma$ lies in~$\sfA_n$.
	
\begin{rem}
	The topological interpretation of Definition~\ref{defn:H} is as follows:
	Let~$\Omega$ be a CW-complex with residually finite fundamental group~$\Gamma$.
	Let~$\Lambda_*$ be a residual chain in~$\Gamma$ and denote by~$\Omega_i$ the covering of~$\Omega$ associated to the subgroup~$\Lambda_i$ of~$\Gamma$.
	Let~$X$ be the cellular $\IZ\Gamma$-chain complex of the universal covering~$\widetilde{\Omega}$.
	Then~$X_{\Lambda_i}$ is isomorphic to the cellular $\IZ$-chain complex of~$\Omega_i$ and hence
	\[
		\widehat{b}_j(X,\Lambda_*;\IF) = \limsup_{i\to \infty}\frac{\dim_\IF H_j(\Omega_i;\IF)}{[\Gamma:\Lambda_i]}.
	\]
	Vanishing of~$\widehat{b}_j(X,\Lambda_*;\IF)$ means that the $j$-th $\IF$-Betti number of the coverings~$(\Omega_i)_i$ grows sublinearly with respect to the index~$[\Gamma:\Lambda_i]$ as~$i\to \infty$.
\end{rem}
		
\begin{lem}\label{lem:Betti}
	Let~$\Gamma$ be a residually finite group and let~$\Lambda_*$ be a residual chain in~$\Gamma$.
	Let~$\IF$ be a field and let~$j\in \IZ$.
	The following hold:
	\begin{enumerate}[label=\enum]
		\item\label{item:Betti cone} Let $f\colon X\to Y$ be a $\IZ\Gamma$-chain map.
		Then
		\[
			\widehat{b}_j\bigl(\Cone(f),\Lambda_*;\IF\bigr)\le \widehat{b}_j(Y,\Lambda_*;\IF) + \widehat{b}_{j-1}(X,\Lambda_*;\IF);
		\]
		\item\label{item:Betti ind} Let~$\Delta$ be a subgroup of~$\Gamma$ and let~$X$ be a $\IZ\Delta$-chain complex.
		Then
		\[
			\widehat{b}_j(\ind_\Delta^\Gamma X,\Lambda_*;\IF) = \widehat{b}_j(X,\Lambda_*\cap \Delta;\IF).
		\]
	\end{enumerate}
	\begin{proof}
		(i) The sequence of $\IZ\Gamma$-chain complexes
		\[
			0 \to Y \to \Cone(f) \to \Sigma X\to 0
		\]
		splits degreewise over~$\IZ\Gamma$  by Lemma~\ref{lem:split ses}~\ref{item:split ses}.
		Then, for every subgroup~$\Lambda$ of~$\Gamma$, the induced sequence of $\IZ$-chain complexes
		\[
			0\to Y_\Lambda\to \Cone(f)_\Lambda\to (\Sigma X)_\Lambda\to 0
		\]	 
		is (degreewise split) exact.
		By the long exact sequence in homology, for all~$j\in \IZ$, we have
		\[
			\dim_\IF H_j\bigl(\IF\otimes_\IZ \Cone(f)_\Lambda\bigr) \le \dim_\IF H_j\bigl(\IF\otimes_\IZ Y_\Lambda\bigr) + \dim_\IF H_j\bigl(\IF\otimes_\IZ (\Sigma X)_\Lambda \bigr).
		\]
		Since~$(\Sigma X)_\Lambda\cong \Sigma(X_\Lambda)$, the claim follows using the suspension isomorphism.
		
		(ii) For a finite index normal subgroup~$\Lambda$ of~$\Gamma$, we have an isomorphism of $\IZ$-chain complexes
		\[
			(\ind_\Delta^\Gamma X)_\Lambda\cong \bigoplus_{\frac{[\Gamma:\Lambda]}{[\Delta:\Lambda\cap \Delta]}}X_{\Lambda\cap \Delta}
		\]
		by Lemma~\ref{lem:coinvariants}~\ref{item:coinvariants chain}.
		By additivity of homology, for all~$j\in \IZ$, we have
		\[
			\dim_\IF H_j\bigl(\IF\otimes_\IZ (\ind_\Delta^\Gamma X)_\Lambda \bigr)= \frac{[\Gamma:\Lambda]}{[\Delta:\Lambda\cap \Delta]} \dim_\IF H_j(\IF\otimes_\IZ X_{\Lambda\cap \Delta})
		\]
		and the claim follows.
	\end{proof}
\end{lem}

\begin{prop}
\label{prop:bootstrappable H}
	Let~$\qm$ be the is-class of residually finite groups and let $\IF$ be a field.
	The family~$\sfH^\qm_*(\IF)$ is an equivariantly bootstrappable property of chain complexes over~$\IZ$.
	\begin{proof}
		First, for every residually finite group~$\Gamma$, the sequence~$\sfH_*^\Gamma(\IF)$ is a bootstrappable property of $\IZ\Gamma$-chain complexes.
		Indeed, the axioms~\ref{ax:pos} and~\ref{ax:susp} clearly hold and axiom~\ref{ax:cone} follows from Lemma~\ref{lem:Betti}~\ref{item:Betti cone}.
		Second, the family~$\sfH_*^\qm(\IF)$ satisfies axiom~\ref{ax:ind}, which follows from Lemma~\ref{lem:Betti}~\ref{item:Betti ind}.
	\end{proof} 
\end{prop}

We obtain the Bootstrapping Theorem~\ref{thm:bootstrapping} for~$\sfH_*(\IF)$.
The family~$\sfH^\qm_*(\IF)$ does not seem to satisfy axiom~\ref{ax:res} due to the non-transitivity of normal subgroups.

Instead of the dimension of homology, we can also consider the cardinality of the torsion subgroup of homology.
It is our primary interest to study the following classes:

\begin{defn}
\label{defn:T}
	Let~$\Gamma$ be a residually finite group and let~$n\in \IZ$.
	The class~$\sfT^\Gamma_n$ consists of all $\IZ\Gamma$-chain complexes~$X$ satisfying for all residual chains~$\Lambda_*$ in~$\Gamma$ and all~$j\le n$ that the torsion homology gradient
	\[
		\widehat{t}_j(X,\Lambda_*)\coloneqq
		\limsup_{i\to \infty}\frac{\log\tors H_j(X_{\Lambda_i})}{[\Gamma:\Lambda_i]}.
	\]
	is equal to~$0$.
	Here we use the convention $\log \infty \coloneqq \infty$.
\end{defn}
Then the class~$\sfT_n$ (Definition~\ref{defn:pre-boot}) consists of all residually finite groups~$\Gamma$ satisfying $\widehat{t}_j(\Gamma,\Lambda_*)=0$ for all residual chains~$\Lambda_*$ in~$\Gamma$ and all~$j\le n$ (Definition~\ref{defn:growth intro}).
Clearly, the class~$\sfT_0$ is the class of all residually finite groups.
It is easy to see that finite groups, free abelian groups, free groups, and surface groups lie in~$\sfT_\infty$.
Residually finite infinite amenable groups of type~$\sfFP_{n}$ lie in~$\sfT_{n-1}$~\cite[Corollary~2]{KKN17} (Corollary~\ref{cor:amenable T}).

For a group~$\Gamma$ of type~$\sfFP_n$, the Universal Coefficient Theorem shows that if~$\Gamma$ lies in~$\sfH_n(\IQ)\cap \sfT_n$, then~$\Gamma$ lies in~$\sfH_n(\IF)$ for every field~$\IF$. 

\begin{rem}\label{rem:T not bootstrappable}
	The sequence~$\sfT_*=(\sfT_n)_{n\in \IZ}$ of classes of groups is not a bootstrappable property of residually finite groups.
	
	Indeed, this follows from the computation of the torsion homology growth of right-angled Artin groups by Okun--Schreve~\cite{Okun-Schreve21}.
	Let $A_L$ be the right-angled Artin group associated to a finite flag simplicial complex~$L$.
	Then, for every~$j\in \IZ$ and every residual chain~$\Lambda_*$ in~$A_L$, we have
	\[
		\widehat{t}_j(A_L,\Lambda_*)= \log \tors H_{j-1}(L;\IZ).
	\]
	Now, let~$L$ be a flag triangulation of~$\IR P^2$ that admits a decomposition $L\cong L_1\cup_{L_0} L_2$ into full subcomplexes, where~$L_1$ is the M\"obius band, $L_0$ is the boundary of the M\"obius band, and~$L_2$ is a 2-dimensional disk.
	Then the group~$A_L$ splits as an amalgamated product $A_L\cong A_{L_1}\ast_{A_{L_0}} A_{L_2}$.
	By construction, we have $A_L\notin \sfT_2$ while $A_{L_i}\in \sfT_\infty$ for~$i\in \{0,1,2\}$ .
	Hence, by Corollary~\ref{cor:bootstrapping}~\ref{item:boot gog}, the sequence~$\sfT_*$ is not a bootstrappable property of groups.
\end{rem}

The failure of~$\sfT_*$ to be a bootstrappable property is due to the fact that torsion does not behave well with taking quotients, which causes complications for axiom~\ref{ax:cone}. 
In the next section, we introduce the sequence~$\sfCR^\Gamma_*$ which is designed to be (degreewise) contained in~$\sfT_{*-1}^\Gamma$ and simultaneously to be a bootstrappable property of $\IZ\Gamma$-chain complexes.

\section{The algebraic cheap rebuilding property}
\label{sec:CR}

We introduce a new bootstrappable property of residually finite groups, called the algebraic cheap rebuilding property (Definition~\ref{defn:CR}), which implies vanishing of torsion homology growth.
Our definitions and results are strongly inspired by the geometric cheap rebuilding property of \ABFG~\cite{ABFG21}.

\subsection{Quantitative homotopy retracts}
\label{sec:rebuilding}

We introduce the notion of rebuildings for chain complexes, which are homotopy retracts satisfying estimates on the ranks of chain modules and norms of chain maps.

We work over the ring~$\IZ$.
A free $\IZ$-module endowed with a $\IZ$-basis is called \emph{based free}.
The basis of a based free $\IZ$-module~$M$ induces an $\ell^1$-norm on~$M$.
For a map $f\colon M\to L$ between based free $\IZ$-modules, we denote by~$\|f\|$ the operator norm with respect to the $\ell^1$-norms on~$M$ and~$L$.
A $\IZ$-chain complex is \emph{based free} if every chain module is based free.

	A \emph{homotopy retract} of chain complexes $(X,X',\xi,\xi',\Xi)$ consists of chain complexes~$X$ and~$X'$, chain maps $\xi\colon X\to X'$ and $\xi'\colon X'\to X$, and a chain homotopy $\Xi\colon \id_X\simeq \xi'\circ \xi$.
	We will sometimes just write~$(X,X')$ instead of~$(X,X',\xi,\xi',\Xi)$, leaving the maps implicit.

\begin{defn}[Rebuilding]
\label{defn:rebuilding}
	Let~$n\in \IZ$ and let~$X$ and~$X'$ be based free $\IZ$-chain complexes such that $X_j$ and~$X'_j$ are finitely generated for all~$j\le n$.
	Let~$T,\kappa\in \IR_{\ge 1}$.
	We say that a homotopy retract $(X,X',\xi,\xi',\Xi)$ is
	\begin{itemize}
		\item an \emph{$n$-domination of~$X$ of quality~$(T,\kappa)$} if for all~$j\le n$
	\begin{align*}
		\rk_{\IZ}(X'_j) &\le \kappa T^{-1}\rk_\IZ(X_j);
	\end{align*}
	
		\item a \emph{weak $n$-rebuilding of~$X$ of quality~$(T,\kappa)$} if for all~$j\le n$
	\begin{align*}
		\rk_\IZ(X'_j) &\le \kappa T^{-1}\rk_\IZ(X_j);
		\\
		\|\partial^{X'}_j\|, \|\xi_j\| &\le 
		\exp(\kappa)T^\kappa;
	\end{align*}
	
		 \item an \emph{$n$-rebuilding of~$X$ of quality~$(T,\kappa)$} if for all~$j\le n$
	\begin{align*}
		\rk_\IZ(X'_j) &\le \kappa T^{-1}\rk_\IZ(X_j);
		\\
		\|\partial^{X'}_j\|, \|\xi_j\|, \|\xi'_j\|, \|\Xi_j\| &\le 
		\exp(\kappa)T^\kappa.
	\end{align*}
	\end{itemize}
\end{defn}

Clearly, an $n$-rebuilding is in particular a weak $n$-rebuilding of the same quality, and a weak $n$-rebuilding is in particular an $n$-domination of the same quality.
For~$T'\le T$, an $n$-domination of quality~$(T,\kappa)$ is in particular of quality~$(T',\kappa)$.
However, the analogous statement for (weak) $n$-rebuildings does not hold.
For (weak) $n$-rebuildings,
in the words of \ABFG~\cite{ABFG21}, the parameter~$T$ captures a tension between ``having small ranks'' and ``maintaining tame norms''.
One wants to decrease the ranks linearly in~$T$, while bounding the norms polynomially in~$T$.

\begin{rem}
\label{rem:explain weak rebuilding}
	Definition~\ref{defn:rebuilding} is designed in order to provide upper bounds on torsion in homology:
	Let~$(X,X',\xi,\xi',\Xi)$ be a weak $n$-rebuilding of quality~$(T,\kappa)$.
	Then~$H_j(X)$ is a retract of~$H_j(X')$ for all~$j\in \IZ$.
	Using Gabber's estimate~\eqref{eqn:Gabber}, for all~$j\le n-1$, we have
	\begin{equation}
	\label{eqn:explain weak rebuilding}
	\begin{aligned}
		\log\tors H_j(X)&\le \log\tors H_j(X')\le \rk_\IZ(X_j')\log_+\|\partial^{X'}_{j+1}\|
		\\
		&\le \kappa T^{-1}\rk_\IZ(X_j)\kappa(1+\log T). 
	\end{aligned}
	\end{equation}
	Now, suppose that~$X$ is a $\IZ$-chain complex and there exists a uniform~$\kappa\in \IR_{\ge 1}$ such that for all~$T\in \IR_{\ge 1}$, $X$ admits a weak $n$-rebuilding of quality~$(T,\kappa)$.
	Then the inequality~\eqref{eqn:explain weak rebuilding} shows that~$\log\tors H_j(X)=0$ because $T^{-1}\log T\to 0$ as $T\to \infty$.
	Later we will employ an asymptotic version of this argument (Lemma~\ref{lem:CWR implies T}).
	The other conditions in Definition~\ref{defn:rebuilding} on the norms of~$\xi$, $\xi'$, and~$\Xi$ ensure that the norm of~$\partial^{X'}$ stays controlled when taking mapping cones (Proposition~\ref{prop:rebuilding cones}).
\end{rem}

The following two examples are algebraic versions of topological homotopy equivalences between circles of different length. 

\begin{ex}
\label{ex:rebuilding circles short}
	For~$d\in \IN$, let $S^{[0,d]}$ be the chain complex with chain modules
	\[
		S^{[0,d]}_j=\begin{cases}
			\bigoplus_{i=0}^{d-1}\IZ\spann{v_i} &\text{for }j=0; \\
			\bigoplus_{i=0}^{d-1}\IZ\spann{e_i} &\text{for }j=1; \\
			0 &\text{otherwise};
		\end{cases}
	\]
	and differential $\partial_1(e_i)=v_{i+1}-v_i$ for~$i\in\{0,\ldots,d-1\}$ considered modulo~$d$.
	
	We construct an $n$-rebuilding~$(S^{[0,d]},S^{[0,1]})$ of quality~$(d,1)$.
	There is a homotopy retract $(S^{[0,d]},S^{[0,1]},\xi,\xi',\Xi)$,
	where the chain maps $\xi\colon S^{[0,d]}\to S^{[0,1]}$ and $\xi'\colon S^{[0,1]}\to S^{[0,d]}$ are given by
	\begin{align*}
		\xi_0(v_i) &= v_0 \quad \text{for all }i;
		\\
		\xi_1(e_i) &=\begin{cases}
			e_0 & \text{if }i=0; \\
			0 & \text{otherwise};
		\end{cases}
		\\
		\xi'_0(v_0) &= v_0;
		\\
		\xi'_1(e_0) &= \sum_{i=0}^{d-1}e_i;
	\end{align*}
	and the chain homotopy $\Xi\colon \id_{S^{[0,d]}}\simeq \xi'\circ \xi$ is given by
	\[
		\Xi_0(v_i)=\begin{cases}
			0 & \text{if }i=0; \\
			-e_i-\cdots-e_{d-1} &\text{if }i\in \{1,\ldots,d-1\}.
		\end{cases}
	\]
	For all~$j\in \IZ$, we have
	\[
		\rk_\IZ(S^{[0,1]}_j) \le d^{-1}\rk_\IZ(S^{[0,d]}_j);
	\]
	\[
		\|\partial^{S^{[0,1]}}_j\| \le 0; \quad
		\|\xi_j\| \le 1; \quad
		\|\xi'_j\| \le d; \quad
		\|\Xi_j\| \le d.
	\]	
	Hence, for all~$n\in \IZ$, the homotopy retract~$(S^{[0,d]},S^{[0,1]})$ is an $n$-rebuilding of quality~$(d,1)$.
	For~$T\le d$, the pair~$(S^{[0,d]},S^{[0,1]})$ provides a weak $n$-rebuilding of quality~$(T,1)$, but in general not an $n$-rebuilding of quality~$(T,1)$.
\end{ex}

The quality of the rebuilding in Example~\ref{ex:rebuilding circles short} is not good enough for later purposes (Example~\ref{ex:integers}).
Instead we will consider the following rebuilding, which is an algebraic reformulation of~\cite[Lemma~10.10]{ABFG21}. 

\begin{ex}
\label{ex:rebuilding circles}
	For~$d\in \IN$ and~$T\in \IR_{\ge 1}$ with~$T\le d$, we construct an $n$-rebuilding of~$S^{[0,d]}$ of quality~$(T,2)$.
	Choose a sequence of integers $0=a_0< a_1<\cdots<a_m=d$ with 
	\[
		T/2\le a_{k+1}-a_k\le T
	\]
	for all~$k\in \{0,\ldots, m-1\}$.
	There is a homotopy retract $(S^{[0,d]}, S^{[0,m]},\xi,\xi',\Xi)$, where the chain maps $\xi\colon S^{[0,d]}\to S^{[0,m]}$ and $\xi'\colon S^{[0,m]}\to S^{[0,d]}$ are given by
	\begin{align*}
		\xi_0(v_i) &=
			v_{k} \quad  \text{if }i\in \{a_{k-1}+1,\ldots, a_{k}\};
		\\
		\xi_1(e_i) &= \begin{cases}
			e_k &\text{if }i=a_k \text{ for some }k; \\
			0 &\text{otherwise};
		\end{cases}
		\\
		\xi'_0(v_i) &= v_{a_i}; 
		\\
		\xi'_1(e_i) &= e_{a_i}+\cdots+e_{a_{i+1}-1};
	\end{align*}
	and the chain homotopy $\Xi\colon \id_{S^{[0,d]}}\simeq \xi'\circ \xi$ is given by
	\[
		\Xi_0(v_i) = \begin{cases}
			0 & \text{if }i=a_k \text{ for some }k; \\
			-e_i-\cdots-e_{a_{k}-1} & \text{if }i\in\{a_{k-1}+1,\ldots,a_{k}-1\}.
		\end{cases}
	\]
	For all~$j\in \IZ$, we have
	\[
		\rk_\IZ(S^{[0,m]}_j) \le 2T^{-1}\rk_\IZ(S^{[0,d]}_j); 
	\]
	\[
		\|\partial^{S^{[0,m]}}_j\| \le 2; \quad
		\|\xi_j\| \le 1; \quad
		\|\xi'_j\| \le T; \quad
		\|\Xi_j\| \le T.
	\]
	Hence, for all~$n\in \IZ$, the homotopy retract $(S^{[0,d]},S^{[0,m]})$ is an $n$-rebuilding of quality~$(T,2)$.
\end{ex}

Directs sums of homotopy retracts are again homotopy retracts.
The same is true for rebuildings and we specify the quality of the resulting rebuilding.

\begin{lem}
\label{lem:rebuilding sums}
	Let $\bfX=(X,X',\xi,\xi',\Xi)$ and $\bfY=(Y,Y',\upsilon,\upsilon',\Upsilon)$ be homotopy retracts of $\IZ$-chain complexes.
	Consider the tuple
	\[
		\bfX\oplus \bfY\coloneqq (X\oplus Y,X'\oplus Y',\xi\oplus \upsilon, \xi'\oplus \upsilon', \Xi\oplus\Upsilon).
	\]
	Let~$n\in \IZ$, let~$T,\kappa_\bfX,\kappa_\bfY\in \IR_{\ge 1}$, and set~$\kappa\coloneqq \max\{\kappa_\bfX,\kappa_\bfY\}$.
	Then the following hold:
	\begin{enumerate}[label=\enum]
		\item The tuple~$\bfX\oplus \bfY$ is a homotopy retract;
		\item If~$\bfX$ is an $n$-domination of quality~$(T,\kappa_\bfX)$ and~$\bfY$ is an $n$-domination of quality~$(T,\kappa_\bfY)$, then~$\bfX\oplus \bfY$ is an $n$-domination of quality~$(T,\kappa)$; 
		\item If~$\bfX$ is a weak $n$-rebuilding of quality~$(T,\kappa_\bfX)$ and~$\bfY$ is a weak $n$-rebuilding of quality~$(T,\kappa_\bfY)$, then~$\bfX\oplus \bfY$ is a weak $n$-rebuilding of quality~$(T,\kappa)$;
		\item If~$\bfX$ is an $n$-rebuilding of quality~$(T,\kappa_\bfX)$ and~$\bfY$ is an $n$-rebuilding of quality~$(T,\kappa_\bfY)$, then~$\bfX\oplus \bfY$ is an $n$-rebuilding of quality~$(T,\kappa)$.
	\end{enumerate}
	\begin{proof}
	Part~(i) is clear.
	
	(ii) Suppose~$\bfX$ and~$\bfY$ are $n$-dominations of quality~$(T,\kappa_\bfX)$ and~$(T,\kappa_\bfY)$, respectively. 
	For all~$j\le n$, we have
	\[
		\rk_\IZ(X'_j\oplus Y'_j) \le \kappa_\bfX T^{-1}\rk_\IZ(X_j) + \kappa_\bfY T^{-1}\rk_\IZ(Y_j) \le \kappa T^{-1}\rk_\IZ(X_j\oplus Y_j).
	\]
	
	(iii) Suppose~$\bfX$ and~$\bfY$ are weak $n$-rebuildings of quality~$(T,\kappa_\bfX)$ and~$(T,\kappa_\bfY)$, respectively. 
	Additionally to part~(ii), for all~$j\le n$, we have
	\[
		\|\partial^{X'}_j\oplus \partial^{Y'}_j\| \le \max\bigl\{\|\partial^{X'}_j\|,\|\partial^{Y'}_j\|\bigr\} \le \max\bigl\{\exp(\kappa_\bfX)T^{\kappa_\bfX}, 
	\exp(\kappa_\bfY)T^{\kappa_\bfY}\bigr\} \le \exp(\kappa)T^\kappa
	\]
	and similarly $\|\xi_j\oplus \upsilon_j\| \le \exp(\kappa)T^\kappa$.
	
	(iv) Suppose~$\bfX$ and~$\bfY$ are $n$-rebuildings of quality~$(T,\kappa_\bfX)$ and~$(T,\kappa_\bfY)$, respectively.
	Additionally to part~(iii), for all~$j\le n$, we have
	\[
		\|\xi'_j\oplus \upsilon'_j\|, \|\Xi_j\oplus\Upsilon_j\| \le \exp(\kappa)T^\kappa
	\]
	by a similar computation.
	\end{proof}
\end{lem}

Mapping cones of homotopy retracts are again homotopy retracts.
The same is true for rebuildings, which is our key stability result for rebuildings.
A homotopy retract~$(X,X',\xi,\xi',\Xi)$ can be viewed as a homotopy commutative square of the form
\[
	\begin{tikzcd}
	X\ar[r, "\xi"]\ar[d, "\id_X" left] & X'\ar[d, "\xi'"]\\
	X\ar[r, "\id_X" below] & X\arrow[from=1-1, to=2-2, phantom, "\text{\Large$\circlearrowleft$}_{\Xi}" description]
	\end{tikzcd}
\]

\begin{prop}[Rebuilding of mapping cones]
\label{prop:rebuilding cones}
	Let~$\bfX=(X,X',\xi,\xi',\Xi)$ and $\bfY=(Y,Y',\upsilon,\upsilon',\Upsilon)$ be homotopy retracts of $\IZ$-chain complexes.
	Let~$f\colon X\to Y$ be a chain map and define the chain map $f'\coloneqq \upsilon\circ f\circ \xi'\colon X'\to Y'$.
	Consider the tuple 
	\[
		\bfC(f)\coloneqq \bigl(\Cone(f),\Cone(f'),(\xi,\upsilon;-\upsilon\circ f\circ \Xi),(\xi',\upsilon';\Upsilon\circ f\circ \xi'),\Psi\bigr),
	\]
	where $\Psi\colon \Cone(f)_*\to \Cone(f)_{*+1}$ is defined as
	\[
		\Psi_j(x,y)\coloneqq \bigl(-\Xi_{j-1}(x),\Upsilon_j(y)+\Upsilon_j\circ f_j\circ \Xi_{j-1}(x)\bigr).
	\]
	Let~$n\in \IN$, let~$T,\kappa_\bfX,\kappa_\bfY\in \IR_{\ge 1}$, and set
	\[
		\kappa\coloneqq \kappa_\bfX+\kappa_\bfY+\log 3+\max\bigl\{\log_+\|f_j\|\bigm\vert j\le n\bigr\},
	\]
	where~$\log_+\coloneqq \max\{\log,0\}$.
	Then the following hold:
	\begin{enumerate}[label=\enum]
	\item The tuple $\bfC(f)$ is a homotopy retract;
	\item If~$\bfX$ is an $(n-1)$-domination of quality~$(T,\kappa_\bfX)$ and~$\bfY$ is an $n$-domination of quality~$(T,\kappa_\bfY)$, then $\bfC(f)$ is an $n$-domination of quality~$(T,\kappa)$;
	\item\label{item:warning} If~$\bfX$ is an $(n-1)$-rebuilding of quality~$(T,\kappa_\bfX)$ and~$\bfY$ is a weak $n$-rebuilding of quality~$(T,\kappa_\bfY)$, then $\bfC(f)$ is a weak $n$-rebuilding of quality~$(T,\kappa)$;
	\item If~$\bfX$ is an $(n-1)$-rebuilding of quality~$(T,\kappa_\bfX)$ and~$\bfY$ is an $n$-rebuilding of quality~$(T,\kappa_\bfY)$, then $\bfC(f)$ is an $n$-rebuilding of quality~$(T,\kappa)$.
	\end{enumerate}
	\textup{We point out that in part~\ref{item:warning}, $\bfX$ is assumed to be an $(n-1)$-rebuilding, and not only a weak $(n-1)$-rebuilding.
	This will have important ramifications later (Proposition~\ref{prop:bootstrappable CR}).}
	\begin{proof}
		(i) The following cube is homotopy commutative in the sense of~\eqref{eqn:cube}:
		\begin{equation}
		\label{eqn:cube retracts}
		\begin{tikzcd}[row sep=large, column sep=huge]
		X'\ar[rrr, "f'"]\ar[ddd, "\xi'" left] 
		& \arrow[from=1-2, to=2-3, phantom, "\text{\Large$\circlearrowleft$}_{-\upsilon\circ f\circ \Xi}" description]
		&& Y'\ar[ddd, "\upsilon'"]
		&& \arrow[from=1-6, to=4-4, phantom, "\text{\Large$\circlearrowleft$}_{\Upsilon\circ f\circ \xi'}" description]
		\\
		\arrow[from=2-1, to=3-2, phantom, "\text{\Large$\circlearrowleft$}_{\Xi}" description]
		& X\ar[r, "f"]\ar[lu, "\xi\ " below]\ar[d, "\id_X" left] 
		\arrow[from=2-2, to=3-3, phantom, "\text{\Large$\circlearrowleft$}_{0}" description]
		& Y\ar[ru, "\ \upsilon" below]\ar[d, "\id_Y"]
		\\
		& X\ar[r, "f" below]\ar[ld, "\id_X\ \ \ " above] & Y\ar[dr, "\ \ \ \id_Y" above]
		& \arrow[from=3-4, to=2-3, phantom, "\text{\Large$\circlearrowleft$}_{\Upsilon}" description]
		\\
		X\ar[rrr, "f" below] 
		& \arrow[from=4-2, to=3-3, phantom, "\text{\Large$\circlearrowleft$}_{0}" description]
		&& Y
		\end{tikzcd}
		\end{equation}
		with $\Phi_*\coloneqq -\Upsilon_{*+1}\circ f_{*+1}\circ \Xi_*\colon X_*\to Y_{*+2}$.
		Indeed, one easily checks that the upper and outer squares are homotopy commutative.
		The other four squares are clearly homotopy commutative.
		We check that~$\Phi$ is a filler of the cube in the sense of equation~\eqref{eqn:filler}: 
		\begin{align*}
			 & \id_{Y_{j+1}}\circ 0 -\Upsilon_j \circ f_j\circ \xi'_j\circ \xi_j + \Upsilon_j\circ f_j - f_{j+1}\circ \Xi_j +0\circ \id_{X_j} + \upsilon'_{j+1}\circ \upsilon_{j+1}\circ f_{j+1}\circ \Xi_j 
			\\
			&= \Upsilon_j\circ f_j - \Upsilon_j \circ f_j\circ \xi'_j\circ \xi_j  - f_{j+1}\circ \Xi_j + \upsilon'_{j+1}\circ \upsilon_{j+1}\circ f_{j+1}\circ \Xi_j 
			\\
			&= \Upsilon_j\circ f_j\circ (\id_{X_j}-\xi'_j\circ \xi_j) - (\id_{Y_{j+1}}-\upsilon'_{j+1}\circ \upsilon_{j+1})\circ f_{j+1}\circ \Xi_j
			\\
			&= \Upsilon_j\circ f_j\circ (\partial^X_{j+1}\circ \Xi_j+\Xi_{j-1}\circ \partial^X_j) - (\partial^Y_{j+2}\circ \Upsilon_{j+1}+\Upsilon_j\circ \partial^Y_{j+1})\circ f_{j+1}\circ \Xi_j
			\\
			&= \Upsilon_j\circ f_j\circ \Xi_{j-1}\circ \partial^X_j - \partial^Y_{j+2}\circ \Upsilon_{j+1}\circ f_{j+1}\circ \Xi_j
			\\
			&= \partial^Y_{j+2}\circ \Phi_j - \Phi_{j-1}\circ \partial^X_j.
		\end{align*}
		Here the third equality uses that~$\Xi$ and~$\Upsilon$ are chain homotopies and the fourth equality uses that~$f$ is a chain map.
		The other steps consist merely of rearranging the terms.
		
		By Lemma~\ref{lem:cube}, the homotopy commutative cube~\eqref{eqn:cube retracts} induces a homotopy commutative square of mapping cones
		\[
		\begin{tikzcd}[column sep=huge]
		\Cone(f)\ar[r,"{(\xi,\upsilon;-\upsilon\circ f\circ \Xi)}"]\ar[d,"{(\id_X,\id_Y;0)}" left]
		& \Cone(f')\ar[d, "{(\xi',\upsilon';\Upsilon\circ f\circ \xi')}"]
		\\
		\Cone(f)\ar[r, "{(\id_X,\id_Y;0)}" below]
		& \Cone(f)
		\arrow[from=2-2, to=1-1, phantom, "\text{\Large$\circlearrowleft$}_{\Psi}" description]
		\end{tikzcd}
		\]
		where 
		\[\Psi_j(x,y) = \bigl(-\Xi_{j-1}(x), \Upsilon_j(y) - \Phi_{j-1}(x)\bigr). \]
		Hence~$\bfC(f)$ is a homotopy retract of chain complexes.
		
		(ii) Suppose that~$\bfX$ is an $(n-1)$-domination of quality~$(T,\kappa_\bfX)$ and~$\bfY$ is an $n$-domination of quality~$(T,\kappa_\bfY)$. 
		For all~$j\le n$, we have
		\begin{align*}
			\rk_\IZ\bigl(\Cone(f')_j\bigr)
			&= \rk_\IZ(X'_{j-1})+\rk_\IZ(Y'_j)
			\\
			&\le \kappa_\bfX T^{-1}\rk_\IZ(X_{j-1}) + \kappa_\bfY T^{-1}\rk_\IZ(Y_j)
			\\
			&\le \kappa T^{-1}\bigl(\rk_\IZ(X_{j-1})+\rk_\IZ(Y_j)\bigr)
			\\
			&= \kappa T^{-1}\rk_\IZ\bigl(\Cone(f)_j\bigr).
		\end{align*}
		
		(iii) Suppose that~$\bfX$ is an $(n-1)$-rebuilding of quality~$(T,\kappa_\bfX)$ and~$\bfY$ is an $n$-weak rebuilding of quality~$(T,\kappa_\bfY)$.
		Additionally to part~(ii), for all~$j\le n$, we have
		\begin{align*}
			\|\partial^{\Cone(f')}_j\|
			&\le \|\partial^{X'}_{j-1}\|+ \|\partial^{Y'}_j\|+ \|\upsilon_{j-1}\|\cdot \|f_{j-1}\|\cdot \|\xi'_{j-1}\|
			\\
			&\le \exp(\kappa_\bfX)T^{\kappa_\bfX} + \exp(\kappa_\bfY)T^{\kappa_\bfY} + 
			\\
			& \qquad \exp(\kappa_\bfY+\log_+\|f_{j-1}\|+\kappa_\bfX)T^{\kappa_\bfY+\log_+\|f_{j-1}\|+\kappa_\bfX}
			\\
			&\le 3\exp(\kappa_\bfX+\kappa_\bfY+\log_+\|f_{j-1}\|)T^{\kappa_\bfX+\kappa_\bfY+\log_+\|f_{j-1}\|}
			\\
			&\le \exp(\kappa_\bfX+\kappa_\bfY+\log_+\|f_{j-1}\|+\log 3)T^{\kappa_\bfX+\kappa_\bfY+\log_+\|f_{j-1}\|+\log 3}
			\\
			&\le \exp(\kappa)T^\kappa
		\end{align*}
		and similarly
		\[
		\|(\xi,\upsilon;-\upsilon\circ f\circ \Xi)_j\| 
			\le \|\xi_{j-1}\| + \|\upsilon_j\|+\|\upsilon_j\|\cdot \|f_j\|\cdot \|\Xi_{j-1}\|
			\le
			\exp(\kappa)T^\kappa.
			\]
           
           (iv) Suppose that~$\bfX$ is an $(n-1)$-rebuilding of quality~$(T,\kappa_\bfX)$ and~$\bfY$ is an $n$-rebuilding of quality~$(T,\kappa_\bfY)$.
           Additionally to part~(iii), for all~$j\le n$, we have
		\begin{align*}
			\|(\xi',\upsilon';\Upsilon\circ f\circ \xi')_j\|
			&\le \|\xi'_{j-1}\|+ \|\upsilon'_j\|+ \|\Upsilon_{j-1}\|\cdot \|f_{j-1}\|\cdot \|\xi'_{j-1}\|
			\le
			\exp(\kappa)T^\kappa;
			\\
			\|\Psi_j\|
			&\le \|\Xi_{j-1}\|+ \|\Upsilon_j\|+ \|\Upsilon_j\|\cdot \|f_j\|\cdot \|\Xi_{j-1}\|
			\le
			\exp(\kappa)T^\kappa;		
			\end{align*}
			by similar computations.
	\end{proof}
\end{prop}

While direct sums can be viewed as mapping cones of the zero map, we formulated Lemma~\ref{lem:rebuilding sums} separately from Proposition~\ref{prop:rebuilding cones} for the optimality of the constant~$\kappa$, which is used later (Proposition~\ref{prop:bootstrappable CR}).
	
The composition of homotopy retracts is again a homotopy retract.
The same is true for dominations and weak rebuildings.
We will not require the corresponding result of rebuildings, though a similar statement (involving a degree shift) also holds.
Lemma~\ref{lem:composition} is an algebraic version of~\cite[Lemma~6.3]{ABFG21}.
\begin{lem}\label{lem:composition}
	Let~$\bfX=(X,X',\xi,\xi',\Xi)$ and~$\bfY=(X',X'',\upsilon,\upsilon',\Upsilon)$ be homotopy retracts of $\IZ$-chain complexes.
	Consider the tuple
	\[
		\bfY\circ \bfX\coloneqq (X,X'',\upsilon\circ \xi,\xi'\circ \upsilon',\Xi+\xi'\circ \Upsilon\circ \xi).
	\]
	Let~$n\in \IN$, let~$T,S,\kappa_\bfX,\kappa_\bfY\in \IR_{\ge 1}$, and set~$\kappa\coloneqq 2\kappa_\bfY\kappa_\bfX$.
	Then the following hold:
	\begin{enumerate}[label=\enum]
		\item The tuple~$\bfY\circ \bfX$ is a homotopy retract;
		\item If~$\bfX$ is an $n$-domination of quality~$(T,\kappa_\bfX)$ and~$\bfY$ is an $n$-domination of quality~$(S,\kappa_\bfY)$, then~$\bfY\circ\bfX$ is an $n$-domination of quality~$(ST,\kappa)$;
		\item\label{item:composition weak rebuilding} If~$\bfX$ is a weak $n$-rebuilding of quality~$(T,\kappa_\bfX)$ and~$\bfY$ is a weak $n$-rebuilding of quality~$(S,\kappa_\bfY)$, then~$\bfY\circ\bfX$ is a weak $n$-rebuilding of quality~$(ST,\kappa)$.
	\end{enumerate}
	\begin{proof}
		(i) A straight-forward calculation shows that~$\Xi+\xi'\circ \Upsilon\circ \xi$ indeed provides a chain homotopy between the chain maps~$\id_X$ and~$\xi'\circ \upsilon'\circ \upsilon\circ \xi$.
		
		(ii) Suppose~$\bfX$ and~$\bfY$ are $n$-dominations of quality~$(T,\kappa_\bfX)$ and~$(T,\kappa_\bfY)$, respectively.
		For all~$j\le n$, we have
		\[
			\rk_\IZ(X''_j)\le \kappa_\bfY S^{-1}\rk_\IZ(X'_j)\le \kappa_\bfY S^{-1}\kappa_\bfX T^{-1}\rk_\IZ(X_j) \le \kappa (ST)^{-1}\rk_\IZ(X_j).
		\]
		
		(iii) Suppose~$\bfX$ and~$\bfY$ are weak $n$-rebuildings of quality~$(T,\kappa_\bfX)$ and~$(T,\kappa_\bfY)$, respectively.
		Additionally to part~(ii), for all~$j\le n$, we have
		\begin{align*}
			\|\partial^{X''}_j\|\ &\le \exp(\kappa_\bfY) T^{\kappa_\bfY} \le \exp(\kappa)T^{\kappa};
			\\
			\|\upsilon_j\circ \xi_j\| &\le \|\upsilon_j\|\cdot \|\xi_j\|\le \exp(\kappa_\bfY)S^{\kappa_{\bfY}}\exp(\kappa_\bfX)T^{\kappa_{\bfX}} 
			 \le \exp(\kappa)(ST)^{\kappa}.
		\end{align*}
		Hence~$\bfY\circ \bfX$ is a weak $n$-rebuilding of quality~$(ST,\kappa)$.
	\end{proof}
\end{lem}
	
\begin{rem}[Geometric rebuilding]
\label{rem:rebuilding}
	Definition~\ref{defn:rebuilding} of $n$-rebuildings for chain complexes is an algebraic version of \ABFG's geometric definition of $n$-rebuildings for CW-complexes~\cite[Definition~1 and Definition~2]{ABFG21}.
	 However, we point out the main differences:
	 We require a homotopy retract in all degrees, as opposed to a truncated homotopy equivalence.
	 We work with the $\ell^1$-norm on based free $\IZ$-modules, instead of the $\ell^2$-norm, because it simplifies some calculations.
	 We also ask for control on the norm of the homotopy~$\Xi_j$ in degrees~$j\le n$, and not only in degrees~$j\le n-1$, because it is needed in the proof of Proposition~\ref{prop:rebuilding cones}.
	 Due to these differences, there is no obvious implication between the algebraic and geometric notions.
	 We chose to formulate Definition~\ref{defn:rebuilding} as above because it simplifies the proofs and still covers the main Example~\ref{ex:rebuilding circles}.
\end{rem}

\subsection{The algebraic cheap (weak) rebuilding property}

We work over the ring~$R=\IZ$.
Recall that for a subgroup~$\Lambda$ of~$\Gamma$ and a $\IZ\Gamma$-chain complex~$X$, we write~$X_\Lambda\coloneqq \IZ\otimes_{\IZ\Lambda} \res^\Gamma_\Lambda X$ for the $\IZ$-chain complex of $\Lambda$-coinvariants.
We consider only residually finite groups.

\begin{defn}
\label{defn:CR}
	Let~$\Gamma$ be a residually finite group and let~$n\in \IZ$.
	The class~$\sfCR^\Gamma_n$ (resp.\ $\sfCWR^\Gamma_n$, $\sfCD^\Gamma_n$) consists of all based free $\IZ\Gamma$-chain complexes~$X$ 
	such that for all~$j\le n$ the $\IZ\Gamma$-module~$X_j$ is finitely generated and~$X$ satisfies the following:
	there exists~$\kappa\in \IR_{\ge 1}$ such that
	for all~$T\in \IR_{\ge 1}$ and all residual chains~$\Lambda_*$ in~$\Gamma$,
	there exists~$i_0\in \IN$ such that for all~$i\ge i_0$, the $\IZ$-chain complex~$X_{\Lambda_i}$ admits an $n$-rebuilding (resp.\ weak $n$-rebuilding, $n$-domination) of quality~$(T,\kappa)$ (Definition~\ref{defn:rebuilding}).
	
	In symbols, 
	\[
		\exists_{\kappa\ge 1} \ 
		\forall_{T\ge 1} \ 
		\forall_{\Lambda_* \text{ res.\ chain in }\Gamma} \
		\exists_{i_0\in \IN} \ 
		\forall_{i\ge i_0} \
		\exists_{\text{$n$-rebuilding of $X_{\Lambda_i}$ of quality } (T,\kappa)}
	\]	
	We define the class~$\sfCR_n$ of residually finite groups~$\Gamma$ for which there exists a projective $\IZ\Gamma$-resolution of~$\IZ$ lying in~$\sfCR^\Gamma_n$.
	We also define the class of groups $\sfCR_\infty\coloneqq \bigcap_{n\in \IZ}\sfCR_n$.
	Similarly, we define~$\sfCWR_n,\sfCWR_\infty,\sfCD_n,\sfCD_\infty$.
	If the group~$\Gamma$ lies in~$\sfCR_n$ (resp.\ $\sfCWR_n$, $\sfCD_n$), we say that~$\Gamma$ satisfies the \emph{algebraic cheap $n$-rebuilding property} (resp.\ \emph{algebraic cheap weak $n$-rebuilding property}, \emph{algebraic cheap $n$-domination property}).
\end{defn}

Recall from Remark~\ref{rem:explain weak rebuilding} that (weak) rebuildings yield upper bounds and vanishing for torsion in homology.
The algebraic cheap (weak) rebuilding property is designed to yield \emph{asymptotic} vanishing for torsion in homology along every residual chain (Lemma~\ref{lem:CWR implies T}).

Clearly, we have inclusions of classes $\sfCR_n^\Gamma\subset  \sfCWR_n^\Gamma\subset \sfCD_n^\Gamma$ for every~$n\in \IZ$.

\begin{lem}
\label{lem:deg0 infinite}
	We have the following equalities of classes:
	\[
		\sfCR_0=\sfCWR_0=\sfCD_0= \textup{class of all residually finite infinite groups}.
	\]
	\begin{proof}
		The inclusions $\sfCR_0\subset \sfCWR_0\subset \sfCD_0$ are clear.
		
		Supposing that~$\Gamma$ lies in~$\sfCD_0$, we show that~$\Gamma$ is infinite.
		Let~$X$ be a free $\IZ\Gamma$-resolution of~$\IZ$ lying in~$\sfCD^\Gamma_0$ witnessed by the constant~$\kappa$.
		For every~$T\in \IR_{\ge 1}$ and every residual chain~$\Lambda_*$ in~$\Gamma$, there exists~$i_0\in \IN$ such that for all~$i\ge i_0$, there exists a 0-domination~$(X_{\Lambda_i},X_{\Lambda_i}')$ of quality~$(T,\kappa)$.
		That is,
		\[
			\rk_\IZ \bigl((X_{\Lambda_i}')_0\bigr)\le \kappa T^{-1}\rk_\IZ\bigl((X_{\Lambda_i})_0\bigr) = \kappa T^{-1}[\Gamma:\Lambda_i]\rk_{\IZ\Gamma}(X_0). 
		\]
		Since~$H_0(X_{\Lambda_i}')$ contains~$H_0(X_{\Lambda_i})\cong \IZ$ as a retract, we have~$\rk_\IZ((X_{\Lambda_i}')_0)\ge 1$.
		Hence
		\[
			\kappa^{-1}T\rk_{\IZ\Gamma}(X_0)^{-1}\le [\Gamma:\Lambda_i].
		\]
		Since~$T\in \IR_{\ge 1}$ was arbitrary, we must have $\lim_{i\to \infty}[\Gamma:\Lambda_i]=\infty$. 
		Hence~$\Gamma$ is infinite.
		
		Supposing that~$\Gamma$ is residually finite infinite, we show that~$\Gamma$ lies in~$\sfCR_0$.
		Fix a based free $\IZ\Gamma$-resolution~$X$ of~$\IZ$ with $X_0=\IZ\Gamma$, $X_1=\bigoplus_\Gamma \IZ\Gamma$, and $\partial^X_1(e_\gamma)=\gamma-1$. 
		Here~$e_\gamma\in X_1$ denotes the $\IZ\Gamma$-basis element corresponding to~$\gamma\in \Gamma$.
		We will show for every finite index subgroup~$\Lambda$ of~$\Gamma$, that~$X_\Lambda$ admits a $0$-rebuilding of quality~$(T,1)$ for all~$T\le [\Gamma:\Lambda]$.
		Then~$X$ lies in~$\sfCR^\Gamma_0$ witnessed by the constant~$\kappa=1$ because for every~$T\in \IR_{\ge 1}$ and every residual chain~$\Lambda_*$ in~$\Gamma$, there exists~$i_0\in \IN$ such that for all~$i\ge i_0$, $[\Gamma:\Lambda_i]\ge T$ since~$\Gamma$ is infinite.
				
		Now, let~$\Lambda$ be a finite index subgroup of~$\Gamma$.		
		Choose a set~$S$ of right-coset representatives for~$\Lambda\backslash \Gamma$ with~$1_\Gamma\in S$.
		Consider the isomorphism of $\IZ\Lambda$-modules
		\begin{align*}
			\res^\Gamma_\Lambda \IZ\Gamma &\cong \bigoplus_S \IZ\Lambda
			\\
			s &\mapsfrom e_s
			\\
			\gamma &\mapsto \gamma t(\gamma)^{-1}e_{t(\gamma)}
		\end{align*}
		where~$e_s$ is the $\IZ\Lambda$-basis element corresponding to~$s\in S$ and~$t(\gamma)\in S$ is such that~$\Lambda\gamma=\Lambda t(\gamma)$.
		Under this isomorphism, the based free $\IZ\Lambda$-resolution~$\res^\Gamma_\Lambda X$ of~$\IZ$ is given in degrees~$1$ and~$0$ by
		\[
			\cdots\to \bigoplus_{\Gamma\times S} \IZ\Lambda\xrightarrow{\res^\Gamma_\Lambda \partial^X_1} \bigoplus_S \IZ\Lambda,	
		\]
		where
		\[
			\res^\Gamma_\Lambda \partial^X_1(e_{(\gamma,s)}) = s\gamma t(s\gamma)^{-1}e_{t(s\gamma)}-e_s.
		\]
		Let~$Y$ be a based free $\IZ\Lambda$-resolution of~$\IZ$ with~$Y_0=\IZ\Lambda$.
		Then there exist mutually $\IZ\Lambda$-chain homotopy inverse $\IZ\Lambda$-chain maps $\xi\colon \res^\Gamma_\Lambda X\to Y$ with~$\xi_0(e_s)=1$ and $\xi'\colon Y\to \res^\Gamma_\Lambda X$ with~$\xi'_0(1)=e_{1_\Gamma}$.
		The $\IZ\Lambda$-chain homotopy $\Xi\colon \id_{\res^\Gamma_\Lambda X}\simeq \xi'\circ \xi$ is given in degree~$0$ by $\Xi_0(e_s)=-e_{(s^{-1},s)}$.
		By construction, we have 
		\begin{align*}
			\rk_{\IZ\Lambda}(Y_0) &\le  [\Gamma:\Lambda]^{-1} \rk_{\IZ\Lambda} (\res^\Gamma_\Lambda X_0);
			\\
			\|\partial^Y_0\|, \|\xi_0\|, \|\xi'_0\|, \|\Xi_0\| & \le 1.
		\end{align*}
		Hence by the following Remark~\ref{rem:equivariant rebuilding}, the $\IZ\Lambda$-homotopy equivalence between~$\res^\Gamma_\Lambda X$ and~$Y$ descends to a $0$-rebuilding~$(X_\Lambda,Y_\Lambda)$ of quality~$(T,1)$ for all~$T\le [\Gamma:\Lambda]$.
	\end{proof}
\end{lem}

\begin{rem}[Equivariant rebuilding]
\label{rem:equivariant rebuilding}
	Let~$\Gamma$ be a group and let~$(X,X',\xi,\xi',\Xi)$ be a $\IZ\Gamma$-homotopy retract of based free $\IZ\Gamma$-chain complexes lying in~$\sfFG_n^\Gamma(\IZ)$.
	There exist~$T,\kappa\in \IR_{\ge 1}$ such that for all~$j\le n$, we have
	\begin{align*}
		\rk_{\IZ\Gamma}(X_j')
		&\le \kappa T^{-1}\rk_{\IZ\Gamma}(X_j);
		\\
		\|\partial^{X'}_j\|, \|\xi_j\|, \|\xi_j'\|, \|\Xi_j\|
		&\le \exp(\kappa)T^\kappa.
	\end{align*}
	Let~$\Lambda$ be a finite index subgroup of~$\Gamma$.
	Then the $\IZ\Gamma$-homotopy retract~$(X,X')$ descends to an $n$-rebuilding~$(X_{\Lambda},X_{\Lambda}')$ of quality~$(T,\kappa)$.
	
	Indeed, for all~$j\le n$, we have
	\[
		\rk_\IZ\bigl((X_{\Lambda}')_j\bigr)
		=[\Gamma:\Lambda]\rk_{\IZ\Gamma}(X_j')\le [\Gamma:\Lambda]\kappa T^{-1} \rk_{\IZ\Gamma}(X_j)=\kappa T^{-1} \rk_\IZ\bigl((X_{\Lambda})_j\bigr).	
	\]
	Since the functor~$(-)_{\Lambda}$ does not increase the norm of maps by Lemma~\ref{lem:coinvariants}~\ref{item:coinvariants map},
	for every~$f\in\{\partial^{X'},\xi,\xi',\Xi\}$ and all~$j\le n$, we have
	\[
		\|(f_{\Lambda})_j\| \le \|f_j\| \le \exp(\kappa) T^\kappa.
	\]
	Similar statements hold for weak rebuildings and dominations.
\end{rem}

The main example of a group for which we can show directly from the definition that it lies in~$\sfCR_\infty$ is the group of integers.

\begin{ex}
\label{ex:integers}
	The group~$\IZ$ lies in~$\sfCR_\infty$.
	
	Indeed, let~$t\in \IZ$ be a generator.
	Consider the projective $\IZ[\spann{t}]$-resolution~$X$ of~$\IZ$
	\[
		0\to X_1 = \IZ[\spann{t}]\xrightarrow{\partial_1} X_ 0=\IZ[\spann{t}] \to \Z \to 0,
	\]
	where $\partial_1(t)=t-1$.
	Set~$\kappa\coloneqq 2$.
	Let~$T\in \IR_{\ge 1}$ and let~$\Lambda_*$ be a residual chain in~$\Gamma$.
	Take~$i_0\in \IN$ such that~$[\Gamma:\Lambda_{i_0}]\ge T$. 
	For all~$i\ge i_0$, the subgroup~$\Lambda_i$ of~$\spann{t}$ is of the form~$\spann{t^{d_i}}$, where $d_i\coloneqq [\Gamma:\Lambda_i]$.
	Then the $\IZ$-chain complex~$X_{\Lambda_i}$ can be identified with~$S^{[0,d_i]}$ from Example~\ref{ex:rebuilding circles}, which admits an $n$-rebuilding of quality~$(T,2)$.
\end{ex}

We will show in Section~\ref{sec:amenable} that residually finite infinite amenable groups of type~$\sfFP_\infty$ lie in~$\sfCWR_\infty$.
We do not know if all these groups lie in the subclass~$\sfCR_\infty$.
However, it will follow from the Bootstrapping Theorem~\ref{thm:bootstrapping} that residually finite infinite \emph{elementary} amenable groups of type~$\sfFP_\infty$ lie in~$\sfCR_\infty$.

The sequences~$\sfCD_*^\Gamma$ and~$\sfCWR_*^\Gamma$, which are defined on the level of chain complexes, are (degreewise) contained in the sequences~$\sfH_*^\Gamma(\IF)$ and~$\sfT_{*-1}^\Gamma$ (Definition~\ref{defn:H} and Definition~\ref{defn:T}), respectively, which are defined on the level of homology.
Our arguments are algebraic versions of~\cite[Theorem~10.20]{ABFG21}.

\begin{lem}
\label{lem:CWR implies T}
	Let~$\Gamma$ be a residually finite group, let~$n\in \IZ$, and let~$X$ be a $\IZ\Gamma$-chain complex.
	The following hold:
	\begin{enumerate}[label=\enum]
		\item If~$X\in \sfCD_n^\Gamma$, then $X\in \sfH_n^\Gamma(\IF)$ for every field~$\IF$;
		\item If~$X\in \sfCWR_n^\Gamma$, then~$X\in \sfT_{n-1}^\Gamma$;
		\item If~$X\in\sfCWR_n^\Gamma$ and~$X\in \sfFG_{n+1}^\Gamma(\IZ)$, then~$X\in \sfT_n^\Gamma$;
		\item\label{item:additional degree}
		If the group~$\Gamma$ lies in~$\sfCWR_n$ and is of type~$\sfFP_{n+1}$, then~$\Gamma\in \sfT_n$.
	\end{enumerate}
	\begin{proof}	
		(i) Suppose that~$X$ lies in~$\sfCD_n^\Gamma$ witnessed by the constant~$\kappa$.
		Let~$\Lambda_*$ be a residual chain in~$\Gamma$ and let~$j\le n$.
		For all~$T\in \IR_{\ge 1}$, there exists~$i_0\in \IN$ such that for all~$i\ge i_0$, the $\IZ$-chain complex~$X_{\Lambda_i}$ admits an $n$-domination~$(X_{\Lambda_i},X'_{\Lambda_i})$ of quality~$(T,\kappa)$.
		That is, for all~$j\le n$, we have
		\[
			\rk_\IZ\bigl((X'_{\Lambda_i})_j\bigr)
			\le \kappa T^{-1}\rk_\IZ\bigl((X_{\Lambda_i})_j\bigr)
			\le \kappa T^{-1}[\Gamma:\Lambda_i]\rk_{\IZ\Gamma}(X_j).
		\]
		Since~$H_j(\IF\otimes_\IZ X_{\Lambda_i})$ is a retract of~$H_j(\IF\otimes_\IZ X'_{\Lambda_i})$, we have
		\[
			\dim_\IF H_j(\IF\otimes_\IZ X_{\Lambda_i})
			\le \dim_\IF H_j(\IF\otimes_\IZ X'_{\Lambda_i}) 	
			\le \rk_\IZ\bigl((X'_{\Lambda_i})_j\bigr).
		\]
		Together, for all~$j\le n$, we have
		\[
			\widehat{b}_j(X,\Lambda_*;\IF)\le \limsup_{i\to \infty}\frac{\kappa T^{-1}[\Gamma:\Lambda_i]\rk_{\IZ\Gamma}(X_j)}{[\Gamma:\Lambda_i]}=\kappa T^{-1}\rk_{\IZ\Gamma}(X_j).
		\]
		Taking~$T\to \infty$, this shows that~$\widehat{b}_j(X,\Lambda_*;\IF)=0$.
		Hence~$X$ lies in~$\sfH_n^\Gamma(\IF)$.
	
		(ii) Suppose that~$X$ lies in~$\sfCWR_n^\Gamma$ witnessed by the constant~$\kappa$.
		Let~$\Lambda_*$ be a residual chain in~$\Gamma$ and let~$j\le n$.
		For all~$T\in \IR_{\ge 1}$, there exists~$i_0\in \IN$ such that for all~$i\ge i_0$, the $\IZ$-chain complex~$X_{\Lambda_i}$ admits a weak $n$-rebuilding~$(X_{\Lambda_i},X_{\Lambda_i}')$ of quality~$(T,\kappa)$.
		Using Gabber's estimate~\eqref{eqn:Gabber}, for all~$j\le n-1$, we have
		\begin{align*}
			\log\tors H_j(X'_{\Lambda_i})
			&\le \rk_\IZ\bigl((X'_{\Lambda_i})_j\bigr)\log_+\|\partial^{X'_{\Lambda_i}}_{j+1}\|
			\\
			&\le \kappa T^{-1} \rk_\IZ\bigl((X_{\Lambda_i})_j\bigr)\kappa (1+\log T)
			\\
			&= \kappa^2 T^{-1} [\Gamma:\Lambda_i]\rk_{\IZ\Gamma}(X_j)(1+\log T).
		\end{align*}
		Since~$H_j(X_{\Lambda_i})$ is a retract of~$H_j(X'_{\Lambda_i})$, we have
		\[
			\log\tors H_j(X_{\Lambda_i})\le \log\tors H_j(X'_{\Lambda_i}).
		\]
		Together, for all~$j\le n-1$, we have
		\begin{align*}
			\widehat{t}_j(X,\Lambda_*)
			&\le \limsup_{i\to \infty} \frac{\kappa^2 T^{-1} [\Gamma:\Lambda_i]\rk_{\IZ\Gamma}(X_j)(1+\log T)}{[\Gamma:\Lambda_i]}
			\\
			&= \kappa^2 T^{-1}\rk_{\IZ\Gamma}(X_j)(1+\log T).
		\end{align*}
		Taking $T\to \infty$, this shows that~$\widehat{t}_j(X,\Lambda_*)=0$.
		Hence~$X$ lies in~$\sfT_{n-1}^\Gamma$.
		
		(iii) Suppose that~$X$ lies in~$\sfCWR_n^\Gamma$ witnessed by the constant~$\kappa$ and that~$X$ lies in~$\sfFG_{n+1}^\Gamma(\IZ)$.
		Let~$\Lambda_*$ be a residual chain in~$\Gamma$.
		In view of part~(ii), it remains to show that $\widehat{t}_n(X,\Lambda_*)=0$.
		
		For all~$T\in \IR_{\ge 1}$, there exists~$i_0\in \IN$ such that for all~$i\ge i_0$, the $\IZ$-chain complex~$X_{\Lambda_i}$ admits a weak $n$-rebuilding~$(X_{\Lambda_i},X_{\Lambda_i}',\xi,\xi',\Xi)$ of quality~$(T,\kappa)$.
		In the following, let~$i\ge i_0$.
		We have $\|\partial^{X_{\Lambda_i}}_{n+1}\|\le \|\partial^X_{n+1}\|$ by Lemma~\ref{lem:coinvariants}~\ref{item:coinvariants map} and~$\|\partial^X_{n+1}\|$ is finite because~$X_{n+1}$ is finitely generated over~$\IZ\Gamma$.
		Define the $\IZ$-chain complex~$X_{\Lambda_i}''$ with chain modules
		\[
			(X_{\Lambda_i}'')_j\coloneqq \begin{cases}
				0 & \text{if }j\ge n+2; \\
				(X_{\Lambda_i})_{n+1} & \text{if }j=n+1; \\
				(X_{\Lambda_i}')_j & \text{if }j\le n;
			\end{cases}
		\]
		and differentials
		\[
			\partial^{X_{\Lambda_i}''}_{j}\coloneqq \begin{cases}
				0 & \text{if }j\ge n+2; \\
				\xi_n\circ \partial^{X_{\Lambda_i}}_{n+1} & \text{if }j=n+1; \\
				\partial^{X_{\Lambda_i}'}_j  & \text{if }j\le n.
			\end{cases}
		\]
		Consider the partial chain maps $\overline{\xi}\colon X_{\Lambda_i}\to X_{\Lambda_i}''$ and $\overline{\xi'}\colon X_{\Lambda_i}''\to X_{\Lambda_i}$ defined up to degree~$n+1$ as follows:
		\begin{align*}
			\overline{\xi}_j &\coloneqq \begin{cases}
				\id_{(X_{\Lambda_i})_{n+1}} & \text{if }j=n+1; \\
				\xi_j & \text{if } j\le n;
			\end{cases}
			\\
			\overline{\xi'}_j &\coloneqq \begin{cases}
				\id_{(X_{\Lambda_i})_{n+1}}-\Xi_n\circ \partial^{X_{\Lambda_i}}_{n+1}&\text{if }j=n+1; \\
				\xi'_j &\text{if }j\le n.
			\end{cases}
		\end{align*}
		We have constructed
		\[\begin{tikzcd}
			(X_{\Lambda_i})_{n+1}\ar{r}{\partial^{X_{\Lambda_i}}_{n+1}}
			\ar[bend left]{dd}{\id}
			& (X_{\Lambda_i})_n\ar{r}{\partial^{X_{\Lambda_i}}_n}
			\ar[bend left]{dd}{\xi_n}
			& (X_{\Lambda_i})_{n-1}\ar{r}
			\ar[bend left]{dd}{\xi_{n-1}}
			& \cdots
			\\
			\\
			(X_{\Lambda_i})_{n+1}\ar{r}[swap]{\xi_n\circ \partial^{X_{\Lambda_i}}_{n+1}}
			\ar[bend left]{uu}{\id-\Xi_n\circ \partial^{X_{\Lambda_i}}_{n+1}}
			& (X_{\Lambda_i}')_n\ar{r}[swap]{\partial^{X_{\Lambda_i}'}_n}
			\ar[bend left]{uu}{\xi'_n}
			& (X_{\Lambda_i}')_{n-1}\ar{r}
			\ar[bend left]{uu}{\xi'_{n-1}}
			& \cdots
		\end{tikzcd}\]
		The $\IZ$-chain homotopy~$\Xi$ provides a partial chain homotopy between~$\id_{X_{\Lambda_i}}$ and~$\overline{\xi'}\circ \overline{\xi}$ up to degree~$n$.
		In particular, $H_n(X_{\Lambda_i})$ is a retract of~$H_n(X_{\Lambda_i}'')$ and hence
		\[	
			\log\tors H_n(X_{\Lambda_i})\le \log\tors H_n(X_{\Lambda_i}'').
		\]
		By construction, we have
		\begin{align*}
			\rk_\IZ\bigl((X_{\Lambda_i}'')_n\bigr)
			&= \rk_\IZ\bigl((X_{\Lambda_i}')_n\bigr)
			\le \kappa T^{-1} \rk_\IZ\bigl((X_{\Lambda_i})_n\bigr)
			\le \kappa T^{-1} [\Gamma:\Lambda_i]\rk_{\IZ\Gamma}(X_n)
		\end{align*}
		and
		\begin{align*}
			\log_+\|\partial^{X_{\Lambda_i}''}_{n+1}\| 
			&= \log_+\|\xi_n\circ \partial^{X_{\Lambda_i}}_{n+1}\|
			\le \log_+\|\xi_n\|+\log_+\|\partial^{X_{\Lambda_i}}_{n+1}\|
			\\
			&\le \kappa(1+\log T)+\log_+\|\partial^{X}_{n+1}\|.
		\end{align*}
		Together, using Gabber's estimate~\eqref{eqn:Gabber} we obtain
		\[
			\log\tors H_n(X_{\Lambda_i}'')\le \kappa T^{-1}[\Gamma:\Lambda_i]\rk_{\IZ\Gamma}(X_n)\bigl(\kappa(1+\log T)+\log_+\|\partial^X_{n+1}\|\bigr).
		\]
		Finally,
		\begin{align*}
			\widehat{t}_n(X,\Lambda_*)
			&\le \limsup_{i\to \infty}\frac{\log\tors H_n(X_{\Lambda_i}'')}{[\Gamma:\Lambda_i]}
			\\
			&\le \kappa T^{-1}\rk_{\IZ\Gamma}(X_n)\bigl(\kappa(1+\log T)+\log_+\|\partial^X_{n+1}\|\bigr).
		\end{align*}
		Since $\|\partial_{n+1}^X\| < \infty$, taking~$T\to \infty$, this shows that~$\widehat{t}_n(X,\Lambda_*)=0$.
		Hence~$X$ lies in~$\sfT_n^\Gamma$.
		
		(iv) Since~$\Gamma$ is of type~$\sfFP_{n+1}$, there exists a based free $\IZ\Gamma$-resolution~$Y$ of~$\IZ$ lying in~$\sfFG_{n+1}^\Gamma(\IZ)$.
		We will show that~$Y$ also lies in~$\sfCWR_n^\Gamma$.
		Then the claim follows from part~(iii).
		
		Since~$\Gamma$ lies in~$\sfCWR_n$, there exists a based free $\IZ\Gamma$-resolution~$X$ of~$\IZ$ lying in~$\sfCWR_n^\Gamma$ witnessed by the constant~$\kappa$.
		We fix a $\IZ\Gamma$-chain homotopy equivalence between~$Y$ and~$X$ given by the Fundamental Lemma of homological algebra.
		Since~$Y$ and~$X$ lie in~$\sfFG_n^\Gamma(\IZ)$, there exists~$\mu\in \IR_{\ge 1}$ such that for every finite index subgroup~$\Lambda$ of~$\Gamma$, the $\IZ\Gamma$-chain homotopy equivalence between~$Y$ and~$X$ descends to a weak $n$-rebuilding~$(Y_\Lambda,X_\Lambda)$ of quality~$(1,\mu)$ (Remark~\ref{rem:equivariant rebuilding}).
		
		Let~$T\in \IR_{\ge 1}$ and let~$\Lambda_*$ be a residual chain in~$\Gamma$.
		Since~$X$ lies in~$\sfCWR_n^\Gamma$, there exists~$i_0\in \IN$ such that for all~$i\ge i_0$, there exists a weak $n$-rebuilding~$(X_{\Lambda_i},X_{\Lambda_i}')$ of quality~$(T,\kappa)$.
		The composition of the weak $n$-rebuildings~$(Y_{\Lambda_i},X_{\Lambda_i})$ and~$(X_{\Lambda_i},X_{\Lambda_i}')$ is a weak $n$-rebuilding~$(Y_{\Lambda_i},X_{\Lambda_i}')$ of quality~$(T,2\kappa\mu)$ by Lemma~\ref{lem:composition}~\ref{item:composition weak rebuilding}.
		Hence~$Y$ lies in~$\sfCWR_n^\Gamma$ witnessed by the constant~$2\kappa\mu$.		
	\end{proof}
\end{lem}

We summarise the relations between the various classes above in a diagram.
	Let~$n\in \IZ$ and let~$\IF$ be a field.
	The following inclusions hold between classes of groups: 
	\[
	\begin{tikzcd}
		\sfCR_n\ar[r, hook] 
		& \sfCWR_n\ar[r, hook]\ar[d, hook]
		\ar[ld,hook, above, sloped, "\sfFP_{n+1}"]
		& \sfCD_n\ar[d, hook]
		&& \text{cheap rebuilding}
		\\
		\sfT_n
		& \sfT_{n-1}
		& \sfH_n(\IF)\ar[r, hook, "\sfFP_{n}"]
		& \sfH_n(\IQ)\ar[d, hook, "\sfFP_{n+1}"]
		& \text{homology growth}
		\\
		\sfI_n\ar[rrr, hook]
		&&& \sfA_n\ar[d, hook, bend right, "\sfFP_{n+1}" left]
		& \ell^2\text{-invariants}
		\\
		\sfC(\le 0^-)_n \ar[u,equal]
		\ar[r, hook]
		& \sfC(\le 0)_n \ar[r, hook]
		& \sfC(< \infty)_n \ar[r, hook]
		& \sfCM_n \ar[u,hook]
	\end{tikzcd}
	\]
	Here some inclusions hold when restricting to groups of type~$\sfFP_{n+1}$ or~$\sfFP_n$, respectively, as indicated by the labels. 
	
Finally, we show that the algebraic cheap rebuilding property is an equivariantly bootstrappable property (Definition~\ref{defn:equivariant bootstrappable}).
The main work has already been done in Proposition~\ref{prop:rebuilding cones}, where we showed that rebuildings are preserved under mapping cones.
	
\begin{prop}
\label{prop:bootstrappable CR}
	Let the question mark symbol~$\qm$ be the is-class of residually finite groups.
	The following hold:
	\begin{enumerate}[label=\enum]
		\item The family~$\sfCR^\qm_*$ is an equivariantly bootstrappable property of chain complexes over~$\IZ$;
		\item\label{item:bootstrappable CWR} Let~$\Gamma$ be a residually finite group.
		The following hold:
		\begin{enumerate}[label={\rm{(\alph*)}}]
			\item The sequence~$\sfCWR_*^\Gamma$ of classes of $\IZ\Gamma$-chain complexes satisfies axioms~\ref{ax:pos} and~\ref{ax:susp};
			\item\label{item:CWR cone} 
			Let~$f\colon X\to Y$ be a $\IZ\Gamma$-chain map.
			If~$X\in \sfCR_{n-1}^\Gamma$ and~$Y\in \sfCWR_n^\Gamma$, then~$\Cone(f)\in \sfCWR_n^\Gamma$;
			\item The family~$\sfCWR^\qm_*$ satisfies axiom~\ref{ax:ind};
		\end{enumerate}
		\item The family~$\sfCD^\qm_*$ is an equivariantly bootstrappable property of chain complexes over~$\IZ$.
	\end{enumerate}
	\begin{proof}
	(i) First, we show that for every residually finite group~$\Gamma$, the sequence~$\sfCR_*^\Gamma$ is a bootstrappable property of $\IZ\Gamma$-chain complexes.
	Indeed, the axioms~\ref{ax:pos} and~\ref{ax:susp} clearly hold.
	For axiom~\ref{ax:cone}, let $f\colon X\to Y$ be a $\IZ\Gamma$-chain map, where~$X$ lies in~$\sfCR_{n-1}^\Gamma$ witnessed by the constant~$\kappa_X$ and~$Y$ lies in~$\sfCR_n^\Gamma$ witnessed by the constant~$\kappa_Y$.
	Let the constant~$\kappa\in \IR_{\ge 1}$ be as in Proposition~\ref{prop:rebuilding cones}, let~$T\in \IR_{\ge 1}$, and let~$\Lambda_*$ be a residual chain in~$\Gamma$.
	There exist~$i_0^X,i_0^Y\in \IN$ such that for all~$i\ge i_0\coloneqq \max\{i_0^X,i_0^Y\}$, the $\IZ$-chain complex~$X_{\Lambda_i}$ admits an $(n-1)$-rebuilding of quality~$(T,\kappa_X)$ and~$Y_{\Lambda_i}$ admits an $n$-rebuilding of quality~$(T,\kappa_Y)$.
	Consider the $\IZ$-chain map $f_{\Lambda_i}\colon X_{\Lambda_i}\to Y_{\Lambda_i}$.
	By Proposition~\ref{prop:rebuilding cones}, the $\IZ$-chain complex~$\Cone(f_{\Lambda_i})$ admits an $n$-rebuilding of quality~$(T,\kappa)$.
	Since~$\Cone(f_{\Lambda_i})\cong \Cone(f)_{\Lambda_i}$, this shows that the $\IZ\Gamma$-chain complex~$\Cone(f)$ lies in~$\sfCR_n^\Gamma$ witnessed by the constant~$\kappa$.
	
	Second, we show that the family~$\sfCR_*^\qm$ satisfies axiom~\ref{ax:ind}.
	Let~$\Delta$ be a subgroup of~$\Gamma$ and let~$X$ be a $\IZ\Delta$-chain complex lying in~$\sfCR_n^\Delta$ witnessed by the constant~$\kappa$.
	Let~$T\in \IR_{\ge 1}$ and let~$\Lambda_*$ be a residual chain in~$\Gamma$.
	There exists~$i_0\in \IN$ such that for all~$i\ge i_0$, the $\IZ$-chain complex~$X_{\Lambda_i\cap \Delta}$ admits an $n$-rebuilding of quality~$(T,\kappa)$.
	Lemma~\ref{lem:coinvariants}~\ref{item:coinvariants chain} provides an isomorphism of $\IZ$-chain complexes
	\[
		(\ind_\Delta^\Gamma X)_{\Lambda_i}\cong \bigoplus_{\frac{[\Gamma:\Lambda_i]}{[\Delta:\Lambda_i\cap \Delta]}} X_{\Lambda_i\cap \Delta}.
	\]
	By Lemma~\ref{lem:rebuilding sums}, the $\IZ$-chain complex~$(\ind_\Delta^\Gamma X)_{\Lambda_i}$ admits an $n$-rebuilding of quality~$(T,\kappa)$.
	This shows that the $\IZ\Gamma$-chain complex~$\ind_\Delta^\Gamma X$ lies in~$\sfCR_n^\Gamma$ witnessed by the constant~$\kappa$.
	
	Parts~(ii) and~(iii) follow from the same proof as part~(i) by replacing ``rebuilding'' with ``weak rebuilding'' or ``domination'' in the appropriate places.
	\end{proof}
\end{prop}

We obtain the Bootstrapping Theorem~\ref{thm:bootstrapping} for~$\sfCR_*$, which is an algebraic version of~\cite[Theorem~10.9]{ABFG21}, and for~$\sfCD_*$.
Since the group~$\IZ$ lies in~$\sfCR_\infty$ (Example~\ref{ex:integers}), it follows that residually finite infinite elementary amenable group of type~$\sfFP_\infty$ lie in~$\sfCR_\infty$ (Example~\ref{ex:elementary:ame}).

\begin{rem}[Arithmetic groups]\label{rem:arithmetic groups}
	Let $G$ be a connected semisimple algebraic $\Q$-group. Let $\Gamma\subset G(\Q)$ be an arithmetic group. 
		Let $r$ be the rational rank of~$G$. It is a classical result by Serre that $\Gamma$ is of type~$\sfFP_\infty$. 
	Since the group~$\IZ$ lies in~$\sfCR_\infty$ (Example~\ref{ex:integers}) and because of Example~\ref{ex:arithmetic groups}, the group~$\Gamma$ lies in $\sfCR_{r-1}$. In particular, $\Gamma\in \sfT_{r-1}$ by Lemma~\ref{lem:CWR implies T}~\ref{item:additional degree}, that is, $\Gamma$ has vanishing torsion homology growth up to degree~$r-1$. Similarly, we obtain that $\Gamma\in\sfH_{r-1}(\IF)$ for every field~$\IF$. This reproves Theorem~\ref{thm:SLd} for~$\Gamma=\mathrm{SL}_d(\IZ)$ and~$r=d-1$. 
\end{rem}

Proposition~\ref{prop:bootstrappable CR}~\ref{item:bootstrappable CWR} means that the family~$\sfCWR^\qm_*$ is close to being an equivariantly bootstrappable property of chain complexes over~$\IZ$.
However, we can show axiom~\ref{ax:cone} only for chain maps under stronger assumptions on the domain, which goes back to Proposition~\ref{prop:rebuilding cones}~\ref{item:warning}.
Nevertheless, we obtain a modified Bootstrapping Theorem for~$\sfCWR_*$ under stronger assumptions on stabilisers of cells of dimension~$\ge 1$. 

\begin{thm}
\label{thm:bootstrapping CWR}
	Let~$\Gamma$ be a residually finite group and let~$n\in \IN$.
	Let~$\Omega$ be a $\Gamma$-CW-complex such that the following hold:
	\begin{enumerate}[label=\enum]
		\item $\Omega$ is $(n-1)$-acyclic;
		\item $\Gamma\backslash \Omega^{(n)}$ is compact;
		\item For every vertex~$v$ of~$\Omega$, the stabiliser~$\Gamma_v$ lies in~$\sfCWR_n$;
		\item For every cell~$\sigma$ of~$\Omega$ with $\dim(\sigma)\in \{1,\ldots,n\}$, the stabiliser~$\Gamma_\sigma$ lies in~$\sfCR_{n-\dim(\sigma)}$.
	\end{enumerate}
	Then $\Gamma\in \sfCWR_n$.
	\begin{proof}
		The proof is analogous to the proof of Theorem~\ref{thm:bootstrapping}.
		The main difference is to replace the use of axiom~\ref{ax:cone} by Proposition~\ref{prop:bootstrappable CR}~\ref{item:bootstrappable CWR}~\ref{item:CWR cone}.
		For instance, the proof of Proposition~\ref{prop:bootstrapping} shows also the following:
		Let~$X$ be an $R\Gamma$-chain complex that is concentrated in degrees~$\ge 0$.
		Suppose that the $R\Gamma$-module~$X_0$ admits a projective resolution lying in~$\sfCWR^\Gamma_n$ and that for all~$j\in \{1,\ldots,n\}$ the $R\Gamma$-module~$X_j$ admits a projective resolution lying in~$\sfCR^\Gamma_{n-j}$. 
		Then there exists a projective $R\Gamma$-chain complex~$\overline{X}\in \sfCWR^\Gamma_n$ that is weakly equivalent to~$X$.
	\end{proof}
\end{thm}

\begin{rem}[Geometric cheap rebuilding property]
\label{rem:CR}
	The algebraic cheap rebuilding property (Definition~\ref{defn:CR}) is an algebraic version of \ABFG's geometric cheap rebuilding property~\cite[Definition~10.6]{ABFG21}.
	However, there is no obvious implication between the algebraic and geometric cheap rebuilding properties due to differences in the definitions of algebraic and geometric rebuildings (Remark~\ref{rem:rebuilding}).
	
	Another significant difference is that in Definition~\ref{defn:CR} we consider only residual chains for simplicity, instead of more general Farber sequences and Farber neighbourhoods.
	The advantage is that axiom~\ref{ax:ind} for~$\sfCR^\qm_*$ is easy to show via Lemma~\ref{lem:coinvariants}~\ref{item:coinvariants chain}.
	However, the family~$\sfCR^\qm_*$ does not seem to satisfy axiom~\ref{ax:res}, due to the non-transitivity of normality.
	In particular, we are not able to show that the class $\sfCR_n$ of groups is closed under commensurability.
	On the other hand, the geometric cheap rebuilding property of~\ABFG\ is invariant under commensurability~\cite[Corollary~10.13~(1)]{ABFG21}.
	It seems plausible that a strengthened version of~$\sfCR^\qm_*$ involving Farber neighbourhoods analogous to~\cite[Definition~10.6]{ABFG21} also satisfies axiom~\ref{ax:res}.
	
	It is not known whether residually finite infinite amenable groups of type~$\sfF_\infty$ satisfy the geometric cheap rebuilding property~\cite[Question~10.21]{ABFG21}.
	However, it follows from~\cite[Section~10]{ABFG21} using the argument in Example~\ref{ex:elementary:ame} that residually finite infinite \emph{elementary} amenable groups of type~$\sfFP_\infty$ satisfy the geometric cheap rebuilding property.
	Similarly, using the argument in Corollary~\ref{cor:commensurability}~\ref{item:boot commensurated}, it follows that the geometric cheap rebuilding property passes from commensurated subgroups to overgroups.

        In contrast to the geometric approach, the algebraic approach also allows
        us to formulate conditions such as $\sfCD_n^\Gamma(\IF)$ over general fields~$\IF$.
        These could lead to further examples of groups in~$\sfH_n(\IF)$. Currently, we are
        not aware of concrete new examples of this type.
\end{rem}	

\begin{rem}	
	The algebraic cheap domination property (Definition~\ref{defn:CR}) is related to Bridson--Kochloukova's geometric notion of ``slowness''~\cite{Bridson-Kochloukova17}\cite[Section~5.2]{Lueck16_survey}.
	Slow groups have vanishing homology growth.
	It is not known whether residually finite amenable groups are slow~\cite[Section~5]{Bridson-Kochloukova17}.
\end{rem}


\subsection{Amenable groups}
\label{sec:amenable}

We show that amenable groups satisfy the algebraic cheap \emph{weak} rebuilding property (Definition~\ref{defn:CR}).
The proof is inspired by the work of Kar--Kropholler--Nikolov~\cite{KKN17}, using F\o lner sequences that are compatible with residual chains.  
We work over the ring~$R=\IZ$.

\begin{lem}
\label{lem:quotient}
	Let $0\to A\xrightarrow{f} X\xrightarrow{g} Y\to 0$ be a short exact sequence of free $\IZ$-chain complexes.
	If the inclusion~$f$ is nullhomotopic, then there exists a $\IZ$-chain map~$h\colon Y\to X$ such that~$\id_X\simeq h\circ g$. 
	\begin{proof}
		The obvious $\IZ$-chain map $q\colon \Cone(f)\to Y$ makes the following diagram commute	
		\[\begin{tikzcd}
			X\ar[hook]{r}\ar{dr}[swap]{g} & \Cone(f)\ar{d}{q} \\
			&Y
		\end{tikzcd}\]
		The map~$q$ is a weak equivalence by the long exact sequences in homology and the Five Lemma.
		Since~$\Cone(f)$ and~$Y$ are free $\IZ$-chain complexes, the weak equivalence~$q$ is a chain homotopy equivalence.
		Let~$r\colon Y\to \Cone(f)$ be a chain homotopy inverse of~$q$.
		Then the composition~$r\circ g$ is chain homotopic to the inclusion~$X\into \Cone(f)$.
		Set~$h\colon Y\to X$ to be the composition
		\[
			h\colon Y\xrightarrow{r} \Cone(f)\xrightarrow{(\id_A,\id_X;H)} \Cone(0)\cong \Sigma A\oplus X\xrightarrow{p_2} X,
		\]
		where~$H$ is a chain homotopy between the nullhomotopic inclusion~$f$ and the zero map~$0$, and~$p_2$ is the projection to the second summand.
		By construction, the composition~$h\circ g$ is homotopic to~$\id_X$.
	\end{proof}
\end{lem}

A \emph{based short exact sequence} of based free $\IZ$-modules is a short exact sequence $0\to K\xrightarrow{f} M\xrightarrow{g} L\to 0$ of based free $\IZ$-modules with bases~$I_K$, $I_M$, and~$I_L$, respectively, such that $f(I_K)\subset I_M$ and $g(I_M\smallsetminus f(I_K))=I_L$.
We say that~$f(K)$ is a \emph{based submodule} of~$M$.
For based free $\IZ$-chain complexes, we obtain corresponding degreewise notions of based subcomplexes and based short exact sequences.

Let~$X$ be a based free $\IZ$-chain complex that is concentrated in degrees~$\ge 0$.
Consider the \emph{standard augmentation map}~$\varepsilon\colon X_0\to \IZ$ that maps each $\IZ$-basis element of~$X_0$ to~$1\in \IZ$.
We say that~$X$ is \emph{augmented over~$\IZ$} if~$\varepsilon\circ \partial^X_1=0$.
If~$X$ is augmented over~$\IZ$, we define the \emph{augmented chain complex~$X^\varepsilon$ associated to~$X$} as the based free $\IZ$-chain complex with chain modules
\[
	X^\varepsilon_j\coloneqq \begin{cases}
		X_j & \text{if }j\ge 0; \\
		\IZ & \text{if }j=-1; \\
	\end{cases}
\]
and differentials
\[
	\partial^{X^\varepsilon}_j\coloneqq \begin{cases}
		\partial^X_j & \text{if }j\ge 1; \\
		\varepsilon & \text{if }j=0.
	\end{cases}
\]
If~$X$ is augmented over~$\IZ$, then so is every based subcomplex of~$X$.
We say that the inclusion~$f\colon A\to X$ of a based subcomplex is \emph{augmentedly nullhomotopic} if the $\IZ$-chain map $f^\varepsilon\colon A^\varepsilon\to X^\varepsilon$ given by
\[
	f^\varepsilon_j\coloneqq \begin{cases}
		f_j & \text{if }j\ge 0; \\
		\id_\IZ & \text{if }j=-1; 
	\end{cases}
\]
is nullhomotopic.
In this case we also say that~$A$ is \emph{augmentedly contractible in~$X$}. 

For a $\IZ$-chain complex~$Y$ that is concentrated in degrees~$\ge 0$, we denote by~$Y^+$ the $\IZ$-chain complex with chain modules
\[
	Y^+_j\coloneqq \begin{cases}
		Y_j & \text{if }j\ge 1; \\
		Y_0\oplus \IZ & \text{if }j=0;
	\end{cases}
\]
and differentials
\[
	\partial^{Y^+}_j\coloneqq \begin{cases}
		\partial^Y_j & \text{if }j\ge 2; \\
		(\partial^Y_1,0) & \text{if }j=1.
	\end{cases}
\]

\begin{lem}
	\label{lem:quotient-augmented}
	Let~$0\to A\xrightarrow{f} X\xrightarrow{g}Y\to 0$ be a based short exact sequence of based free $\IZ$-chain complexes that are concentrated in degrees~$\ge 0$.
	If~$X$ is augmented over~$\IZ$ and the inclusion~$f$ is augmentedly nullhomotopic, then there exist $\IZ$-chain maps $g^+\colon X\to Y^+$ and $h^+\colon Y^+\to X$ such that~$\id_X\simeq h^+\circ g^+$. 
	
	Moreover, let~$n\in \IN$ and suppose that for all~$j\le n$ the $\IZ$-module~$X_j$ is finitely generated.
	Denote 
	\[
		T'\coloneqq \min\biggl\{\frac{\rk_\IZ(X_j)}{\rk_\IZ(Y^+_j)}\ \biggm| 0\le j\le n\biggr\},
		\quad
		\kappa\coloneqq \max\Bigl\{1, \max\bigl\{\log_+\|\partial^X_j\|\bigm\vert 0\le j\le n\bigr\}\Bigr\}.
	\]
	Then, for all~$T\in \IR_{\ge 1}$ with~$T\le T'$, there exists a weak $n$-rebuilding~$(X,Y^+)$ of quality~$(T,\kappa)$.
	\begin{proof}
		We have a short exact sequence of free $\IZ$-chain complexes
		\[
			0\to A^\varepsilon\xrightarrow{f^\varepsilon} X^\varepsilon\xrightarrow{g} Y\to 0,
		\]
		where, by abuse of notation, $g\colon X^\varepsilon\to Y$ is the $\IZ$-chain map given by~$g\colon X\to Y$ in degrees~$\ge 0$ and by zero in degree~$-1$.
		Lemma~\ref{lem:quotient} yields a $\IZ$-chain map~$h\colon Y\to X^\varepsilon$ such that there exists a chain homotopy~$H\colon \id_{X^\varepsilon}\simeq h\circ g$.
		We define the chain maps~$g^+\colon X\to Y^+$ and~$h^+\colon Y^+\to X$ by
		\begin{align*}
			g^+_j & \coloneqq \begin{cases}
				g_j & \text{if } j\ge 1; \\
				(g_0,\varepsilon) & \text{if } j=0;
			\end{cases}
			\\
			h^+_j & \coloneqq \begin{cases}
				h_j & \text{if } j\ge 1; \\
				h_0\oplus H_{-1} & \text{if } j=0.
			\end{cases}
		\end{align*}
		We have constructed
		\[\begin{tikzcd}
			\cdots\ar{d} && \cdots\ar{d} && \cdots\ar{d}
			\\
			X_2\ar{rr}{g_2}\ar{d}[swap]{\partial^X_2}
			&& Y_2\ar{rr}{h_2}\ar{d}{\partial^Y_2}
			&& X_2\ar{d}{\partial^X_2}
			\\
			X_1\ar{rr}{g_1}\ar{d}[swap]{\partial^X_1} && Y_1\ar{rr}{h_1}\ar{d}{(\partial^Y_1,0)} && X_1\ar{d}{\partial^X_1}
			\\
			X_0\ar{rr}{(g_0,\varepsilon)}
			&& Y_0\oplus \IZ\ar{rr}{h_0\oplus H_{-1}}
			&& X_0
		\end{tikzcd}\]
		Then one easily checks that~$H_{*\ge 0}$ provides a chain homotopy $\id_X\simeq h^+\circ g^+$.
				
		Now, suppose that~$X_j$ is finitely generated for all~$j\le n$.
		For all~$j\le n$, we have
		\begin{align*}
			\rk_\IZ(Y^+_j) &\le T'^{-1} \rk_\IZ(X_j) \le T^{-1}\rk_\IZ(X_j);
			\\
			\|\partial^{Y^+}_j\| &\le \|\partial^X_j\| \le \exp(\kappa);
			\\
			\|g^+_j\| &\le 2.
		\end{align*}
		Hence~$(X,Y^+)$ is a weak $n$-rebuilding of quality~$(T,\kappa)$.
	\end{proof}
\end{lem}

In the following, when we write~$\lim_{i\to \infty}\frac{a_i}{b_i}$ for~$a_i,b_i\in \IZ$, we assume implicitly that~$b_i\neq 0$ for all large enough~$i\in \IN$.

\begin{prop}
\label{prop:boundary}
	Let~$\Gamma$ be a residually finite infinite group and let~$n\in \IN$.
	Let~$X$ be a based free $\IZ\Gamma$-resolution of~$\IZ$ (with the standard augmentation) such that the $\IZ\Gamma$-module~$X_j$ is finitely generated for all~$j\le n$.
	Suppose that for all residual chains~$\Lambda_*$ in~$\Gamma$ and all large enough~$i\in \IN$, the based free $\IZ$-chain complex~$X^i\coloneqq X_{\Lambda_i} =\IZ\otimes_{\IZ[\Lambda_i]}\res^\Gamma_{\Lambda_i}X$ admits a based subcomplex~$A^i$
	that is augmentedly contractible in~$X^i$ such that for all~$j\in \{0,\ldots,n\}$, we have
	\begin{equation}
	\label{eqn:boundary}
		\lim_{i\to \infty}\frac{\rk_\IZ(X^i_j)-\rk_\IZ(A^i_j)}{\rk_\IZ(X^i_j)}=0.
	\end{equation}
	Then~$\Gamma\in \sfCWR_n$.
	\begin{proof}
		Set~$\kappa\coloneqq \max\{1,\max\{\log_+\|\partial^X_j\|\mid j\le n\}\}$.
		Let~$T\in \IR_{\ge 1}$ and let~$\Lambda_*$ be a residual chain in~$\Gamma$.
		For~$i$ large enough, we denote by~$Y^i\coloneqq X^i/A^i$ the quotient chain complex. 
		For all~$j\in \{0,\ldots,n\}$, we have
		\[
			\lim_{i\to \infty}\frac{\rk_\IZ(Y^i_j)}{\rk_\IZ(X^i_j)}=0
		\]
		by assumption.
		Moreover, in degree~$0$ we have 
		\[
			\lim_{i\to \infty}\frac{\rk_\IZ(Y^i_0)+1}{\rk_\IZ(X^i_0)}=0
		\]
		using that~$\lim_{i\to \infty}\rk_\IZ(X^i_0)=\infty$ because~$\Gamma$ is infinite. 
		Together, we have
		\[
			\lim_{i\to \infty}\frac{\rk_\IZ(X^i_j)}{\rk_\IZ\bigl((Y^i)^+_j\bigr)}=\infty.
		\]
		Hence there exists~$i_0\in \IN$ such that for all~$i\ge i_0$
		\[
			\min\biggl\{\frac{\rk_\IZ(X^i_j)}{\rk_\IZ\bigl((Y^i)^+_j\bigr)}\ \biggm| 0\le j\le n\biggr\}\ge T.
		\]
		By Lemma~\ref{lem:coinvariants}~\ref{item:coinvariants map}, we have
		\[
			\|\partial^{X^i}_j\|\le \|\partial^X_j\|\le \exp(\kappa).
		\]
		The based free $\IZ\Gamma$-resolution~$X$ of~$\IZ$ is in particular augmented over~$\IZ$ 
		and thus, so is the based free $\IZ$-chain complex~$X^i$.
		Hence Lemma~\ref{lem:quotient-augmented} applies and yields a weak $n$-rebuilding~$(X^i,(Y^i)^+)$ of quality~$(T,\kappa)$.
	\end{proof}
\end{prop}

Condition~\eqref{eqn:boundary} means that asymptotically the subcomplex~$A^i$ is ``large'' in~$X^i$.
We will apply Proposition~\ref{prop:boundary} in the context of amenable groups to appropriate subcomplexes coming from F\o lner sequences. 

Let~$\Gamma$ be a finitely generated group and let~$S$ be a finite generating set of~$\Gamma$.
For a subset~$F$ of~$\Gamma$, we write $\partial_SF\coloneqq \{\gamma\in F\mid \gamma \cdot s\in \Gamma\smallsetminus F \text{ for some }s\in S\}$.
We also write~$\interior_S(F)\coloneqq F\smallsetminus \partial_S F$.
Recall that a \emph{F\o lner sequence} in~$\Gamma$ with respect to~$S$ is a sequence~$F_*=(F_i)_{i\in \IN}$ of non-empty finite subsets of~$\Gamma$ satisfying
\[
	\lim_{i\to \infty}\frac{\#\partial_S F_i}{\# F_i}=0.
\]
A finitely generated group is \emph{amenable} if and only if it admits a F\o lner sequence.
Residually finite amenable groups admit F\o lner sequences satisfying additional properties:

\begin{thm}[{Weiss~\cite{Weiss01}\cite[Theorem~7]{KKN17}}]
\label{thm:Folner}
	Let~$\Gamma$ be a finitely generated amenable group that is residually finite and infinite.
	 Let~$S$ be a finite generating set of~$\Gamma$ and let~$\Lambda_*$ be a residual chain in~$\Gamma$.
	Then there exists a F\o lner sequence~$F_*$ in~$\Gamma$ with respect to~$S$ such that for all~$i\in \IN$, the set~$F_i$ is a set of right-coset representatives for~$\Lambda_i\backslash \Gamma$.
\end{thm}

\noindent
We say that a F\o lner sequence~$F_*$ as in Theorem~\ref{thm:Folner} is \emph{compatible} with the residual chain~$\Lambda_*$.

\begin{thm}
\label{thm:amenable CWR}
	Let $n \in \N$ and let $\Gamma$ be a residually finite infinite amenable group of type~$\sfFP_n$. 
	Then~$\Gamma\in \sfCWR_n$.
	\begin{proof}
		Since~$\Gamma$ is of type~$\sfFP_n$, there exists a based free $\IZ\Gamma$-resolution~$X$ of~$\IZ$ such that the free $\IZ\Gamma$-module~$X_j$ is finitely generated for all~$j\le n$.
		The differentials~$\partial_1,\ldots,\partial_n$ are given by (right-)multiplication with matrices~$M_1,\ldots,M_n$, respectively, whose entries are in~$\IZ\Gamma$.
		Let~$S$ be a finite generating set of~$\Gamma$ containing all group elements appearing in the entries of the matrices~$M_1,\ldots,M_n$.
		For~$j\ge 0$, denote by~$S^j$ the finite set of all words with letters in~$S$ of length~$\le j$.
		Let~$\Lambda_*$ be a residual chain in~$\Gamma$.
		By Theorem~\ref{thm:Folner}, there exists a F\o lner sequence~$F_*$ in~$\Gamma$ with respect to the finite generating set~$S^{n+1}$ that is compatible with the residual chain~$\Lambda_*$.
		For all~$i\in \IN$, we consider the based free $\IZ$-chain complex~$X^i\coloneqq X_{\Lambda_i}= \IZ \otimes_{\IZ[\Lambda_i]} \res^\Gamma_{\Lambda_i} X$.
		For~$i$ large enough, we will construct a based subcomplex~$A^i$ of~$X^i$ that is augmentedly contractible in~$X^i$ such that for all~$j\in \{0,\ldots,n\}$, we have
		\[
			\lim_{i\to \infty}\frac{\rk_\IZ(X^i_j)-\rk_\IZ(A^i_j)}{\rk_\IZ(X^i_j)}=0.
		\]
		Then Proposition~\ref{prop:boundary} yields that~$\Gamma$ lies in~$\sfCWR_n$.	
		
		Now, since the $\IZ\Gamma$-resolution~$X$ is based free, for all~$j$ we have an isomorphism of $\IZ\Gamma$-modules $X_j\cong \bigoplus_{I_j}\IZ\Gamma$, where~$I_j$ is the set of basis elements.
		The differential~$\partial_j\colon X_j\to X_{j-1}$ is given by (right-)multiplication with the $(I_j\times I_{j-1})$-matrix~$M_j=(m_j^{kl})_{k\in I_j,l\in I_{j-1}}$, where~$m_j^{kl}=\sum_{s\in S}\lambda_j^{kl}(s)s\in \IZ\Gamma$ with~$\lambda_j^{kl}(s)\in \IZ$.
		Consider the isomorphism of $\IZ[\Lambda_i]$-modules
		\begin{align*}
			\res^\Gamma_{\Lambda_i} \IZ\Gamma &\cong \bigoplus_{F_i} \IZ[\Lambda_i] \\
			f &\mapsfrom e_f \\
			\gamma &\mapsto \gamma t(\gamma)^{-1}e_{t(\gamma)}
		\end{align*}
		where~$e_f$ is the $\IZ[\Lambda_i]$-basis element corresponding to~$f\in F_i$, and~$t(\gamma)\in F_i$ is such that~$\Lambda_i\gamma=\Lambda_i t(\gamma)$. 
		We can identify the $\IZ$-chain modules of~$X^i$ as
		\[
			X^i_j= \IZ\otimes_{\IZ[\Lambda_i]}\res^\Gamma_{\Lambda_i} X_j \cong \IZ\otimes_{\IZ[\Lambda_i]} \bigoplus_{I_j}\bigoplus_{F_i}\IZ[\Lambda_i] \cong \bigoplus_{I_j\times F_i}\IZ.
		\]
		Under this identification, the differential $\partial^i_j\colon X^i_j\to X^i_{j-1}$ maps the $\IZ$-basis element~$e_{(k,f)}\in X^i_j$ for~$(k,f)\in I_j\times F_i$ to
		\[
			\sum_{l\in I_{j-1}}\sum_{s\in S}\lambda_j^{kl}(s)e_{(l,t(fs))}\in X^i_{j-1}.
		\] 
		Define the based subcomplex~$A^i$ of~$X^i$ by 
		\[
			A^i_j\coloneqq 
			\begin{cases}
				\bigoplus_{I_j\times \interior_{S^{j+1}} (F_i)} \IZ & \text{if }j\in\{0,\ldots,n\}; \\
				0 & \text{if }j\ge n+1.
		\end{cases}\]
		The chain complex~$A^i$ is well-defined, as the differential $\partial^i_j\colon X^i_j\to X^i_{j-1}$ restricts to a differential $A^i_j\to A^i_{j-1}$.
		Indeed, if~$f\in \interior_{S^{j+1}}(F_i)$, then~$fs\in \interior_{S^j}(F_i)\subset F_i$ and hence~$t(fs)=fs\in \interior_{S^j}(F_i)$.
		For~$j\le n$, the index set~$I_j$ is finite by assumption.
		Hence, for all~$j\le n$, we have
		\[
			\lim_{i\to \infty}\frac{\rk_\IZ(X^i_j)-\rk_\IZ(A^i_j)}{\rk_\IZ(X^i_j)}
			= \lim_{i\to \infty}\frac{\#I_j\cdot \#\partial_{S^{j+1}}(F_i)}{\#I_j\cdot \#F_i}
			\le \lim_{i\to \infty}\frac{\#\partial_{S^{n+1}}(F_i)}{\# F_i}
			=0.
		\]
		It remains to show that the inclusion~$A^i\to X^i$ is augmentedly nullhomotopic.
		The inclusion~$A^i\to X^i$ factors through the projection~$X\to X^i$ via the $\IZ$-chain map~$A^i\to X$ that maps the $\IZ$-basis element~$e_{(k,f)}\in A^i_j\cong \bigoplus_{I_j\times \interior_{S^{j+1}(F_i)}}\IZ$ to~$f\cdot e_k\in X_j\cong \bigoplus_{I_j}\IZ\Gamma$.
		In particular, we obtain a factorisation of the augmented inclusion~$(A^i)^\varepsilon\to (X^i)^\varepsilon$ through~$X^\varepsilon$.
		Since~$X$ is a free $\IZ\Gamma$-resolution of~$\IZ$, the augmented $\IZ$-chain complex~$X^\varepsilon$ is contractible.
		Hence the augmented inclusion~$(A^i)^\varepsilon\to (X^i)^\varepsilon$ is nullhomotopic.
	\end{proof}
\end{thm}

By Lemma~\ref{lem:CWR implies T}, we recover the following known result~\cite{Cheeger-Gromov86, LLS11, KKN17}.

\begin{cor}
\label{cor:amenable T}
	Let $n \in \N$ and let $\Gamma$ be a residually finite infinite amenable group of type~$\sfFP_n$. 
	Then~$\Gamma$ lies in~$\sfT_{n-1}$ and~$\sfH_n(\IF)$ for every field~$\IF$.
\end{cor}

As a consequence of the modified Bootstrapping Theorem~\ref{thm:bootstrapping CWR} for~$\sfCWR_*$ we obtain:

\begin{cor}
\label{cor:amenable graph}
	Let~$\Gamma$ be a residually finite fundamental group of a finite graph of groups and let~$n\in \IZ$.
	If all vertex groups are infinite amenable of type~$\sfFP_n$ and all edge groups are infinite elementary amenable of type~$\sfFP_\infty$, then $\Gamma\in \sfCWR_n$ (and hence $\Gamma\in \sfT_{n-1}$).
\end{cor}

The group~$\Gamma$ in Corollary~\ref{cor:amenable graph}, also more generally with (not necessarily elementary) amenable edge groups, lies in~$\sfH_n(\IF)$ for every field~$\IF$.
This already follows from the Bootstrapping Theorem~\ref{thm:bootstrapping} for~$\sfH_*(\IF)$ (Proposition~\ref{prop:bootstrappable H}) and Corollary~\ref{cor:amenable T}.

\begin{rem}[$\ell^2$-Torsion]
	Non-trivial amenable groups of type~$\sfFL$ have vanishing $\ell^2$-torsion~\cite[Theorem~1.3]{Li-Thom}.
	Hence the fundamental group of a finite graph of non-trivial amenable groups that are of type~$\sfFL$ has vanishing $\ell^2$-torsion~\cite[Theorem~3.93]{Lueck02}.
	
	Let~$\Gamma$ be a residually finite fundamental group of a finite graph of groups with non-trivial amenable vertex groups of type~$\sfFL$ and non-trivial elementary amenable edge groups of type~$\sfFL$.
	Then the $\ell^2$-torsion~$\rho^{(2)}(\Gamma)$ vanishes by the above, and~$\Gamma$ lies in~$\sfT_\infty$ by Corollary~\ref{cor:amenable graph}.
	In particular, we have
	\[
		\rho^{(2)}(\Gamma)=0=\sum_{j\ge 0}(-1)^j\cdot \widehat{t}_j(\Gamma,\Lambda_*)
	\]
	for all residual chains~$\Lambda_*$ in~$\Gamma$.
	This confirms a conjecture of L\"uck~\cite[Conjecture~1.11~(3)]{Lueck13} for the group~$\Gamma$.
\end{rem}

\noindent
 \textbf{Acknowledgements.}
 We thank Miko\l aj Fr\k{a}czyk for answering our questions about the article~\cite{ABFG21}
 and Uli Bunke for suggesting to consider Novikov--Shubin invariants.
 We thank the referee for many helpful comments.

 K.L., C.L., and M.U.\ were supported by the CRC~1085 ``Higher Invariants'' (Universit\"at Regensburg, funded by the DFG).

 M.M.\ was supported by the ERC ``Definable Algebraic Topology'' DAT -- Grant Agreement number~101077154. 

 M.M.\ was funded by the European Union -- NextGenerationEU under the National Recovery and Resilience Plan (PNRR) -- Mission 4 Education and research -- Component 2 From research to business -- Investment~1.1 Notice Prin~2022 -- DD~N.~104 del 2/2/2022, from title ``Geometry and topology of manifolds'', proposal code 2022NMPLT8 -- CUP~J53D23003820001.

 R.S.~was funded by the Deutsche Forschungsgemeinschaft (DFG, German Research Foundation) -- project number~338540207.


\bibliographystyle{amsalphaabbrv}
\bibliography{bibcones}

\setlength{\parindent}{0cm} 
\enlargethispage{1cm}

\end{document}